\documentclass[a4paper,9pt]{article}

\usepackage[centertags]{amsmath}

\usepackage{authblk}
\usepackage{multirow}
\usepackage{amssymb}
\usepackage{graphicx}
\usepackage{listing}
\usepackage{mathrsfs}
\usepackage{amsmath,amsthm,amssymb,fontenc,color,textcomp,amsfonts,graphicx,graphics}
\usepackage[colorlinks=true]{hyperref}
\usepackage{hyperref}% http://ctan.org/pkg/hyperref
\hypersetup{colorlinks = true,linkcolor = black,citecolor = blue}
\usepackage{xcolor}
\usepackage{booktabs}

\theoremstyle{definition}
\newtheorem{defn}{Definition}[section]
\newtheorem{thm}{Theorem}[section]

\newtheorem{Ex}{Example}[section]
\newtheorem{prop}{Proposition}[section]

\usepackage[linesnumbered,ruled,vlined]{algorithm2e}

\usepackage{caption}
\usepackage{subcaption}

\newcommand{\be}{\begin{equation}}
\newcommand{\ee}{\end{equation}}

\numberwithin{equation}{section}
\begin{document}
                                            
\title{Tensor Arnoldi-Tikhonov and GMRES-type methods for ill-posed problems with a 
t-product structure\\}

\author{{\Large Lothar Reichel\thanks{\,e-mail: reichel@math.kent.edu} \; and\; Ugochukwu O. Ugwu\thanks{\,e-mail: uugwu@kent.edu}} \vspace{.3cm}
\mbox{\ }\\ 
{\normalsize \large Department of Mathematical Sciences, Kent State University, Kent, 
OH 44240, USA}}

\date{}

%\date{\normalsize{ \large Jan 2020}}

\maketitle \vspace*{-0.5cm}

\thispagestyle{empty}

%\newpage
%==================================================================
\begin{abstract}
This paper describes solution methods for linear discrete ill-posed problems defined by 
third order tensors and the t-product formalism introduced in [M. E. Kilmer and C. D. 
Martin, Factorization strategies for third order tensors, Linear Algebra Appl., 435 
(2011), pp. 641--658]. A t-product Arnoldi (t-Arnoldi) process is defined and applied 
to reduce a large-scale Tikhonov regularization problem for third order tensors to a 
problem of small size. The data may be represented by a laterally oriented matrix or a 
third order tensor, and the regularization operator is a third order tensor. The 
discrepancy principle is used to determine the regularization parameter and the number of
steps of the t-Arnoldi process. Numerical examples compare results for several solution 
methods, and illustrate the potential superiority of solution methods that tensorize over 
solution methods that matricize linear discrete ill-posed problems for third order tensors.  
\vspace{.3cm}
 
  \vspace{-.2cm}
\end{abstract}
%==================================================================

%\label{sec:in}

\section{Introduction}
We are concerned with the solution of large-scale least squares problems of the form
\begin{equation}
\min_{\mathcal{\vec{X}}\in\mathbb{R}^{m\times 1\times n}} 
\|\mathcal{A*\vec{X} - \vec{B}} \|_F,
\label{1}
\end{equation}
where $\mathcal{A} = [a_{ijk}]_{i,j,k=1}^{m,m,n}\in\mathbb{R}^{m\times m\times n}$ is a third order 
tensor of ill-determined tubal rank, i.e., the Frobenius norm of the singular tubes of 
$\mathcal{A}$, which are analogues of the singular values of a matrix, decay rapidly to 
zero with increasing index, and there are many nonvanishing singular tubes of tiny 
Frobenius norm of different orders of magnitude (cf. Definition \ref{ipd} below). Least squares problems with a tensor of 
this kind are referred to as linear discrete ill-posed problems. The tensors 
$\mathcal{\vec{X}}\in\mathbb{R}^{m\times 1\times n}$ and 
$\mathcal{\vec{B}}\in\mathbb{R}^{m\times 1\times n}$ in \eqref{1} are lateral slices of 
third order tensors, and the operator $\ast$ denotes the tensor t-product introduced in 
the seminal work by Kilmer and Martin \cite{KM}. We will review the t-product in Section 
\ref{sec2}. 

An advantage of the formulation \eqref{1} with the t-product, when compared
to other products, is that the t-product avoids loss of information inherent in the 
flattening of a tensor; see Kilmer et al. \cite{KBHH}. The t-product preserves the 
natural ordering and higher correlations embedded in the data, and has been found useful 
in many application areas, including completion of seismic data \cite{EANK}, image 
deblurring problems \cite{GIJS,KBHH,KM,LU}, facial recognition \cite{HKBH}, tomographic 
image reconstruction \cite{SKH}, and tensor compression \cite{ZEAHK}.

Throughout this paper, $\|\cdot\|_F$ denotes the Frobenius norm of a third order tensor,
which for $\mathcal{A} = [a_{ijk}]_{i,j,k=1}^{m,m,n}\in\mathbb{R}^{m\times m\times n}$ is defined by
\begin{equation*}
\|\mathcal{A}\|_F = \sqrt{\sum_{i=1}^m \sum_{j=1}^{m} \sum_{k=1}^n a^2_{ijk}}.
\end{equation*} 

In applications of interest to us, such as image and video restoration, the 
data tensor $\mathcal{\vec{B}} \in \mathbb{R}^{m \times 1 \times n}$ is contaminated by 
measurement error (noise) that is represented by a tensor 
$\mathcal{\vec{E}}\in\mathbb{R}^{m\times 1\times n}$. Thus,
\begin{equation}\label{unava}
\mathcal{\vec{B}} = \mathcal{\vec{B}}_\text{true} + \mathcal{\vec{E}},
\end{equation}
where $\mathcal{\vec{B}}_\text{true} \in \mathbb{R}^{m \times 1 \times n}$ represents the 
unavailable error-free data tensor that is associated with the known data tensor
$\mathcal{\vec{B}}$. We assume the unavailable linear system of equations
\[
\mathcal{A*\vec{X} = \vec{B}}_\text{true}
\]
to be consistent and let $\mathcal{\vec{X}_\text{true}}$ denote its (unknown) exact 
solution of minimal Frobenius norm. 

We would like to compute an accurate approximation of $\mathcal{\vec{X}_\text{true}}$. 
Straightforward solution of \eqref{1} typically does not yield a meaningful approximation 
of $\mathcal{\vec{X}}_\text{true}$, because the severe ill-conditioning of $\mathcal{A}$ 
and the error in $\mathcal{\vec{B}}$ result in a large propagated error in the computed
solution. We remedy this difficulty by replacing \eqref{1} by a 
nearby problem that is less sensitive to perturbations of the right-hand side 
$\mathcal{\vec{B}}$, i.e., we solve the penalized least squares problem 
\begin{equation}\label{LSQ}
\min_{\mathcal{\vec{X}}\in\mathbb{R}^{m\times 1\times n}} 
\left\{\|\mathcal{A*\vec{X}-\vec{B}}\|_F^2+\mu^{-1}\|\mathcal{L*\vec{X}}\|_F^2\right\},
\end{equation}
where $\mathcal{L}\in\mathbb{R}^{s\times m\times n}$ is a regularization operator and 
$\mu>0$ is a regularization parameter. This replacement is commonly referred to as 
Tikhonov regularization. Let $\mathcal{N}(\mathcal{M})$ denote the null space of the 
tensor $\mathcal{M}$ under $*$ and assume that $\mathcal{L}$ satisfies 
\begin{equation}\label{Null}
\mathcal{N}(\mathcal{A})\cap\mathcal{N}(\mathcal{L})=\{\mathcal{\vec{O}}\},
\end{equation}
where $\mathcal{\vec{O}}$ denotes an $m \times n$ zero matrix oriented laterally; see 
below. Then \eqref{LSQ} has a unique solution 
$\mathcal{\vec{X}}_\mu\in\mathbb{R}^{m\times 1\times n}$ for any $\mu>0$ (cf. Theorem \ref{th01}). The closeness of
$\mathcal{\vec{X}}_\mu$ to $\mathcal{\vec{X}}_\text{true}$ and the sensitivity of 
$\mathcal{\vec{X}}_\mu$ to the error $\mathcal{E}$ in $\mathcal{B}$ depends on the value 
of $\mu$. We determine $\mu$ by the discrepancy principle, which is described and analyzed in, 
e.g., \cite{EHN}. Application of the discrepancy principle requires that a bound 
\begin{equation}\label{errbd}
\|\mathcal{\vec{E}}\|_F \leq \delta
\end{equation}
be available. The parameter $\mu>0$ then is determined so that $\mathcal{\vec{X}}_\mu$ 
satisfies
\begin{equation}\label{discr}
\|\mathcal{\vec{B}}-\mathcal{A}*\mathcal{\vec{X}}_\mu \|_F=\eta \delta, 
\end{equation}
where $\eta>1$ is a user specified constant independent of $\delta>0$. It can be shown that 
$\mathcal{\vec{X}}_\mu\rightarrow\mathcal{\vec{X}}_\text{true}$ as $\delta\searrow 0$; see
\cite{EHN} for a proof in a Hilbert space setting.

Many other methods, including generalized cross validation (GCV) and the L-curve 
criterion, also can be used to determine the regularization parameter; see, e.g., 
\cite{CMRS,FRR,GHW,H2n,Ki,KR,RR} for discussions and illustrations for the situation when
$\mathcal{A}$ is a matrix and $\mathcal{\vec{B}}$ is a vector.

It is well known that a few steps of the (standard) Arnoldi process can be used to reduce 
a large matrix to a matrix of small size. The small matrix so obtained can be used to 
define a small 
Tikhonov regularization problems that is easy to solve; see \cite{CMRS,DMR,GNR,LR} for
discussions and illustrations. It is the purpose of the present paper to extend the 
(standard) matrix version of the Arnoldi process, described, e.g., in \cite{Sa}, to third 
order tensors using the t-product formalism. This gives us the t-Arnoldi process. 
Application of $\ell\geq 1$ steps of this process, generically, furnishes an orthonormal 
basis for the $\ell$-dimensional tensor Krylov (t-Krylov) subspace
\begin{equation}\label{tKry}
\mathbb{K}_\ell(\mathcal{A},\mathcal{\vec{B}})={\rm \text{t-span}}\left\{\mathcal{\vec{B}}, 
\mathcal{A}*\mathcal{\vec{B}},\mathcal{A}^2*\mathcal{\vec{B}},\dots,
\mathcal{A}^{\ell-1}*\mathcal{\vec{B}}\right\}.
\end{equation}
%\uu{
%The subspace \eqref{tKry} is generated by $\mathcal{A}$ and $\mathcal{\vec{B}}$, and is composed of $\ell$ lateral slices with 
%\begin{equation*}
%\mathcal{A}^{(i-1)}*\mathcal{\vec{B}} = \mathcal{A}^{(i-2)}*\mathcal{A*\vec{B}}, ~i = 2,3,\dots,\ell,
%\end{equation*}
%where $\mathcal{A}^0$ is the identity tensor denoted by $\mathcal{I}$; see \cite{GIJS,KBHH}. The 
%meaning of {\rm t-span} depends on the Arnoldi-type process used. We will comment on this in 
%Sections \ref{sec3} and \ref{sec4}.} 
The meaning of t-span is discussed in Sections \ref{sec3} and \ref{sec4}.
Each step of the t-Arnoldi process requires one tensor-matrix product evaluation with $\mathcal{A}$. Often fewer tensor-matrix product evaluations are required to solve Tikhonov minimization problems \eqref{LSQ} than when the t-product Golub-Kahan bidiagonalization (tGKB) process, described by Kilmer et al. \cite{KBHH} is used, because
each step of the latter demands two tensor-matrix product evaluations, one with 
$\mathcal{A}$ and one with $\mathcal{A}^T$, where the superscript $^T$ denotes transposition. 

We refer to our solution scheme for \eqref{LSQ} as the t-product Arnoldi-Tikhonov (tAT) 
regularization method. It is based on reducing the tensor
$\mathcal{A}\in\mathbb{R}^{m\times m\times n}$ to a small upper Hessenberg tensor. We also
describe a global tAT (G-tAT) method for the solution of \eqref{LSQ}. This method works 
with a data tensor slice $\mathcal{\vec{B}}\in\mathbb{R}^{m\times 1\times n}$ and is
closely related to the T-global Arnoldi-Tikhonov regularization method recently 
described by El Guide et al. \cite{GIJS}, which takes $\mathcal{L}$ equal to the identity 
tensor denoted by $\mathcal{I}$, determines the regularization parameter by the GCV method, 
and works with a general data tensor $\mathcal{B}\in\mathbb{R}^{m\times p\times n}$, $p>1$. 
Differently from the tAT method, the G-tAT and the T-global Arnoldi-Tikhonov 
regularization methods involve matricization of the tensor $\mathcal{A}$. Specifically, 
the G-tAT method first reduces $\mathcal{A}$ in \eqref{LSQ} to an upper Hessenberg matrix 
by carrying out a few steps of the global t-Arnoldi (G-tA) process. This process furnishes 
an orthonormal basis for a t-Krylov subspace \eqref{tKry}. It differs from the t-Arnoldi 
process in the choice of inner product. Algorithm \ref{Alg: 14} in Section \ref{sec5} 
provides the details of the G-tA process. Numerical examples with the t-Arnoldi and G-tA 
processes are presented in Section \ref{sec6}. The tAT and G-tAT methods based on these
processes determine the regularization parameter by the discrepancy principle. 

We also describe an extension of the (standard) generalized minimal residual (GMRES)
method proposed by Saad and Schultz \cite{SS} to third order tensors based on the 
t-product formalism. This extension will be referred to as the t-product GMRES (tGMRES) 
method. The tGMRES method for the solution of \eqref{1} computes iterates in t-Krylov 
subspaces of the form \eqref{tKry}; the $\ell$th approximate solution 
$\mathcal{\vec{X}}_\ell$ determined by tGMRES with initial approximate solution 
$\mathcal{\vec{X}}_0=\mathcal{\vec{O}}$ satisfies
\begin{equation}\label{GEM}
\|\mathcal{A}*\mathcal{\vec{X}}_\ell-\mathcal{\vec{B}} \|_F=
\min_{\mathcal{\vec{X}}\in\mathbb{K}_\ell(\mathcal{A},\mathcal{\vec{B}})}\|
\mathcal{A*\vec{X} - \vec{B}} \|_F, \;\;\; \ell = 1,2,\dots~.
\end{equation}

Another extension of the (standard) GMRES method by Saad and Schulz \cite{SS} for the 
solution of tensor equations is provided by the global tGMRES (G-tGMRES) method,
which is described in Subsection \ref{sec5.2}. This method is closely related to the 
T-global GMRES method recently presented by El Guide et al. \cite{GIJS}. The methods
differ in that the data for the G-tGMRES method is represented by a lateral slice 
$\mathcal{\vec{B}}$, while the data for T-global GMRES method is a general third order 
tensor $\mathcal{B}\in\mathbb{R}^{m\times p\times n}$, $p>1$. Moreover, our implementation
of the t-GMRES and G-tGMRES methods uses the discrepancy principle to determine when to 
terminate the iterations. Differently from 
the tGMRES method, the G-tGMRES and T-global GMRES methods involve matricization of the
tensor $\mathcal{A}$. While the tGMRES method is based on the t-Arnoldi process described 
in Section \ref{sec3}, the G-tGMRES method is based on the global t-Arnoldi (G-tA) 
process. 

Many other methods for solving \eqref{LSQ} and \eqref{GEM} that do not apply the t-product have been described in the literature; see, e.g., \cite{BNR,GIJB,IGJ,SE2}. These methods replace matrix-vector products by tensor-matrix products and involve matricization. A careful comparison of all these methods is outside the scope of the present paper. Here we note that computed examples of Section \ref{sec6} indicate that methods that avoid matricization often determine approximate solutions of higher quality than methods that involve matricization.

We also are interested in solving minimization problems analogous to \eqref{1}, in which
$\mathcal{\vec{B}}$ is replaced by a general third order tensor $\mathcal{B}$. This leads
to the Tikhonov minimization problem
\begin{equation}\label{multrhs}
\min_{\mathcal{X}\in\mathbb{R}^{m\times p\times n}}
\left\{\|\mathcal{A}*\mathcal{X}-\mathcal{B}\|_F^2+\mu^{-1}\|\mathcal{L*X}\|_F^2\right\}, 
\;\;\; \mathcal{B} \in \mathbb{R}^{m \times p \times n}, \;\;\; p>1.
\end{equation}
Besides our work \cite{LU}, no literature is available on solution methods for 
\eqref{LSQ} and \eqref{multrhs} for $\mathcal{L \neq I}$. The present paper focuses on 
developing tensor Arnoldi-Tikhonov-type methods for this situation.

Four methods for the solution of \eqref{multrhs} will be described. Three of them are 
based on the tAT and G-tAT methods applied to the lateral slices $\mathcal{\vec{B}}_j$, 
$j=1,2,\dots,p$, of $\mathcal{B}$, independently. The other method generalizes the 
T-global Arnoldi-Tikhonov regularization method recently presented by El Guide et al. 
\cite{GIJS} to allow for $\mathcal{L}\neq\mathcal{I}$. This method works with the lateral 
slices of the data tensor $\mathcal{B}$ simultaneously, and will be referred to as the
generalized global tAT (GG-tAT) method. 

A comparison of the solution methods for \eqref{multrhs} is presented in Section 
\ref{sec6}. Computed examples show the GG-tAT method to require less CPU time, but the
G-tAT method may yield higher accuracy. The fact that the GG-tAT requires less CPU time 
is to be expected since it uses larger chunks of data at a time. 

We remark that the G-tAT and GG-tAT methods belong to the AT$\_$BTF (Arnoldi-Tikhonov Based 
Tensor Format) family of methods recently described by Beik et al. \cite{BNR}. They 
involve flattening and require additional product definitions to the t-product. 

Finally, we will discuss a variant of the T-global GMRES method that recently has been
described by El Guide et al. \cite{GIJS} and is based on t-product formalism. We will 
refer to our variant as the generalized global tGMRES (GG-tGMRES) method. This method replaces the data tensor $\mathcal{\vec{B}}$ in \eqref{GEM} by a general third
order tensor $\mathcal{B}$ and determines iterates in t-Krylov subspaces 
$\mathbb{K}_\ell(\mathcal{A},\mathcal{B})$. The $\ell$th iterate $\mathcal{X}_\ell$ 
determined by the GG-tGMRES method with initial iterate 
$\mathcal{X}_0=\mathcal{O}\in\mathbb{R}^{m\times p\times n}$ solves
\begin{equation}\label{GGEM}
\|\mathcal{A}*\mathcal{X}_\ell-\mathcal{B} \|_F=
\min_{\mathcal{X}\in\mathbb{K}_\ell(\mathcal{A},\mathcal{B})}
\|\mathcal{A*\mathcal{X}-\mathcal{B}}\|_F, \;\;\; \ell = 1,2,\dots~.
\end{equation}
In the T-global GMRES method by El Guide et al. \cite{GIJS}, the iterations are terminated 
based on a residual Frobenius norm and a set tolerance that is independent of the error in 
$\mathcal{B}$. Differently from the T-global GMRES method, our approach for solving 
\eqref{GGEM} uses the discrepancy principle to determine the number of iterations to carry 
out with the GG-tGMRES method.

This paper is organized as follows. Section \ref{sec2} introduces notation and 
preliminaries associated with the t-product. Methods based on the t-Arnoldi process are 
described in Section \ref{sec3}. This includes Tikhonov regularization methods, one of 
which is based on a nested t-Krylov subspace, and
GMRES-type methods for the computation of approximate solutions of \eqref{1} and the
analogous minimization problem obtained by replacing the tensor slice $\mathcal{\vec{B}}$ by 
a third order tensor $\mathcal{B}$. Thus, we can consider color image and video restoration 
problems. For the former, $\mathcal{B}$ represents a blurred and noisy RGB image of dimension 
$m\times p\times 3$, while for gray-scale video restoration problems, $\mathcal{B}$ is of 
dimension $m \times p \times n$ with a sequence of $n$ consecutive blurred and noisy video 
frames. Section \ref{sec4} describes algorithms that are
based on the generalized global t-Arnoldi (GG-tA) process with data tensor $\mathcal{B}$. 
The algorithms of Section \ref{sec5} are obtained by modifying algorithms of Section 
\ref{sec4} to be applicable to each lateral slice of $\mathcal{B}$ separately. This allows us to consider, for instance, the restoration of gray-scale images. Section \ref{sec6} presents some numerical examples that illustrate the performance of 
these methods. Concluding remarks can be found in Section \ref{sec7}.

\section{Notation and Preliminaries}\label{sec2}
This section reviews results on the t-product introduced by Kilmer et al. \cite{KBHH,KM} and 
defines notation from \cite{KM,KB} to be used in the sequel.
In this paper, a \textit{tensor} is of third order, i.e., a three-dimensional array of 
real scalars denoted by the calligraphic script letters, say, 
$\mathcal{A}=[a_{ijk}]_{i,j,k=1}^{\ell,m,n}\in\mathbb{R}^{\ell\times m\times n}$ with real 
entries $a_{ijk}$. Matrices and vectors are second and first order tensors, respectively. 
We use capital letters to denote matrices, lower case letters to denote vectors, and bold 
face lower case letters to denote tube fibers (tubal scalars or tubes). A fiber of a 
third order tensor is a 1D section obtained by fixing two of the indices. Using MATLAB 
notation, $\mathcal{A}(:,j,k)$, $\mathcal{A}(i,:,k)$, and $\mathcal{A}(i,j,:)$ denote 
mode-1, mode-2, and mode-3 fibers, respectively. A slice of a third order tensor is a 2D 
section obtained by fixing one of the indices. With MATLAB notation, $\mathcal{A}(i,:,:)$,
$\mathcal{A}(:,j,:)$, and $\mathcal{A}(:,:,k)$ denote the $i$th horizontal, $j$th lateral,
and $k$th frontal slices, respectively. The $j$th lateral slice is also denoted by 
$\mathcal{\vec{A}}_j$. It is a tensor and will be referred to as a tensor column. 
Moreover, the $k$th frontal slice, which also will be denoted by $\mathcal{A}^{(k)}$, is a matrix. 

Given $\mathcal{A}\in\mathbb{R}^{\ell\times m\times n}$ with $\ell\times m$ frontal slices
$\mathcal{A}^{(i)}$, $i = 1,2,\dots,n$, the operator $\mathtt{unfold}(\mathcal{A})$ 
returns a block $\ell n \times m$ matrix made up of the faces $\mathcal{A}^{(i)}$ of 
$\mathcal{A}$. The $\mathtt{fold}$ operator folds back the unfolded $\mathcal{A}$, i.e.,
\begin{equation*}
\mathtt{unfold}(\mathcal{A}) = \begin{bmatrix}
\mathcal{A}^{(1)}\\
\mathcal{A}^{(2)}\\
\vdots\\
\mathcal{A}^{(n)}
\end{bmatrix}, \;\;\;\;\; \mathtt{fold(unfold(\mathcal{A})) = \mathcal{A}}.
\end{equation*}
The operator $\mathtt{bcirc}(\mathcal{A})$ generates an $\ell n\times mn$ block circulant
matrix with $\mathtt{unfold}(\mathcal{A})$ forming the first block column, 
\begin{equation*}
\mathtt{bcirc}(\mathcal{A}) = \begin{bmatrix}
\mathcal{A}^{(1)} & \mathcal{A}^{(n)} & \dots &\mathcal{A}^{(2)}\\
\mathcal{A}^{(2)} &\mathcal{A}^{(1)} & \dots &\mathcal{A}^{(3)}\\
\vdots & \vdots & \ddots & \vdots\\
\mathcal{A}^{(n)} & \mathcal{A}^{(n-1)} & \dots & \mathcal{A}^{(1)}
\end{bmatrix}.
\end{equation*} 

\begin{defn}{(t-product \cite{KM})}\label{defn 2.1}
Let $\mathcal{A}\in\mathbb{R}^{\ell\times m\times n}$ and 
$\mathcal{B}\in\mathbb{R}^{m\times p\times n}$. Then the t-product $\mathcal{A*B}$ is the
tensor $\mathcal{C}\in\mathbb{R}^{\ell\times p\times n}$ defined by
\begin{equation}\label{tenA}
\mathcal{C}:=\mathtt{fold}(\mathtt{bcirc}(\mathcal{A})\cdot\mathtt{unfold}(\mathcal{B})),
\end{equation}
where ``$\cdot$'' denotes the standard matrix-matrix product.
\end{defn}

We can view $\mathcal{C}$ in \eqref{tenA} as an $\ell\times p$ matrix of tubes oriented 
along the third dimension with its $(i,j)$th tube given by 
\begin{equation*}
\mathcal{C}(i,j,:)=\sum_{k=1}^p \mathcal{B}(i,k,:)\ast\mathcal{C}(k,j,:).
\end{equation*}
This shows that the t-product is analogous to matrix multiplication, except that 
multiplication between scalars is replaced by circular convolution between tubes.

The matrix $\mathtt{bcirc}(\mathcal{A})$ can be block diagonalized by the discrete Fourier
transform (DFT) matrix combined with the Kronecker product. Suppose that
$\mathcal{A}\in\mathbb{R}^{\ell\times m\times n}$ and let $F_n\in{\mathbb C}^{n\times n}$ 
denote the unitary DFT matrix. Then 
\begin{equation}
\bar{A} := \mathtt{blockdiag}(\widehat{\mathcal{A}}^{(1)},\widehat{\mathcal{A}}^{(2)},
\dots,\widehat{\mathcal{A}}^{(n)})=(F_n\otimes I_\ell)\cdot\mathtt{bcirc}
(\mathcal{A})\cdot(F_n^*\otimes I_m),
\label{flop}
\end{equation}
where $\otimes$ is the Kronecker product and $F_n^*$ denotes the conjugate transpose of 
$F_n$. The matrix $\bar{A}$ is an $\ell n\times mn$ block diagonal matrix with 
$\ell\times m$ blocks $\widehat{\mathcal{A}}^{(i)}$, $i=1,2,\dots,n$. The matrices 
$\widehat{\mathcal{A}}^{(i)}$ are the frontal slices of the tensor 
$\widehat{\mathcal{A}}$ obtained by applying the discrete Fourier transform along each 
tube of $\mathcal{A}$. 
%They satisfy the symmetry conditions:
%\begin{equation}\label{floor}
% \left\{
%\begin{array}{ll}
%\widehat{\mathcal{A}}^{(1)} \in \mathbb{R}^{m \times m}\\
%{\tt conj}(\widehat{\mathcal{A}}^{(i)}) = \widehat{\mathcal{A}}^{(n-i+2)},\;\;\; 
%i = 2,3,\dots,\lfloor\frac{n+1}{2}\rfloor,
%\end{array}
%\right.
%\end{equation}
%where $\lfloor t\rfloor$ denotes the nearest integer less than or equal to $t$; see 
%\cite{KBHH,FCLLY}. 
We remark that 
\[
\| \mathcal{A}\|_F= \frac{1}{\sqrt{n}} \|\bar{A}\|_F.
\]

The t-product is a natural extension of matrix multiplication for third order tensors 
\cite{KM}. Higher order tensors allow the definition of analogues of the t-product; 
see \cite{MSL}. Matrix algorithms for QR and SVD factorizations have analogues for 
third order tensors; see Kilmer et al. \cite{KBHH}.

We may choose to evaluate $\mathcal{A*B}$ according to Definition 
\ref{defn 2.1} if the tensors $\mathcal{A}$ and $\mathcal{B}$ are sparse. For general 
tensors $\mathcal{A}\in\mathbb{R}^{\ell\times m\times n}$ and 
$\mathcal{B}\in\mathbb{R}^{m\times p\times n}$, the t-product $\mathcal{A*B}$ can be 
computed efficiently by using the transformation \eqref{flop}, i.e., 
\begin{equation}\label{altp}
\mathcal{A*B}=\mathtt{fold}\left((F_n^* \otimes I_\ell)\bar{A}(F_n\otimes I_m) \cdot
\mathtt{unfold}(\mathcal{B})\right).
\end{equation}
The right-hand side of \eqref{flop} can be evaluated in $\mathcal{O}(\ell mn\log_2(n))$ 
arithmetic floating point operations (flops) using the fast Fourier transform (FFT); see 
\cite{KM}. 

The t-product is readily computed in MATLAB. We often will use the superscript\; $\widehat{}$\; to denote 
objects that are obtained by taking the FFT along the third dimension. Using MATLAB 
notation, let $\mathcal{\widehat{C}}:=\mathtt{fft}(\mathcal{C},[\;],3)$ be the tensor 
obtained by applying the FFT to $\mathcal{C}$ along the third dimension. Then the 
t-product $\mathcal{A}*\mathcal{B}$ can be computed by first taking the FFT along the 
tubes of $\mathcal{A}$ and $\mathcal{B}$ to get 
$\mathcal{\widehat{A}}=\mathtt{fft}(\mathcal{A},[\;],3)$ and 
$\mathcal{\widehat{B}} = \mathtt{fft}(\mathcal{B},[\;],3)$, followed by a
matrix-matrix product of each pair of the frontal slices of $\mathcal{\widehat{A}}$ and $\mathcal{\widehat{B}}$,
\[
\mathcal{\widehat{C}}(:,:,i)=\mathcal{\widehat{A}}(:,:,i)\cdot\mathcal{\widehat{B}}(:,:,i),
\;\; i=1,2,\dots,n,
\]
and then taking the inverse FFT along the third dimension to obtain 
$\mathcal{C} = \mathtt{ifft}(\mathcal{\widehat{C}},[\;],3)$. 
%We remark that by using 
%\eqref{floor}, we can save computational time in the Fourier domain by evaluating 
%\eqref{element} as follows:
%\begin{equation}\label{toolbox}
%\mathcal{\widehat{C}}(:,:,i) =  \left\{
%\begin{array}{ll}
%\mathcal{\widehat{A}}(:,:,i) \cdot \mathcal{\widehat{B}}(:,:,i), ~~~ i=2,3,\dots,
%\lceil \frac{n+1}{2} \rceil\\
%{\tt conj}(\widehat{\mathcal{C}}^{(n-i+2)}),~~~~~ i=\lceil\frac{n+1}{2}\rceil + 1,\dots,n,\\
%\end{array},
%\right.
%\end{equation}
%where $\lceil t \rceil$ denotes the nearest integer greater than or equal to $t$. 
The t-product \eqref{altp} can be computed by using the MATLAB tensor-tensor product 
toolbox\footnote{\url{https//github.com/canyilu/tproduct}}; see \cite{FCLLY}. Certain
symmetry properties can be utilized during the computations. This is done in the 
computations reported in Section \ref{sec6}.

%Express the tensor $\mathcal{A}$ in terms of its tensor columns, i.e., 
%\[
%\mathcal{A} = [\mathcal{\vec{A}}_1, \mathcal{\vec{A}}_2, \dots, \mathcal{\vec{A}}_m] \in \mathbb{R}^{l \times m\times n}, \;\; \mathcal{\vec{A}}_j \in \mathbb{R}^{l \times 1 \times n}, \;\; j=1,2,\dots, m. 
%\]
%Then following Newman et al. \cite{NKH}, we define the range of the tensor $\mathcal{A}$ denoted by $\mathcal{R(A)}$ as the \textit{t-linear span} of the tensor columns of $\mathcal{A}$ given by
%\begin{equation*}
%\mathcal{R(A)} = \{ \mathcal{\vec{A}}_1*\mathbf{x}_1 + \cdots + \mathcal{\vec{A}}_m*\mathbf{x}_m\; | \;\mathbf{x}_j \in \mathbb{R}^{1 \times 1 \times n} \};
%\end{equation*}
%see Kilmer et al. \cite{KBHH} for more details on the range and null space of $\mathcal{A}$. The action
%of a third order tensor on oriented matrices, using the t-product, defines a linear 
%transformation. In particular, if $T\mathcal{(\vec{X}) = A*\vec{X}}$, then 
%$T:\mathbb{R}^{m \times 1 \times n} \rightarrow \mathbb{R}^{l \times 1 \times n}$ is a
%t-linear operator with the property that
%\[
%T(\mathcal{\vec{X}}*\mathbf{c}+ \mathcal{\vec{Y}}*\mathbf{d}) = 
%T(\mathcal{\vec{X}})*\mathbf{c}+T(\mathcal{\vec{Y}})*\mathbf{d}
%\]
%for arbitrary tubal scalars $\mathbf{c, d}$ and arbitrary 
%$\mathcal{\vec{X}},\mathcal{\vec{Y}}\in \mathbb{R}^{m\times 1\times n}$; see Braman 
%\cite{B} for a proof. If $l = m$, then $T$ will be an invertible liner operator when $\mathcal{A}$ is invertible.

Let $\mathcal{A}\in\mathbb{R}^{\ell\times m\times n}$. The tensor transpose 
$\mathcal{A}^T\in\mathbb{R}^{m\times\ell\times n}$ is the tensor obtained by transposing 
each one of the frontal slices of $\mathcal{A}$, and then reversing the order of the 
transposed frontal slices 2 through $n$; see \cite{KM}. The tensor transpose has similar 
properties as the matrix transpose. For instance, if $\mathcal{A}$ and $\mathcal{B}$ are 
two tensors such that $\mathcal{A*B}$ and $\mathcal{B}^T*\mathcal{A}^T$ are defined, then 
$(\mathcal{A*B})^T=\mathcal{B}^T*\mathcal{A}^T$.

The identity tensor $\mathcal{I}\in\mathbb{R}^{m\times m\times n}$ is a tensor, whose 
first frontal slice, $\mathcal{I}^{(1)}$, is the $m\times m$ identity matrix and all other
frontal slices, $\mathcal{I}^{(i)}$, $i = 2,3, \dots, n$, are zero matrices; see 
\cite{KM}. 

The concept of orthogonality is well defined under the t-product formalism; see Kilmer and
Martin \cite{KM}. A tensor $\mathcal{Q}\in\mathbb{R}^{m\times m\times n}$ is said to be 
orthogonal if $\mathcal{Q}^T*\mathcal{Q} = \mathcal{Q}*\mathcal{Q}^T = \mathcal{I}$. 
Analogously to the columns of an orthogonal matrix, the lateral slices of an orthogonal
tensor $\mathcal{Q}$ are orthonormal, i.e.,
\begin{equation*}
\mathcal{Q}^T(:,i,:)*\mathcal{Q}(:,j,:) = \left\{
        \begin{array}{ll}
            {\bf e}_1 & i=j,\\
            {\bf 0} & i\neq j,
        \end{array}
    \right.
\end{equation*}
where ${\bf e}_1\in\mathbb{R}^{1\times 1\times n}$ is a tubal scalar whose $(1,1,1)$ entry
equals $1$ and the remaining entries vanish. It is shown in \cite{KM} that if 
$\mathcal{Q}$ is an orthogonal tensor, then 
\begin{equation}\label{lem: 2.1}
\|\mathcal{Q*A}\|_F=\|\mathcal{A}\|_F.
\end{equation}
The tensor $\mathcal{Q}\in\mathbb{R}^{\ell\times m\times n}$ with $\ell>m$ is said to be
partially orthogonal if $\mathcal{Q}^T*\mathcal{Q}$ is well defined and equal to the 
identity tensor $\mathcal{I}\in\mathbb{R}^{m\times m\times n}$; see \cite{KM}.

A tensor $\mathcal{A}\in\mathbb{R}^{m\times m\times n}$ is said to have an inverse, 
denoted by $\mathcal{A}^{-1}$, provided that $\mathcal{A}*\mathcal{A}^{-1}=\mathcal{I}$ 
and $\mathcal{A}^{-1}*\mathcal{A}=\mathcal{I}$. Moreover, a tensor is said to be f-diagonal 
if each frontal slice of the tensor is a diagonal matrix; see \cite{KM}. 

The tensor singular value decomposition (tSVD) of 
$\mathcal{A}\in\mathbb{R}^{\ell \times m\times n}$, introduced by Kilmer and Martin \cite{KM},
is given by
\begin{equation*}
\mathcal{A} = \mathcal{U}*\mathcal{S}*\mathcal{V}^T,
\end{equation*}
where $\mathcal{U}\in\mathbb{R}^{\ell \times \ell \times n}$ and 
$\mathcal{V}\in\mathbb{R}^{m\times m\times n}$ are orthogonal tensors, and the tensor 
\[
\mathcal{S}={\rm diag}[\mathbf{s}_1,\mathbf{s}_2,\dots,\mathbf{s}_{\min\{\ell ,m\}}]\in
\mathbb{R}^{\ell \times m\times n}
\]
is f-diagonal with singular tubes ${\mathbf{s}}_j\in\mathbb{R}^{1\times 1\times n}$, 
$j =1,2,\dots,\min\{\ell, m\}$, ordered according to
\[
\|\mathbf{s}_1\|_F\geq\|\mathbf{s}_2\|_F\geq\cdots\geq\|\mathbf{s}_{\min\{\ell,m\}}\|_F.
\]
The number of nonzero singular tubes of $\mathcal{A}$ is referred to as the tubal rank of $\mathcal{A}$; see Kilmer et al. \cite{KBHH}. The singular tubes of $\mathcal{A}$ are analogues of the singular values of a matrix $A$. In linear discrete ill-posed problem that require the solution of a linear system
of equations or least squares problem with a matrix $A$, this matrix has many singular values of 
different orders of magnitude close to zero. Definition \ref{ipd} describes linear discrete 
ill-posed tensor problems.

\begin{defn}\label{ipd}
The tensor least squares problems \eqref{1} is said to be a linear discrete ill-posed problem for third order tensors under $*$ if $\mathcal{A}$ has ill-determined tubal rank, i.e., the Frobenius norm of the singular tubes of $\mathcal{A}$ decays rapidly to zero with increasing index, and there are many nonvanishing singular tubes of tiny Frobenius norm of different orders of magnitude.
\end{defn}

We remark that this definition is not in terms of the frontal slices 
$\mathcal{A}^{(i)}, i = 1,2,\dots, n,$ of $\mathcal{A}$, but describes a property of the whole tensor $\mathcal{A}$, i.e., of the singular tubes of $\mathcal{A}$. The singular tubes are computed by finding the singular value decomposition of each frontal slice $\mathcal{\widehat{A}}^{(i)}$, $i = 1,2,\dots,n$, of $\mathcal{\widehat{A}}$ in the Fourier domain; see \cite{KM} for details.

The norm of a nonzero tensor column $\mathcal{\vec{X}}\in\mathbb{R}^{m\times 1\times n}$ 
is defined as
\[
\|\mathcal{\vec{X}}\| := \frac{\|\vec{\mathcal{X}}^T*\vec{\mathcal{X}}\|_F}
{\|\vec{\mathcal{X}}\|_F},
\]
and $\|\mathcal{\vec{X}}\| = 0$ if $\mathcal{\vec{X}} = \mathcal{\vec{O}}$; see 
\cite{KBHH} for details. The Frobenius norm of a tensor column 
$\mathcal{\vec{X}}\in\mathbb{R}^{m\times 1\times n}$ is given by
\[
\|\mathcal{\vec{X}}\|_F^2=\left(\mathcal{\vec{X}}^T*\mathcal{\vec{X}}\right)_{(:,:,1)};
\]
see \cite{KBHH}. Thus, the square of the Frobenius norm of $\mathcal{\vec{X}}$ is the 
first frontal face of the tube 
$\mathcal{\vec{X}}^T*\mathcal{\vec{X}}\in\mathbb{R}^{1\times 1\times n}$.

Algorithm \ref{Alg: 01}, which takes a nonzero tensor 
$\mathcal{\vec{X}} \in \mathbb{R}^{m\times 1\times n}$ and returns a normalized tensor 
$\mathcal{\vec{V}}\in\mathbb{R}^{m\times 1\times n}$ and a tubal scalar 
$\mathbf{a}\in\mathbb{R}^{1\times 1\times n}$ such that 
\begin{equation*}
 \mathcal{\vec{X}} = \mathcal{\vec{V}}* \mathbf{a} \;\;\; \text{and} \;\;\; 
 \|\mathcal{\vec{V}}\| = 1,
\end{equation*}
is important in the sequel. Note that the tubal scalar $\mathbf{a}$ might not be 
invertible; see \cite{KBHH} for details. We mention that $\mathbf{a}$ is invertible if 
there is a tubal scalar $\mathbf{b}$ such that $\mathbf{a*b} = \mathbf{b*a} = {\bf e}_1$. 
The scalar $\mathbf{a}^{(j)}$ is the $j$th face of the $1\times 1\times n$ tubal scalar 
$\mathbf{a}$, while $\mathcal{\vec{V}}^{(j)}$ is a vector with $m$ entries, and is the 
$j$th frontal face of $\mathcal{\vec{V}}\in\mathbb{R}^{m\times 1\times n}$. The call of 
the MATLAB function $\mathtt{randn}(m,1)$ in Algorithm \ref{Alg: 01} generates a 
pseudo-random $m$-vector with normally distributed entries with zero mean and variance 
one. In Algorithm \ref{Alg: 01} and elsewhere in this paper, $\|\cdot\|_2$ denotes the 
Euclidean vector norm.

\vspace{.3cm}
\begin{algorithm}[H]
\SetAlgoLined
\KwIn{$\mathcal{\vec{X}} \in \mathbb{R}^{m \times 1 \times n} \neq \mathcal{\vec{O}}$}
\KwOut{$\mathcal{\vec{V}}$, $\mathbf{a}$ with $\|\mathcal{\vec{V}}\| = 1$}
$\mathcal{\vec{V}} \leftarrow \mathtt{fft}(\mathcal{\vec{X}},[\;],3)$\\
\For{$j = 1$ \bf{to} $n$}{
$\mathbf{a}^{(j)} \gets \|\mathcal{\vec{V}}^{(j)}\|_2 \;\;\;$ ($\mathcal{\vec{V}}^{(j)}$ is a vector)\\
\eIf{$\mathbf{a}^{(j)} > \mathtt{tol}$}{
$\mathcal{\vec{V}}^{(j)} \gets \frac{1}{\mathbf{a}^{(j)}} \mathcal{\vec{V}}^{(j)}$\\
}{
$\mathcal{\vec{V}}^{(j)} \gets \mathtt{randn}(m,1); \;\; \mathbf{a}^{(j)} \gets \|\mathcal{\vec{V}}^{(j)}\|_2;\;\; \mathcal{\vec{V}}^{(j)} \gets \frac{1}{\mathbf{a}^{(j)}} \mathcal{\vec{V}}^{(j)}; \;\; \mathbf{a}^{(j)} \gets 0$\\
 }
}
$\mathcal{\vec{V}} \gets \mathtt{ifft}(\mathcal{\vec{V}},[\;],3); \;\; \mathbf{a} \gets \mathtt{ifft}(\mathbf{a},[\;],3)$
\caption{Normalize \cite{KBH}}
 \label{Alg: 01}
\end{algorithm}\vspace{.3cm}

The t-product based tensor QR (tQR) factorization implemented by Algorithm \ref{Alg: 2} is
described by Kilmer et al. \cite{KBHH}. Let 
$\mathcal{A}\in\mathbb{R}^{\ell\times m\times n}$. Then its tQR factorization is given by
\begin{equation*}
\mathcal{A = Q*R},  
\end{equation*} 
where the tensor $\mathcal{Q}\in\mathbb{R}^{\ell\times m\times n}$ is partially orthogonal
and the tensor $\mathcal{R}\in\mathbb{R}^{m\times m\times n}$ is f-upper triangular 
(i.e., each face is upper triangular). 

\vspace{.3cm}
\begin{algorithm}[H]
\SetAlgoLined
\KwIn{$\mathcal{A}\in\mathbb{R}^{\ell\times m\times n}$, $\ell\geq m$}
\KwOut{$\mathcal{Q}\in\mathbb{R}^{\ell\times m\times n},\; \mathcal{R}\in
\mathbb{R}^{m\times m\times n}$ such that $\mathcal{A}=\mathcal{Q*R}$}
$\widehat{\mathcal{A}}\leftarrow\mathtt{fft}(\mathcal{A},[\;],3)$\\
 \For {$i=1$\bf{to} $n$}{
 Factor $\widehat{\mathcal{A}}(:,:,i)=QR$, where $Q$ is unitary\\
 $\widehat{\mathcal{Q}}(:,:,i) \leftarrow Q, \;\;\; \widehat{\mathcal{R}}(:,:,i) \leftarrow R$
 }
 $\mathcal{Q}\leftarrow \mathtt{ifft}(\widehat{\mathcal{Q}},[\;],3), \;\;\; \mathcal{R}\leftarrow \mathtt{ifft}(\widehat{\mathcal{R}},[\;],3)$
 \caption{{tQR factorization \cite{KBHH}}}
 \label{Alg: 2}
\end{algorithm}\vspace{.3cm}

We introduce additional definitions used by El Guide et al. \cite{GIJS}. They will be 
needed when discussing the G-tAT, GG-tAT, G-tGMRES and GG-tGMRES methods in Sections \ref{sec4} and \ref{sec5}. Let
\[
\mathbb{C}_k:=[\mathcal{C}_1,\mathcal{C}_2,\dots,\mathcal{C}_k]\in
\mathbb{R}^{m\times kp\times n}, \;\;\;\;\; 
\mathcal{C}_k:=[\mathcal{\vec{C}}_1,\mathcal{\vec{C}}_2,\dots,\mathcal{\vec{C}}_k]\in
\mathbb{R}^{m\times k\times n},
\]
where $\mathcal{C}_j\in\mathbb{R}^{m\times p\times n}$ and 
$\mathcal{\vec{C}}_j\in\mathbb{R}^{m\times 1\times n}$. Suppose that 
$y = [y_1, \dots, y_k]^T\in\mathbb{R}^k$. Then El Guide et al. defined the product $\circledast$ as
\[
\mathbb{C}_k\circledast y=\sum_{j=1}^k y_j\mathcal{C}_j, \;\; \;\;\; 
\mathcal{C}_k\circledast y=\sum_{j=1}^k y_j \mathcal{\vec{C}}_j.
\]
It can be shown that for orthogonal tensors $\mathbb{Q}\in\mathbb{R}^{m\times kp\times n}$
and $\mathcal{Q}\in\mathbb{R}^{m\times k\times n}$, one has
\begin{equation}\label{norm2F}
\|\mathbb{Q}\circledast y\|_F=\|y\|_2, \;\;\;\;\; \|\mathcal{Q}\circledast y\|_F=\|y\|_2;
\end{equation}
see \cite{GIJS} for details.

Consider the tensors $\mathcal{C}=[c_{ijk}]$ and $\mathcal{D}=[w_{ijk}]\in
\mathbb{R}^{m\times p\times n}$ with lateral slices 
$\mathcal{\vec{C}}=[c_{i1k}]$ and 
$\mathcal{\vec{D}}=[d_{i1k}]\in\mathbb{R}^{m\times 1\times n}$, respectively. Define the 
scalar products
\begin{equation*}
\langle\mathcal{C},\mathcal{D}\rangle=\sum_{i=1}^m\sum_{j=1}^p\sum_{k=1}^nc_{ijk}d_{ijk},
\;\;\;\;\;\; \langle\mathcal{\vec{C}},\mathcal{\vec{D}}\rangle=
\sum_{i=1}^m\sum_{k=1}^n c_{i1k}d_{i1k}.
\end{equation*}
Let
\begin{equation}\label{nota}
\begin{split}
\mathbb{A}:=[\mathcal{A}_1,\mathcal{A}_2,\dots,\mathcal{A}_m]\in
\mathbb{R}^{\ell\times km\times n}, \;\;\;\;\; 
\mathbb{B}:=[\mathcal{B}_1,\mathcal{B}_2,\dots,\mathcal{B}_p]\in
\mathbb{R}^{\ell\times kp\times n},\\ 
\mathcal{A}:=[\mathcal{\vec{A}}_1,\mathcal{\vec{A}}_2,\dots,\mathcal{\vec{A}}_m]\in
\mathbb{R}^{\ell\times m\times n}, \;\;\;\;\;  
\mathcal{B}:=[\mathcal{\vec{B}}_1,\mathcal{\vec{B}}_2,\dots,\mathcal{\vec{B}}_p]\in
\mathbb{R}^{\ell\times p\times n},
\end{split}
\end{equation}
where $\mathcal{A}_i\in\mathbb{R}^{\ell\times k\times n}$, 
$\mathcal{\vec{A}}_i\in\mathbb{R}^{\ell\times 1\times n}$, $i=1,2,\dots,m$, and  
$\mathcal{B}_j\in\mathbb{R}^{\ell\times k\times n}$, 
$\mathcal{\vec{B}}_j\in\mathbb{R}^{\ell\times 1\times n}$, $j=1,2,\dots,p$. Following El 
Guide et al. \cite{GIJS}, we define the T-diamond products 
$\mathbb{A}^T\Diamond\mathbb{B}$ and $\mathcal{A}^T\Diamond\mathcal{B}$. They yield 
$m\times p$ matrices with entries
\begin{equation*}
[\mathbb{A}^T\Diamond\mathbb{B}]_{ij}=\langle \mathcal{A}_i, \; \mathcal{B}_j\rangle, 
\;\;\;\;\; 
[\mathcal{A}^T\Diamond\mathcal{B}]_{ij}=\langle \mathcal{\vec{A}}_i, \; 
\mathcal{\vec{B}}_j\rangle, \;\;\;  i = 1,2,\dots,m, \;\; j= 1,2,\dots,p.
\end{equation*}

The generalized global tensor QR (GG-tQR) factorization is described in \cite{LU} and 
implemented by Algorithm \ref{Alg: tggqr}. Given $\mathbb{A}$ in \eqref{nota}, this
factorization is defined by 
\[
\mathbb{A}=\mathbb{Q}\circledast R,
\]
where $R\in\mathbb{R}^{m \times m}$ is an upper triangular matrix, and the tensor
$\mathbb{Q}\in\mathbb{R}^{\ell\times km\times n}$ with $\ell \geq k$ has $k$ partially orthogonal tensor columns such 
that 
\[
\mathbb{Q}^T \Diamond \mathbb{Q} = I_m,
\]
where $I_m$ is the $m \times m$ identity matrix.

\vspace{.3cm}
\begin{algorithm}[H]
\SetAlgoLined
\KwIn{$\mathbb{A} := [\mathcal{A}_1,\mathcal{A}_2,\dots,\mathcal{A}_m]\in 
\mathbb{R}^{\ell\times km\times n}$, $\mathcal{A}_j\in\mathbb{R}^{\ell\times k\times n}$,
$j=1,2,\dots,m$, $\ell\geq k$}
\KwOut{$\mathbb{Q} := [\mathcal{Q}_1,\mathcal{Q}_2,\dots,\mathcal{Q}_m]\in 
\mathbb{R}^{\ell\times km \times n}$, $R = (r_{ij})\in\mathbb{R}^{m \times m}$ such that 
$\mathbb{A}=\mathbb{Q}\circledast R$ and $\mathbb{Q}^T\Diamond\mathbb{Q} = I_m$}
Set $r_{11} \leftarrow\langle\mathcal{A}_1,\mathcal{A}_1\rangle^{1/2}$,
$\mathcal{Q}_1\leftarrow\frac{1}{r_{11}}\mathcal{A}_1$\\
 \For {$j  = 1,2,\dots,m$}{
$\mathcal{W}\rightarrow \mathcal{A}_j$\\
 \For {$i  = 1,2,\dots,j-1$}{
$r_{ij}\leftarrow\langle\mathcal{Q}_i,\mathcal{W}\rangle$\\
$\mathcal{W}\leftarrow\mathcal{W}-r_{ij}\mathcal{Q}_i$\\
 }
$r_{jj}\leftarrow\langle\mathcal{W},\mathcal{W}\rangle^{1/2}$\\
$\mathcal{Q}_j \leftarrow\mathcal{W}/r_{jj}$}
\caption{Generalized global tQR (GG-tQR) factorization \cite{LU}}
\label{Alg: tggqr}
\end{algorithm}\vspace{.3cm}

We also will need a special case of the GG-tQR factorization, which works with each
lateral slice $\mathcal{\vec{A}}_i$, $i=1,2,\dots,m$, of tensor $\mathcal{A}$ in \eqref{nota}. 
This
factorization method is implemented by Algorithm \ref{Alg: tgqr}; it is also described
in \cite{LU}, and is there referred to as the global tQR (G-tQR) factorization method. 

\vspace{.3cm}
\begin{algorithm}[H]
\SetAlgoLined
\KwIn{$\mathcal{A}:=[\mathcal{\vec{A}}_1,\mathcal{\vec{A}}_2,\dots,\mathcal{\vec{A}}_m]\in
\mathbb{R}^{\ell\times m\times n}$, $\mathcal{\vec{A}}_j\in
\mathbb{R}^{\ell\times 1\times n}$, $j = 1,2,\dots,m$, $\ell\geq m$}
\KwOut{$\mathcal{Q}:=[\mathcal{\vec{Q}}_1,\mathcal{\vec{Q}}_2,\dots,\mathcal{\vec{Q}}_m]
\in\mathbb{R}^{\ell\times m\times n}$, $\mathcal{\vec{Q}}_j\in
\mathbb{R}^{\ell\times 1\times n}$, $\bar{R}=[r_{ij}]\in\mathbb{R}^{m\times m}$ such that
$\mathcal{A}=\mathcal{Q}\circledast\bar{R}$ and $\mathcal{Q}^T\Diamond\mathcal{Q}=I_m$}
$r_{11}\leftarrow\langle\mathcal{\vec{A}}_1,\mathcal{\vec{A}}_1\rangle^{1/2}$, 
$\mathcal{\vec{Q}}_1\leftarrow\frac{1}{r_{11}}\mathcal{\vec{A}}_1$\\
\For {$j=1,2,\dots,m$}{
$\mathcal{\vec{W}}\leftarrow\mathcal{\vec{A}}_j$\\
\For {$i=1,2,\dots,j-1$}{
$r_{ij}\leftarrow\langle\mathcal{\vec{Q}}_i,\mathcal{\vec{W}}\rangle$\\
$\mathcal{\vec{W}}\leftarrow\mathcal{\vec{W}} - r_{ij} \mathcal{\vec{Q}}_i$\\
}
$r_{jj}\leftarrow\langle\mathcal{\vec{W}},\mathcal{\vec{W}}\rangle^{1/2}$\\
$\mathcal{\vec{Q}}_j\leftarrow\mathcal{\vec{W}}/r_{jj}$}
\caption{Global tQR (G-tQR) factorization \cite{LU}}
 \label{Alg: tgqr}
\end{algorithm}\vspace{.3cm}

We conclude this section with the definition of some tensor operators that are convenient 
to apply in Section \ref{sec6}. The matrix $X\in\mathbb{R}^{m\times n}$ is associated with
the tensor $\mathcal{\vec{X}}\in\mathbb{R}^{m\times 1\times n}$ by the $\mathtt{squeeze}$ 
and $\mathtt{twist}$ operators, defined by Kilmer et al. \cite{KBHH}, i.e.,
\[
\mathcal{\vec{X}}=\mathtt{twist}(X)\;\; {\rm and} \;\; 
X = \mathtt{squeeze}(\mathcal{\vec{X}}).
\]
Note that the $\mathtt{squeeze}$ operator is identical to the MATLAB squeeze function.

We also define the $\mathtt{multi}\_\mathtt{squeeze}$ and $\mathtt{multi}\_\mathtt{twist}$
operators that enable us to squeeze and twist a general third order tensor. The tensor
$\mathcal{C}\in\mathbb{R}^{m\times p\times n}$ is associated with 
$\mathcal{D}\in\mathbb{R}^{m\times n\times p}$ by 
\[
\mathcal{D} = \mathtt{multi}\_\mathtt{twist}(\mathcal{C}) ~~ {\rm and} ~~ 
\mathcal{C}=\mathtt{multi}\_\mathtt{squeeze}(\mathcal{D}),
\] 
where $\mathtt{multi}\_\mathtt{twist}(\mathcal{C})$ twists each of the frontal slices 
$\mathcal{C}^{(i)}$, $i=1,2,\dots,n$, of $\mathcal{C}$ by using the $\mathtt{twist}$ 
operator, and stacks them as lateral slices $\mathcal{\vec{D}}_i$, $i=1,2,\dots,n$, of 
$\mathcal{D}$. Moreover, the operator $\mathtt{multi}\_\mathtt{squeeze}(\mathcal{D})$ 
squeezes the lateral slices of $\mathcal{D}$ using the $\mathtt{squeeze}$ operator and 
stacks them as faces of $\mathcal{C}$.

\section{Methods based on the t-Arnoldi process}\label{sec3}
We first describe an algorithm for the t-Arnoldi process. This algorithm is applied in 
Subsections \ref{sec3.1} and \ref{sec3.2} to reduce the large-scale problem \eqref{1} to
a problem of small size.

Let $\mathcal{A}\in\mathbb{R}^{m\times m\times n}$. The t-Arnoldi process described by 
Algorithm \ref{Alg: tArn} (cf. the matrix version in \cite[Chapter 5]{Sa}) reduces the 
tensor $\mathcal{A}$ to an upper Hessenberg tensor (t-Hessenberg), whose every face is an upper Hessenberg matrix. 

\vspace{.3cm} 
\begin{algorithm}[H]
\SetAlgoLined
\KwIn{$\mathcal{A}\in\mathbb{R}^{m\times m\times n},\; 
\mathcal{\vec{B}}\in\mathbb{R}^{m\times 1\times n}\neq\mathcal{\vec{O}}$}
$[\mathcal{\vec{Q}}_1,{\bf z_1}]\leftarrow \text{Normalize}(\mathcal{\vec{B}})$ with 
${\bf z_1}$ invertible, and such that $\mathcal{\vec{B}}=\mathcal{\vec{Q}}_1*{\bf z_1}$ and
$\|\mathcal{\vec{Q}}_1\|=1$ \\
\For {$j=1,2,\dots,\ell$}{
$\mathcal{\vec{W}}\leftarrow\mathcal{A}*\mathcal{\vec{Q}}_j$\\
\For {$i= 1,2,\dots,j$}{
$\mathbf{h}_{ij}\leftarrow\mathcal{\vec{Q}}_i^T*\mathcal{\vec{W}}$\\
$
 \left\{
        \begin{array}{ll}
\mathcal{\vec{W}}\leftarrow\mathcal{\vec{W}}-\mathcal{\vec{Q}}_i*\mathbf{h}_{ij} \; 
(\text{no reorthogonalization})\\
\mathcal{\vec{W}}\leftarrow\mathcal{\vec{W}}-\mathcal{\vec{Q}}_i*\mathbf{h}_{ij}, \; 
\mathcal{\vec{W}}\leftarrow\mathcal{\vec{W}}-\sum\limits_{k=1}^{i} 
\mathcal{\vec{Q}}_k*(\mathcal{\vec{Q}}_k^T*\mathcal{\vec{W}}) \; 
(\text{with reorthogonalization})\\
    \end{array}
    \right.
$\\
}
 $[\mathcal{\vec{Q}}_{j+1},\mathbf{h}_{j+1, j}]\leftarrow$ 
 Normalize$(\mathcal{\vec{W}})$ with $\mathbf{h}_{j+1, j}$ invertible
}
\caption{The t-Arnoldi process}
\label{Alg: tArn}
\end{algorithm}\vspace{.3cm}

The t-Arnoldi process is said to {\it break down} if any of the subdiagonal tubal scalars 
$\mathbf{h}_{j+1,j}$ for $j=1,2,\ldots,\ell$ is not invertible. This is analogous to a
break down of the (standard) Arnoldi process. We will assume that the number of steps,
$\ell$, of the t-Arnoldi process is small enough to avoid break down, i.e., that $\ell$ 
is chosen small enough so that every subdiagonal tubal scalar $\mathbf{h}_{j+1,j}$ is
invertible for $j=1,2,\ldots,\ell$. This means, in particular, that the transformed tubal 
scalars $\widehat{\mathbf{h}}_{j+1,j}$ of $\mathbf{h}_{j+1,j}$ do not have zero Fourier 
coefficients. 

Algorithm \ref{Alg: tArn} produces the partial t-Arnoldi decomposition
\begin{equation}
\mathcal{A} * \mathcal{Q}_\ell = \mathcal{Q}_{\ell+1} * \mathcal{\bar{H}}_\ell, 
\label{3.1}
\end{equation}
where
\begin{center}
$\mathcal{\bar{H}}_\ell = \begin{bmatrix}
\mathbf{h}_{11}     & 			&			  &       \dots             &                                   \mathbf{h}_{1\ell}\\ 
\mathbf{h}_{21}     & \mathbf{h}_{22} & \\
 			 &\mathbf{h}_{32} & \mathbf{h}_{33}    &  &  \vdots \\
			&		          & 	 \ddots  	 & \ddots                  &                                         \\
			&		          &		            &		           	 \mathbf{h}_{\ell, \ell-1} & \mathbf{h}_{\ell,\ell}\\
			&	                     &			&		          	 		     & \mathbf{h}_{\ell+1,\ell}
\end{bmatrix}\in\mathbb{R}^{(\ell+1)\times\ell\times n}$
\end{center}
is of upper t-Hessenberg form. The lateral slices $\mathcal{\vec{Q}}_j$, 
$j=1,2,\dots,\ell$, of $\mathcal{Q}_\ell \in\mathbb{R}^{m\times \ell \times n}$ form an 
orthonormal tensor basis for the t-Krylov subspace \eqref{tKry}, where {\rm t-span} 
refers to the set of all tensor linear (t-linear) combinations, whose coefficients are 
tubal scalars, ${\bf c}_i \in \mathbb{R}^{1 \times 1 \times n}$, $i = 1,2,\dots,\ell$. 
Thus,
\be \label{ttKry}
\mathbb{K}_\ell(\mathcal{A},\mathcal{\vec{B}}) = 
\bigg\{ \mathcal{\vec{Z}} \in \mathbb{R}^{m \times 1 \times n},~ \mathcal{\vec{Z}} =
\sum_{i=1}^\ell (\mathcal{A}^{(i-1)}*\mathcal{\vec{B}})*{\bf c}_i \bigg\}, ~~ \mathcal{A}^{0} = \mathcal{I}.
\ee 
%where $\mathcal{\vec{B}}$ is normalized with a normalization invertible tubal scalar 
%(cf. Algorithm \ref{Alg: 01}). 
The t-Arnoldi process  generates an orthonormal tensor basis for the t-Krylov subspace 
\eqref{ttKry} by  applying the standard Arnoldi process to each frontal slice 
$\mathcal{\widehat{A}}^{(i)}$, $i=1,2,\dots,n$, of $\mathcal{\widehat{A}}$ simultaneously.
This process applies the normalization Algorithm 1 to the data tensor $\mathcal{\vec{B}}$.

We comment on the complexity of the standard Arnoldi and t-Arnoldi processes. Let 
$A \in \mathbb{R}^{m \times m}$ be a dense matrix and $1\leq\ell\ll m$ the number of steps
carried out by the standard Arnoldi process. Then this process requires 
$\mathcal{O}(\ell^2 m + \ell m^2)$ flops, since $\ell$ matrix-vector product with $A$ cost
$\mathcal{O}(\ell m^2)$ flops and $\mathcal{O}(\ell^2 m)$ flops are required for 
orthogonalization.

We implement the t-Arnoldi process with transformations to and from the Fourier 
domain.  For a dense tensor $\mathcal{A}\in \mathbb{R}^{m \times m \times n}$, application
of $1\leq\ell\ll m$ steps of this process requires application of $\ell$ steps of the 
standard (matrix) Arnoldi process to the frontal slices $\mathcal{A}^{(i)}$, 
$i=1,2,\dots,n$, of $\mathcal{A}$ simultaneously in the Fourier domain, and 
orthogonalization. Each transformation of $\mathcal{A}$ and $\mathcal{\vec{Q}}_j$ to and 
from the Fourier domain in step 3 of Algorithm \ref{Alg: tArn} costs 
$\mathcal{O}(m^2n\log(n))$ and $\mathcal{O}(mn\log(n))$ flops, respectively. Moreover, 
$\ell$ matrix-vector products between the faces of $\mathcal{A}$ and $\mathcal{\vec{Q}}_j$
in the Fourier domain cost $\mathcal{O}(\ell m^2)$ flops. For $n$ frontal slices, it has a
complexity of $\mathcal{O}(\ell m^2 n)$ flops in the Fourier domain. Similarly, the 
orthogonalization steps $4$-$7$ in the Fourier domain cost $\mathcal{O}(\ell^2 mn)$ flops 
for $n$ frontal slices. Note that it costs $\mathcal{O}(n\log(n))$ flops to transform each 
tubal scalar ${\bf h}_{ij}$ to and from the Fourier domain. Hence, the total flop count for 
carrying out $\ell$ steps of the t-Arnoldi process  in the Fourier domain  is $\mathcal{O}((\ell m^2+\ell^2 m)n)$
flops. The cost is the same for the G-tA process implemented by Algorithm \ref{Alg: 14} in
Section \ref{sec5}.

We will use the decomposition \eqref{3.1} to determine an approximate solution of the 
Tikhonov minimization problems \eqref{LSQ} and \eqref{multrhs} in Subsection \ref{sec3.1},
and of the minimization problems \eqref{GEM} and \eqref{GGEM} in Subsection \ref{sec3.2}.

\subsection{Tensor Arnoldi-Tikhonov Regularization Methods}\label{sec3.1}
This subsection discusses the computation of an approximate solution of the tensor 
Tikhonov regularization problem \eqref{LSQ} with the aid of the t-Arnoldi process. We 
describe how this process can be used in conjunction with the discrepancy principle 
\eqref{discr}, and show that the penalized least squares problem (1.3) has a unique 
solution $\mathcal{\vec{X}}_\mu$; see, e.g., \cite{CR} for a proof of the matrix case.

\begin{thm}\label{th01}
Let $\mu >0$ be the regularization parameter. The minimization problem \eqref{LSQ} has a 
unique solution 
\begin{equation}\label{Xmu}
\mathcal{\vec{X}}_\mu=(\mathcal{A}^T*\mathcal{A}+\mu^{-1}\mathcal{L}^T*\mathcal{L} )^{-1}*
\mathcal{A}^T*\mathcal{\vec{B}}
\end{equation}
that satisfies the normal equations 
\be \label{normeq}
(\mathcal{A}^T*\mathcal{A}+\mu^{-1}\mathcal{L}^T*\mathcal{L})*\mathcal{\vec{X}}= 
\mathcal{A}^T*\mathcal{\vec{B}}.
\ee
\end{thm}
\noindent
{\it Proof:}
The function 
\[
\mathcal{J}_\mu(\mathcal{\vec{X}}) := \|\mathcal{A}*\mathcal{\vec{X}} - \mathcal{\vec{B}}\|^2_F + \mu^{-1}\| \mathcal{L*\vec{X}}\|_F^2
\]
can be written as
\[
\mathcal{J}_\mu(\mathcal{\vec{X}}) = 
\bigg\| \begin{bmatrix} \mathcal{A}\\
\mu^{-1/2}\mathcal{L} \end{bmatrix}* \mathcal{\vec{X}} - \begin{bmatrix}
\mathcal{\vec{B}}\\
\mathcal{\vec{O}}
\end{bmatrix}\bigg\|_F^2,
\]
where 
\[
\begin{bmatrix} \mathcal{A}\\ \mu^{-1/2}\mathcal{L} \end{bmatrix} \in 
\mathbb{R}^{(m+s) \times m \times n},\quad 
\begin{bmatrix} \mathcal{\vec{B}}\\ \mathcal{\vec{O}} \end{bmatrix} \in 
\mathbb{R}^{(m+s) \times 1 \times n},\quad 
\mathcal{\vec{O}}\in \mathbb{R}^{s \times 1 \times n}.
\]
Thus, $\mathcal{\vec{X}}_\mu$ is a minimizer of $\mathcal{J}_\mu(\mathcal{\vec{X}})$ 
if and only if $\mathcal{\vec{X}}_\mu$ is the solution to the normal equations
\[ 
\begin{bmatrix} \mathcal{A}\\
\mu^{-1/2}\mathcal{L} \end{bmatrix}^T*\begin{bmatrix} \mathcal{A}\\
\mu^{-1/2}\mathcal{L} \end{bmatrix}*\mathcal{\vec{X}} = \begin{bmatrix} \mathcal{A}\\
\mu^{-1/2}\mathcal{L} \end{bmatrix}^T*\begin{bmatrix}
\mathcal{\vec{B}}\\
\mathcal{\vec{O}}
\end{bmatrix},
\]
which can be written as \eqref{normeq}. Due to \eqref{Null} the solution is unique.~~$\Box$

A similar formulation of \eqref{normeq} when $\mathcal{L = I}$ has been described by 
Kilmer et al. \cite{KBHH} and Martin et al. \cite{MSL}.

When the regularization operator $\mathcal{L}$ is the identity tensor, the solution 
\eqref{Xmu} simplifies to 
\begin{equation}\label{9}
\mathcal{\vec{X}}_{\mu}=(\mathcal{A}^T*\mathcal{A}+\mu^{-1}\mathcal{I})^{-1}*
\mathcal{A}^T*\mathcal{\vec{B}}.
\end{equation}
Using this expression for $\mathcal{\vec{X}}_\mu$, define the function
\begin{equation}\label{11}
\phi(\mu):=\|\mathcal{A}*\mathcal{\vec{X}}_\mu - \mathcal{\vec{B}}\|^2_F.
\end{equation}
Then equation \eqref{discr} (for $\mathcal{L}=\mathcal{I}$) can be written as
\begin{equation}\label{zerof}
\phi(\mu) = \eta^2 \delta^2.
\end{equation}
A zero-finder, such as bisection, Newton's method, or a related method \cite{BPR,RS}, can
be used to solve \eqref{zerof} for $\mu_{\rm discr}=\mu>0$. We assume here and below that 
$\delta>0$. Then $\mathcal{\vec{X}}_{\mu_{\rm discr}}$ satisfies the discrepancy principle
\eqref{discr} (when $\mathcal{L}=\mathcal{I}$). 

The following properties of $\phi$ are shown in \cite{LU}. We remark that while the 
solution \eqref{9} is meaningful for $\mu>0$ only, we may define $\phi(\mu)$ for 
$\mu\geq 0$ by continuity. 

\begin{prop}\label{prop1}
Assume that $\mathcal{A}^T*\mathcal{\vec{B}} \neq \mathcal{\vec{O}}$ and let $\phi(\mu)$ be given by 
\eqref{11} with $\mathcal{\vec{X}}_\mu$ defined by \eqref{9}. Then 
\[
\phi(\mu) = \left(\mathcal{\vec{B}}^T*(\mu\mathcal{A}*\mathcal{A}^T+\mathcal{I})^{-2}*
\mathcal{\vec{B}}\right)_{(:,:,1)},\qquad \mu>0,
\]
and $\phi(0)=\|\mathcal{\vec{B}}\|^2_F$. Moreover,
\[
\phi^\prime(\mu) < 0\;\;\; {\rm and} \;\;\; \phi^{\prime\prime}(\mu) > 0
\]
for $\mu>0$.
\end{prop}

\subsubsection{The tAT methods for the solution of \eqref{LSQ}}\label{sec3.11}
We develop the t-product Arnoldi-Tikhonov (tAT) regularization method for the 
approximate solution of least squares problems of the form \eqref{LSQ}. The method will be
used to illustrate the potential superiority of tensorizing as opposed to vectorizing or 
matricizing ill-posed tensor equations in general. This method will be generalized in 
Subsection \ref{sec3.12} to the least squares problems \eqref{multrhs} with a general data
tensor $\mathcal{B}$.

Let $\mathcal{\vec{X}}=\mathcal{Q}_\ell*\mathcal{\vec{Y}}$ for some 
$\mathcal{\vec{Y}}\in\mathbb{R}^{\ell\times 1\times n}$ and substitute the 
decomposition \eqref{3.1} into \eqref{LSQ}. This yields 
\begin{equation}
\min_{\mathcal{\vec{Y}} \in \mathbb{R}^{\ell \times 1 \times n}}
\{ \|\mathcal{\bar{H}}_\ell*\mathcal{\vec{Y}} - \mathcal{Q}_{\ell+1}^T*\mathcal{\vec{B}}\|^2_F 
+ \mu^{-1}\|\mathcal{L}*\mathcal{Q}_\ell*\mathcal{\vec{Y}}\|^2_F\}.
\label{4.2}
\end{equation}
Using the fact that $\mathcal{\vec{B}} = \mathcal{\vec{Q}}_1*\mathbf{z}_1$ (cf. Algorithm
\ref{Alg: tArn}), we obtain 
\begin{equation}\label{QTB}
\mathcal{Q}_{\ell+1}^T*\mathcal{\vec{B}}=\vec{e}_1*\mathbf{z}_1\in
\mathbb{R}^{(\ell+1)\times 1\times n},
\end{equation}
where the $(1,1,1)$th entry of $\vec{\mathit{e}}_1\in\mathbb{R}^{m\times 1\times n}$ 
equals $1$ and the remaining entries vanish. Substitute \eqref{QTB} into \eqref{4.2} to 
obtain
\begin{equation}
\min_{\mathcal{\vec{Y}}\in \mathbb{R}^{\ell\times 1\times n}}
\{ \|\mathcal{\bar{H}}_\ell*\mathcal{\vec{Y}} - \vec{\mathit{e}}_1*\mathbf{z}_1\|^2_F + 
\mu^{-1}\|\mathcal{L}*\mathcal{Q}_\ell*\mathcal{\vec{Y}}\|^2_F\}.
\label{4.4}
\end{equation}

In the computed examples of Section \ref{sec6}, we use the regularization operators 
$\mathcal{L}_1\in\mathbb{R}^{(m-2)\times m\times n}$ and 
$\mathcal{L}_2\in\mathbb{R}^{(m-1)\times m\times n}$, where the tensor $\mathcal{L}_1$ has
the tridiagonal matrix 
\begin{equation}\label{regop1}
\mathcal{L}^{(1)}_1 = \frac{1}{4}\begin{bmatrix}
-1 & 2 & -1 \\
   & \ddots & \ddots & \ddots \\
   & & -1 & 2 & -1
\end{bmatrix}\in\mathbb{R}^{(m-2) \times m}
\end{equation}
as its first frontal slice, and the remaining frontal slices 
$\mathcal{L}_1^{(i)}\in\mathbb{R}^{(m-2)\times m}$, $i=2,3,\ldots,n$, are zero matrices. 
The first face of the tensor $\mathcal{L}_2$ is the bidiagonal matrix  
\begin{equation}\label{regop2}
\mathcal{L}^{(1)}_2 = \frac{1}{2}\begin{bmatrix}
1 & -1 \\
 & 1 & -1 \\
   & & \ddots & \ddots \\
   & &  & 1 & -1
\end{bmatrix}\in \mathbb{R}^{(m-1) \times m},
\end{equation} 
and the remaining faces $\mathcal{L}_2^{(i)}\in\mathbb{R}^{(m-1)\times m}$, 
$i=2,3,\dots,n$, are zero matrices. 

Our approach of handling these regularization operators is analogous to the technique used 
in \cite{HRY}. It can be applied to many other regularization operators as well. We use 
Algorithm \ref{Alg: 2} to compute the tQR factorization  
\[
\mathcal{L}*\mathcal{Q}_\ell=\mathcal{Q}_{\mathcal{L},\ell}*
\mathcal{R}_{\mathcal{L},\ell},
\]
where the tensor $\mathcal{Q}_{\mathcal{L},\ell}\in\mathbb{R}^{s\times\ell\times n}$ has 
$\ell$ orthonormal tensor columns and the tensor
$\mathcal{R}_{\mathcal{L},\ell}\in\mathbb{R}^{\ell\times\ell\times n}$ is f-upper 
triangular. In view of \eqref{lem: 2.1}, the minimization problem
\eqref{4.4} simplifies to
\begin{equation}\label{4.6}
\min_{\mathcal{\vec{Y}}\in \mathbb{R}^{\ell \times 1 \times n}}\{ \|\mathcal{\bar{H}}_\ell*\mathcal{\vec{Y}} - \vec{\mathit{e}}_1*\mathbf{z}_1\|^2_F + \mu^{-1}\|\mathcal{R}_{\mathcal{L},\ell}*\mathcal{\vec{Y}}\|^2_F\}.
\end{equation}

For the regularization operators \eqref{regop1} and \eqref{regop2}, as well as for many
other regularization operators $\mathcal{L}$, the tensor $\mathcal{R}_{\mathcal{L},\ell}$ 
is invertible and not very ill-conditioned. In this situation, we may form
\begin{equation}\label{4.7}
\mathcal{\vec{Z}} = \mathcal{R}_{\mathcal{L},\ell}*\mathcal{\vec{Y}}, \;\;\;\;\; 
\mathcal{\widetilde{H}}_\ell=\mathcal{\bar{H}}_\ell*\mathcal{R}_{\mathcal{L},\ell}^{-1},
\end{equation}
where $\mathcal{\widetilde{H}}_\ell$ is computed by solving $\ell$ systems of equations.
Substituting the above expressions into \eqref{4.6} yields
\begin{equation}\label{4.99}
\min_{\mathcal{\vec{Z}}\in \mathbb{R}^{\ell \times 1 \times n}}
\left\{\|\mathcal{\widetilde{H}}_\ell*\mathcal{\vec{Z}}-\vec{\mathit{e}}_1*\mathbf{z}_1\|^2_F
+ \mu^{-1}\|\mathcal{\vec{Z}}\|^2_F\right\}.
\end{equation}
The minimization problem \eqref{4.99} can be solved fairly stably by computing the solution of
\begin{equation}\label{4.9}
\min_{\mathcal{\vec{Z}}\in \mathbb{R}^{\ell \times 1 \times n}}
\bigg\| \begin{bmatrix} \mathcal{\widetilde{H}}_\ell \\
\mu^{-1/2}\mathcal{I} \end{bmatrix}* \mathcal{\vec{Z}} - \begin{bmatrix}
\vec{\mathit{e}}_1*\mathbf{z}_1\\
\mathcal{\vec{O}}
\end{bmatrix}\bigg\|_F
\end{equation}
using Algorithm \ref{Alg: 7} below. The solution of \eqref{4.9} can be expressed as
\begin{equation}\label{4.10}
\mathcal{\vec{Z}}_{\mu,\ell}=(\mathcal{\widetilde{H}}_\ell^T*\mathcal{\widetilde{H}}_\ell 
+\mu^{-1}\mathcal{I})^{-1}*\mathcal{\widetilde{H}}_\ell^T*
{\vec{\mathit{e}}}_1*\mathbf{z}_1,
\end{equation}
and the associated approximate solution of \eqref{LSQ} is given by 
\[
\mathcal{\vec{X}}_{\mu,\ell}=\mathcal{Q}_\ell *\mathcal{R}_{\mathcal{L},\ell}^{-1}*
(\mathcal{\widetilde{H}}_\ell^T*\mathcal{\widetilde{H}}_\ell+ 
\mu^{-1} \mathcal{I})^{-1}*\mathcal{\widetilde{H}}_\ell^T*\vec{\mathit{e}}_1*\mathbf{z}_1.
\]

\vspace{.3cm}
\begin{algorithm}[H]
\SetAlgoLined
\KwIn{$\mathcal{C}\in\mathbb{R}^{\ell\times m\times n}$, where its Fourier transform has 
nonsingular frontal slices; $\mathcal{\vec{D}}\in\mathbb{R}^{\ell\times 1\times n}$, 
$\mathcal{\vec{D}}\neq \mathcal{\vec{O}}$}
\KwOut{The solution $\mathcal{\vec{Y}}\in\mathbb{R}^{m\times 1\times n}$ of 
$\min_{\mathcal{\vec{Y}}\in\mathbb{R}^{m\times 1\times n}}
\|\mathcal{C*\vec{Y}}-\mathcal{\vec{D}}\|_F$}
$\mathcal{{C}} \leftarrow \mathtt{fft}(\mathcal{C},[\;],3)$\\
$\mathcal{{\vec{D}}} \leftarrow  \mathtt{fft}(\mathcal{\vec{D}},[\;],3)$\\
\For {$i =1$ \bf{to} $ n$ }{
$\mathcal{{\vec{Y}}}(:,:,i)=\mathcal{{C}}(:,:,i)\backslash\mathcal{{\vec{D}}}(:,:,i)$, 
where $\backslash$ denotes MATLAB's backslash operator
}
$\mathcal{\vec{Y}} \leftarrow \mathtt{ifft}(\mathcal{{\vec{Y}}},[\;],3)$ 
 \caption{Solution of a generic tensor least squares problem \cite{LU}}
 \label{Alg: 7}
\end{algorithm}
\vspace{.3cm}

We use the discrepancy principle \eqref{discr} to determine the regularization parameter 
$\mu>0$ and the required number of steps of the t-Arnoldi process as follows. Define the 
function
\begin{equation}\label{41.22}
\phi_\ell(\mu):=\|\mathcal{\widetilde{H}}_\ell*\mathcal{\vec{Z}}_{\mu,\ell}-
\vec{\mathit{e}}_1*\mathbf{z}_1\|_F^2,
\end{equation}
which is analogous to \eqref{11}. Substituting \eqref{4.10} into \eqref{41.22}, and using the 
identity 
\[
\mathcal{I} - \widetilde{\mathcal{H}}_\ell*(\widetilde{\mathcal{H}}_\ell^T*\widetilde{\mathcal{H}}_\ell+\mu^{-1}\mathcal{I})^{-1}*\widetilde{\mathcal{H}}_\ell^T =
(\mu\widetilde{\mathcal{H}}_\ell*\widetilde{\mathcal{H}}_\ell^T+\mathcal{I})^{-1},
\]
we obtain 
\be \label{phiprop}
\phi_\ell(\mu) = \left((\vec{\mathit{e}}_1*\mathbf{z}_1)^T*(\mu\mathcal{\widetilde{H}}_\ell*
\mathcal{\widetilde{H}}_\ell^T+\mathcal{I})^{-2}*\vec{\mathit{e}}_1*
\mathbf{z}_1\right)_{(:,:,1)}.
\ee

The following proposition shows that we can apply the discrepancy principle \eqref{discr} 
to the reduced problem to determine $\mu>0$, i.e., we require $\mu$ to be such that
\[
\|\mathcal{\widetilde{H}}_\ell*\mathcal{\vec{Z}}_{\mu,\ell}-
\vec{\mathit{e}}_1*\mathbf{z}_1\|_F = \eta\delta.
\]

\begin{prop}\label{prop32}
Let $\mu = \mu_\ell$ solve $\phi_\ell(\mu)=\eta^2\delta^2$ and let 
$\mathcal{\vec{Z}}_{\mu,\ell}$ solve \eqref{4.9}. Let $\mathcal{\vec{Y}}_{\mu,\ell}$ and
$\mathcal{\vec{Z}}_{\mu,\ell}$ be related by \eqref{4.7}. Then the associated approximate 
solution $\mathcal{\vec{X}}_{\mu,\ell} = \mathcal{Q}_\ell*\mathcal{\vec{Y}}_{\mu,\ell}$ of
\eqref{1} satisfies
\[
\|\mathcal{A}*\mathcal{\vec{X}}_{\mu,\ell} - \mathcal{\vec{B}}\|_F^2= 
\big((\vec{\mathit{e}}_1*\mathbf{z}_1)^T*(\mu\mathcal{\widetilde{H}}_\ell*
\mathcal{\widetilde{H}}_\ell^T + \mathcal{I})^{-2}*\vec{\mathit{e}}_1*
\mathbf{z}_1\big)_{(:,:,1)}.
\]
\end{prop}

\noindent
\textit{Proof}: Substituting 
$\mathcal{\vec{X}}_{\mu,\ell}=\mathcal{Q}_\ell*\mathcal{\vec{Y}}_{\mu,\ell}$ into 
\eqref{discr} and using the decomposition of \eqref{3.1}, as well as \eqref{QTB} and 
\eqref{lem: 2.1}, gives
\begin{equation*} 
\|\mathcal{A}*\mathcal{\vec{X}}_{\mu,\ell}-\mathcal{\vec{B}}\|_F^2=
\|\mathcal{Q}_{\ell+1}*\mathcal{\bar{H}}_\ell*\mathcal{\vec{Y}}_{\mu,\ell}-\mathcal{\vec{B}}\|_F^2 =
\|\mathcal{\bar{H}}_\ell*\mathcal{\vec{Y}}_{\mu,\ell}-\vec{e}_1*\mathbf{z}_1\|_F^2=
\|\mathcal{\widetilde{H}}_\ell*\mathcal{\vec{Z}}_{\mu,\ell}-\vec{e}_1*\mathbf{z}_1\|_F^2. ~~~\Box
\end{equation*}

It can be shown analogously as Proposition \ref{prop1} that the function $\phi_\ell(\mu)$ 
is decreasing and convex with $\phi_\ell(0)=\|\vec{\mathit{e}}_1*\mathbf{z}_1\|_F^2$. 
Therefore, Newton's method can be used for the solution of 
\begin{equation}\label{4.14}
\phi_\ell(\mu)-\eta^2\delta^2 = 0
\end{equation}
without safeguarding for any initial approximate solution $\mu_0\geq 0$ smaller than the 
solution of \eqref{4.14}. In particular, we may use $\mu_0 = 0$ when $\phi_\ell(\mu)$ and 
$\phi^\prime_\ell(\mu)$ are suitably defined at $\mu = 0$. Note that when the 
regularization parameter $\mu>0$ in \eqref{LSQ} is replaced by $1/\mu$, the analogue of 
the function $\phi_\ell$ obtained is not guaranteed to be convex. Then Newton's method has
to be safeguarded. An algorithm for Newton's method can be found in \cite{LU}.

We refer to the solution method for \eqref{4.2} described above as the tAT method. It is 
implemented by Algorithm \ref{Alg: 8} with $p=1$. It follows from Proposition \ref{prop1},
with $\phi$ replaced by $\phi_\ell$, that $\phi_\ell(\mu)$ is a decreasing function of 
$\mu$. A lower bound for $\phi_\ell(\mu)$ on the right-hand side of \eqref{minphik} can be
established similarly as in the proof of \cite[Proposition 3.6]{LU}. 
\begin{prop}\label{prop3.5}
Let $\phi_\ell(\mu)$ be given by \eqref{phiprop}. Then 
\begin{equation}\label{minphik}
\lim_{\mu \rightarrow \infty} \phi_\ell(\mu) = \Big(\mathbf{z}_1^T*\mathcal{U}(1,:,:)*\mathcal{D}*\mathcal{U}(1,:,:)^T*\mathbf{z}_1 \Big)_{(:,:,1)},
\end{equation}
where $\mathcal{D} \in \mathbb{R}^{(\ell+1) \times (\ell+1) \times n}$ is a tensor whose first frontal slice $\mathcal{D}^{(1)}$ has entry $1$ at the $(\ell+1, \ell+1)$st position, and the remaining frontal slices $\mathcal{D}^{(i)}$, $i=2, \dots, n$, are zero matrices. The tensor $\mathcal{U} \in \mathbb{R}^{(\ell+1) \times (\ell+1) \times n}$ is the left singular tensor of $\mathcal{\widetilde{H}}_\ell$.
\end{prop}

The values
\[
\ell\rightarrow\lim_{\mu\rightarrow\infty}\phi_\ell(\mu)
\]
typically decrease quite rapidly as $\ell$ increases, because making $\ell$ larger 
increases the dimension of the subspace over which the least squares problem \eqref{4.2}
is minimized. Therefore, generally, only a fairly small number of steps of Algorithm 
\ref{Alg: 8} are required to satisfy \eqref{4.14} for some $0<\mu<\infty$.

\subsubsection{tAT methods for the solution of \eqref{multrhs}}\label{sec3.12}
This subsection generalizes the solution methods of Subsection \ref{sec3.11} to the 
solution of least squares problems of the form \eqref{multrhs}. The methods of this 
subsection can be applied to color image and video restorations. Several matrix-based 
methods for the solution of these restoration problems have recently been described by 
Beik et al. \cite{BJNR,BNR} and El Guide et al. \cite{GIJB,GIJS}. 

We present two algorithms for the solution of \eqref{multrhs}. They both consider 
\eqref{multrhs} as $p$ separate Tikhonov minimization problems 
\begin{equation}\label{33n}
\min_{\mathcal{\vec{X}}_j \in \mathbb{R}^{m\times 1\times n}}
\{\|\mathcal{A}*\mathcal{\vec{X}}_j-\mathcal{\vec{B}}_j\|_F^2+
\frac{1}{\mu}\|\mathcal{L}*\mathcal{\vec{X}}_j\|_F^2\}, \;\;\;\; j = 1,2,\dots,p,
\end{equation}
where $\mathcal{\vec{B}}_1,\mathcal{\vec{B}}_2,\dots,\mathcal{\vec{B}}_p$ are tensor
columns of the data tensor $\mathcal{B}$ in \eqref{multrhs}. Both algorithms are based on
the t-Arnoldi process and the tAT method described in Subsection \ref{sec3.11}. 

Let $\mathcal{\vec{B}}_{j,{\rm true}}$ denote the unknown error-free tensor (slice) 
associated with the available error-contaminated tensor (slice) $\mathcal{\vec{B}}_j$, and
assume that bounds $\delta_j$ for the norm of the errors
\[
\mathcal{\vec{E}}_j:=\mathcal{\vec{B}}_j-\mathcal{\vec{B}}_{j,{\rm true}},\;\;\; 
j=1,2\dots,p,
\]
are available or can be estimated, i.e., 
\begin{equation}\label{deltaj}
\|\mathcal{\vec{E}}_j\|_F\leq\delta_j,\quad j=1,2,\ldots,p,
\end{equation}
cf. \eqref{unava} and \eqref{errbd}. Algorithm \ref{Alg: 8} solves each one of the $p$ 
least squares problems \eqref{33n} independently. 

\vspace{.3cm}
\begin{algorithm}[H]
\SetAlgoLined
\KwIn{$\mathcal{A}$, $p$, $\mathcal{\vec{B}}_1,\mathcal{\vec{B}}_2,\dots,
\mathcal{\vec{B}}_p$, $\delta_1,\delta_2,\dots,\delta_p$, $\mathcal{L}$, $\eta>1$, 
$\ell_\text{init}= 2$}
\For{$j=1,2,\dots,p$}{
$\ell\leftarrow\ell_\text{init}$, $[\mathcal{\vec{Q}}_1,\mathbf{z}_1]\leftarrow 
\mathtt{Normalize}(\mathcal{\vec{B}}_j)$.\\
Compute $\mathcal{Q}_\ell,\mathcal{Q}_{\ell+1}$, and $\mathcal{\bar{H}}_\ell$ by Algorithm
\ref{Alg: tArn}\\
Construct $\mathcal{R}_{\mathcal{L},\ell}$ by computing the tQR factorization of 
$\mathcal{L}*\mathcal{Q}_\ell$ using Algorithm \ref{Alg: 2}\\
Compute $\mathcal{\widetilde{H}}_\ell\leftarrow 
\mathcal{\bar{H}}_\ell*\mathcal{R}_{\mathcal{L},\ell}^{-1}$ and let
$\vec{\mathit{e}}_1\leftarrow\mathcal{I}(:,1,:)$\\
Solve the minimization problem 
\begin{equation*}
\min_{\mathcal{\vec{Z}}\in\mathbb{R}^{\ell\times 1\times n}} 
\|\mathcal{\widetilde{H}}_\ell*\mathcal{\vec{Z}}-\vec{\mathit{e}}_1*\mathbf{z}_1 \|_F
\end{equation*}
for $\mathcal{\vec{Z}}_\ell$ by using Algorithm \ref{Alg: 7}\\

\While{$\|\mathcal{\widetilde{H}}_\ell*\mathcal{\vec{Z}}_\ell- 
\vec{\mathit{e}}_1*\mathbf{z}_1\|_F  \geq \eta \delta_j$}{
$\ell \leftarrow \ell+1$\\
$\mathtt{Go \; to \; step \; 3}$}
Determine the regularization parameter by the discrepancy principle, i.e., compute the 
zero $\mu_\ell>0$ of 
\[
\xi_\ell(\mu):=\|\mathcal{\widetilde{H}}_\ell*\mathcal{\vec{Z}}_{j,\mu_\ell}-
\vec{\mathit{e}}_1*\mathbf{z}_1\|_F^2-\eta^2 \delta_j^2
\]  
and the associated solution $\mathcal{\vec{Z}}_{j,\mu_\ell}$ of 
 \begin{equation*}
 \min_{\mathcal{\vec{Z}}\in \mathbb{R}^{\ell \times 1 \times n}}\bigg\| \begin{bmatrix}
\mathcal{\widetilde{H}}_\ell \\
\mu_\ell^{-1/2}\mathcal{I}
\end{bmatrix}* \mathcal{\vec{Z}} - \begin{bmatrix}
\vec{\mathit{e}}_1*\mathbf{z}_1\\
\mathcal{\vec{O}}
\end{bmatrix}\bigg\|_F
\end{equation*}
by  using Algorithm \ref{Alg: 7}\\
Compute $\mathcal{\vec{Y}}_{j,\mu_\ell}\leftarrow 
\mathcal{R}_{\mathcal{L},\ell}^{-1}*\mathcal{\vec{Z}}_{j,\mu_\ell}, \;\; 
\mathcal{\vec{X}}_{j,\mu_\ell} \leftarrow \mathcal{Q}_\ell*\mathcal{\vec{Y}}_{j,\mu_\ell}$
\\
}
\caption{The tAT$_p$ method for the solution of \eqref{multrhs} by solving the $p$ problems 
\eqref{33n} independently}
 \label{Alg: 8}
\end{algorithm} \vspace{.3cm}

Algorithm \ref{Alg: 9} generates a t-Krylov subspace 
$\mathbb{K}_\ell(\mathcal{A},\mathcal{\vec{B}}_1)$ of sufficiently large dimension $\ell$
to contain accurate enough approximate solutions of all the $p$ least squares problems
\eqref{33n}. Thus, we first solve the least squares problem \eqref{33n} for $j=1$ by 
Algorithm \ref{Alg: 9}, and then seek to solve the least squares problem \eqref{33n} for 
$j=2$ using the same t-Krylov subspace 
$\mathbb{K}_\ell(\mathcal{A},\mathcal{\vec{B}}_1)$. If the discrepancy principle cannot 
be satisfied, then the dimension $\ell$ of the t-Krylov subspace is increased until the 
discrepancy principle can be satisfied. Having solved this least squares problem, we 
proceed similarly to solve the problems \eqref{33n} for $j=3,4,\ldots,p$. The details are 
described by Algorithm \ref{Alg: 9}. The t-Arnoldi process is implemented with 
reorthogonalization when applied in Algorithm \ref{Alg: 9} to ensure that the quantities 
$\mathcal{Q}_{\ell +1}^T*\mathcal{\vec{B}}_j$ are evaluated with sufficient accuracy. When 
the required number of t-Arnoldi steps, $\ell$, for solving the least squares problem 
is large, it may be beneficial to restart Algorithm \ref{Alg: 9} with the tensor 
$\mathcal{\vec{B}}_j$. Restarting was not required in the computations reported in 
Section \ref{sec6}.

\vspace{.3cm}
\begin{algorithm}[H]
\SetAlgoLined
\KwIn{$\mathcal{A}$, $p$, $\mathcal{\vec{B}}_1,\mathcal{\vec{B}}_2,\dots,
\mathcal{\vec{B}}_p$, $\delta_1,\delta_2,\dots,\delta_p$, $\mathcal{L}$, $\eta>1$, 
$\ell_\text{init}= 2$}
$\ell\leftarrow\ell_\text{init}$, $[\mathcal{\vec{Q}}_1,\sim]\leftarrow 
\mathtt{Normalize}(\mathcal{\vec{B}}_1)$  by Algorithm \ref{Alg: 01}\\
Compute $\mathcal{Q}_\ell,\mathcal{Q}_{\ell+1}$ and $\mathcal{\bar{H}}_\ell$ by Algorithm
\ref{Alg: tArn} with reorthogonalization of the tensor columns of $\mathcal{Q}_\ell$ and 
$\mathcal{Q}_{\ell +1}$\\
Construct $\mathcal{R}_{\mathcal{L},\ell}$ by computing the tQR factorization of 
$\mathcal{L}*\mathcal{Q}_\ell$ by using Algorithm \ref{Alg: 2}\\
Compute $\mathcal{\widetilde{H}}_\ell\leftarrow
\mathcal{\bar{H}}_\ell*\mathcal{R}_{\mathcal{L}, \ell}^{-1}$\\
Solve the minimization problem 
\begin{equation*}
\min_{\mathcal{\vec{Z}}\in \mathbb{R}^{\ell \times 1 \times n}} 
\|\mathcal{\widetilde{H}}_\ell*\mathcal{\vec{Z}}-\mathcal{Q}_{\ell+1}^T*
\mathcal{\vec{B}}_1\|_F
\end{equation*}
for $\mathcal{\vec{Z}}_\ell$ by using Algorithm \ref{Alg: 7}

\While{$\|\mathcal{\mathcal{\widetilde{H}}}_\ell*\mathcal{\vec{Z}}_\ell-
\mathcal{Q}_{\ell +1}^T*\mathcal{\vec{B}}_1 \|_F\geq\eta\delta_1$}{
 $\ell \leftarrow \ell+1$\\
 $\mathtt{Go \; to \; step \; 2}$}

Determine the regularization parameter $\mu_\ell$ by the discrepancy principle, i.e., 
compute the zero $\mu_\ell>0$ of 
\[
\xi_\ell(\mu):=\|\mathcal{\mathcal{\widetilde{H}}}_\ell*\mathcal{\vec{Z}}_{1,\mu_\ell}-
\mathcal{Q}_{\ell +1}^T*\mathcal{\vec{B}}_1 \|_F^2-\eta^2 \delta_1^2
\]
Compute the associated solution $\mathcal{\vec{Z}}_{1,\mu_\ell}$ of 
\begin{equation*}
 \min_{\mathcal{\vec{Z}}_1\in \mathbb{R}^{\ell \times 1 \times n}}
 \left\|\begin{bmatrix}
\mathcal{\widetilde{H}}_\ell \\
\mu_\ell^{-1/2}\mathcal{I}
\end{bmatrix}* \mathcal{\vec{Z}}_1 - \begin{bmatrix}
 \mathcal{Q}_{\ell +1}^T* \mathcal{\vec{B}}_1\\
\mathcal{\vec{O}}
\end{bmatrix}\right\|_F
\end{equation*}
by using Algorithm \ref{Alg: 7}\\
Compute $\mathcal{\vec{Y}}_{1,\mu_\ell}\leftarrow\mathcal{R}_{\mathcal{L},\ell}^{-1}*
\mathcal{\vec{Z}}_{1,\mu_\ell}, \;\; \mathcal{\vec{X}}_{1,\mu_\ell}\leftarrow 
\mathcal{Q}_\ell*\mathcal{\vec{Y}}_{1,\mu_\ell}$\\

\For{$j = 2,\dots, p$}{
$[\mathcal{\vec{Q}}_1, \sim] \leftarrow \mathtt{Normalize}(\mathcal{\vec{B}}_j)$\\
\While{$\|\mathcal{\mathcal{\widetilde{H}}}_\ell*\mathcal{\vec{Z}}_\ell-
\mathcal{Q}_{\ell+1}^T*\mathcal{\vec{B}}_j\|_F\geq\eta\delta_j$}{
$\ell\leftarrow \ell+1$\\
\mbox{Repeat~steps~2-5 with the present tensors $\mathcal{\widetilde{H}}_\ell$, 
$\mathcal{Q}_{\ell+1}^T$, and $\mathcal{\vec{B}}_j$}}
\mbox{Repeat~step~10} with the present $\delta_j$ and the tensors 
$\mathcal{\widetilde{H}}_\ell$, $\mathcal{Q}_{\ell+1}^T$, and $\mathcal{\vec{B}}_j$ to 
compute $\mathcal{\vec{Z}}_{j,\mu_\ell}$\\
Compute $\mathcal{\vec{Y}}_{j,\mu_\ell}\leftarrow
\mathcal{R}_{\mathcal{L},\ell}^{-1}*\mathcal{\vec{Z}}_{j,\mu_\ell}, \;\; 
\mathcal{\vec{X}}_{j,\mu_\ell}\leftarrow\mathcal{Q}_\ell*\mathcal{\vec{Y}}_{j,\mu_\ell}$\\
}
\caption{The {\tt nested}$\_$tAT$_p$ method for the solution of \eqref{multrhs} by solving the $p$ problems 
\eqref{33n} using a nested t-Krylov subspace}\label{Alg: 9}
\end{algorithm} \vspace{.3cm}

\subsection{tGMRES methods for the solution of \eqref{GEM} and \eqref{GGEM}}
\label{sec3.2}
We first describes the t-product GMRES (tGMRES) method for the approximate solution of 
\eqref{GEM}. This method subsequently will be generalized to the solution of 
problems of the form \eqref{GGEM}. We remark that the tGMRES method is analogous to the 
(standard) GMRES method introduced by Saad and Schultz \cite{SS}. Regularizing properties 
of the (standard) GMRES method for the situation when $\mathcal{A}$ is a matrix are 
discussed in \cite{CLR,Ne}. 

Substituting $\mathcal{\vec{X}}=\mathcal{Q}_\ell*\mathcal{\vec{Y}}$ into the right-hand 
side of \eqref{GEM}, using \eqref{3.1} as well as \eqref{QTB} and 
\eqref{lem: 2.1}, gives the reduced minimization problem 
\[
\min_{\mathcal{\vec{Y}}\in\mathbb{R}^{\ell\times 1\times n}}
\|\mathcal{\bar{H}}_{\ell}*\mathcal{\vec{Y}} - \vec{e}_1*\mathbf{z}_1\|_F.
\]
We refer to this solution method for \eqref{GEM} as the tGMRES method. It is implemented
by Algorithm \ref{Alg: 10} with $p=1$. The number of t-Arnoldi steps required by the 
tGMRES method is determined by the discrepancy principle
\begin{equation}\label{gdiscr}
\|\mathcal{\bar{H}}_{\ell}*\mathcal{\vec{Y}}-\vec{e}_1*\mathbf{z}_1\|_F\leq\eta\delta
\end{equation} 
in Algorithm \ref{Alg: 10}, where $\eta>1$ is a user-specified constant that is 
independent of $\delta$; cf. \eqref{discr}. Thus, we terminate the tGMRES iterations as 
soon as an iterate $\mathcal{\vec{Y}}=\mathcal{\vec{Y}}_\ell$ that satisfies 
\eqref{gdiscr} has been found. Generally, only fairly few iterations are needed. 
Restarting tGMRES therefore typically is not required.

We turn to a tGMRES method for the solution of \eqref{GGEM}, which we refer to as the 
tGMRES$_p$ method. This method, implemented by Algorithm \ref{Alg: 10}, considers 
\eqref{GGEM} as $p$ separate minimization problems 
\begin{equation}\label{GMRE}
\|\mathcal{A}*\mathcal{\vec{X}}_{j,\ell}-\mathcal{\vec{B}}_j \|_F= 
\min_{\mathcal{\vec{X}}_j \in \mathbb{K}_\ell(\mathcal{A}, \mathcal{\vec{B}}_j)} 
\|\mathcal{A}*\mathcal{\vec{X}}_j-\mathcal{\vec{B}}_j\|_F, \;\;\; \ell = 1,2,\dots,\;\; 
j = 1, 2, \dots, p,
\end{equation}
where $\mathcal{\vec{B}}_1,\mathcal{\vec{B}}_2,\dots,\mathcal{\vec{B}}_p$ are tensor
columns of the data tensor $\mathcal{B}$ in \eqref{GGEM}. The input parameters $\delta_j$
for Algorithm \ref{Alg: 10} are defined by \eqref{deltaj}. The number of steps $\ell$ is
chosen large enough to satisfy the discrepancy principle.

\vspace{.3cm}
\begin{algorithm}[H]
\SetAlgoLined
\KwIn{$\mathcal{A}$, $p$, $\mathcal{\vec{B}}_1,\mathcal{\vec{B}}_2,\dots,
\mathcal{\vec{B}}_p$, $\delta_1,\delta_2,\dots \delta_p$, $\mathcal{L}$, $\eta>1$, 
$\ell_\text{init}= 2$}
\For{$j = 1,2,\dots, p$}{
$\ell \leftarrow \ell_\text{init}$, $[\mathcal{\vec{Q}}_1,\mathbf{z}_1]\leftarrow 
\mathtt{Normalize}(\mathcal{\vec{B}}_j)$\\
Compute $\mathcal{Q}_\ell,\mathcal{Q}_{\ell+1}$ and $\mathcal{\bar{H}}_\ell$ by Algorithm 
\ref{Alg: tArn}\\
Construct $\vec{\mathit{e}}_1\leftarrow\mathcal{I}(:,1,:)$\\
Solve the minimization problem
\begin{equation*}
\min_{\mathcal{\vec{Y}}_j\in \mathbb{R}^{\ell \times 1 \times n}} 
\|\mathcal{\bar{H}}_\ell*\mathcal{\vec{Y}}_j -  \vec{\mathit{e}}_1*\mathbf{z}_1 \|_F
\end{equation*}
for $\mathcal{\vec{Y}}_{j,\ell}$ by using Algorithm \ref{Alg: 7}\\
\While{$\|\mathcal{\bar{H}}_\ell*\mathcal{\vec{Y}}_{j,\ell}-
\vec{\mathit{e}}_1*\mathbf{z}_1\|_F\geq\eta\delta_j$}{
 $\ell\leftarrow\ell+1$\\
 $\mathtt{Go \; to \; step \; 3}$}
Compute $\mathcal{\vec{X}}_{j,\ell}\leftarrow\mathcal{Q}_\ell*
\mathcal{\vec{Y}}_{j,\ell}$\\
}
\caption{The tGMRES$_p$ method for the solution of \eqref{GGEM}}\label{Alg: 10}
\end{algorithm} \vspace{.3cm}

\section{Methods Based on the Generalized Global t-Arnoldi Process}\label{sec4}
This section discusses the computation of an approximate solution of the tensor Tikhonov 
regularization problem \eqref{multrhs} and the minimization problem \eqref{GGEM} with 
the aid of the T-global Arnoldi process recently described by El Guide et al. \cite{GIJS}. 
Application of a few, say $1\leq\ell\ll m$, steps of the T-global Arnoldi process to the
tensor $\mathcal{A}\in\mathbb{R}^{m\times m\times n}$, reduces this tensor to a small 
upper Hessenberg matrix $\bar{H}_\ell\in\mathbb{R}^{(\ell+1)\times\ell}$. We refer to
this process as the generalized global t-Arnoldi (GG-tA) process. It is implemented by 
Algorithm \ref{Alg: 11}. We assume that the number of 
steps, $\ell$, is small enough to avoid breakdown. Then application of the GG-tA process to 
$\mathcal{A}$ with initial tensor $\mathcal{B}$ yields the decomposition
\begin{equation}\label{decom}
\mathcal{A}*\mathbb{Q}_\ell = \mathbb{Q}_{\ell+1} \circledast \bar{H}_\ell, 
\end{equation}
where  
\begin{equation*}
\mathbb{Q}_j := [\mathcal{Q}_1,\mathcal{Q}_2,\dots,\mathcal{Q}_j]\in
\mathbb{R}^{m\times pj \times n},\;\; j\in\{\ell,\ell+1\},
\end{equation*}
and
\begin{equation}\label{spec}
\begin{array}{rcl}
\mathcal{A}*\mathbb{Q}_\ell&=&[\mathcal{A}*\mathcal{Q}_1,\mathcal{A}*\mathcal{Q}_2,\dots,
\mathcal{A}*\mathcal{Q}_\ell]\in\mathbb{R}^{m \times \ell p \times n},\\
\mathbb{Q}_{\ell+1}\circledast\bar{H}_\ell&=&[\mathbb{Q}_{\ell+1}\circledast \bar{H}_\ell(:,1), 
\mathbb{Q}_{\ell+1}\circledast \bar{H}_\ell(:,2),\dots,\mathbb{Q}_{\ell+1}\circledast
\bar{H}_\ell(:,\ell)]\in\mathbb{R}^{m \times \ell p\times n}.
\end{array}
\end{equation}
The tensors $\mathcal{Q}_j\in\mathbb{R}^{m\times p\times n}$, $j=1,2,\dots,\ell$,
generated by Algorithm \ref{Alg: 11} form an orthonormal tensor basis for the t-Krylov 
subspace $\mathbb{K}_\ell(\mathcal{A},\mathcal{B})$, which is analogous to the space
\eqref{tKry}, 
\be \label{gKry}
\mathbb{K}_\ell(\mathcal{A},\mathcal{B}) = \bigg\{ \mathcal{Z} \in \mathbb{R}^{m \times p \times n},~ \mathcal{Z} = \sum_{i=1}^\ell \alpha_i (\mathcal{A}^{(i-1)}*\mathcal{B}),~\alpha_i\in\mathbb{R}\bigg\}.
\ee
Additional property of the t-Krylov subspace \eqref{gKry} is summed up in the following proposition; see Trefethen and Bau \cite{TB} for the matrix case.
\begin{prop}
Any $\mathcal{Z} \in \mathbb{K}_\ell(\mathcal{A},\mathcal{B})$ is equal to $p(\mathcal{A})*\mathcal{B}$ for some polynomial $p$ of degree $\leq \ell-1$.
\end{prop}
\noindent
{\it Proof:} $\forall$ $\mathcal{Z} \in \mathbb{K}_\ell(\mathcal{A},\mathcal{B})$, 
\begin{equation*}
\mathcal{Z} = \alpha_0 \mathcal{B} + \alpha_1 \mathcal{A*B} + \dots + \alpha_j \mathcal{A}^j*\mathcal{B} = (\alpha_0 + \alpha_1 \mathcal{A} + \dots + \alpha_j \mathcal{A}^j)*\mathcal{B}, ~j\leq \ell-1
\end{equation*}
Following the definition of standard tensor function, see, e.g., \cite{Lund, MQW1, MQW2}, and letting 
\[
p(\mathcal{A}) = \alpha_0 + \alpha_1 \mathcal{A} + \dots + \alpha_j \mathcal{A}^j = \sum_{j=0}^{\ell-1} \alpha_j \mathcal{A}^j, ~{\rm gives}~ \mathcal{Z} = p(\mathcal{A})*\mathcal{B}. ~~~~~\Box
\]

The upper Hessenberg matrix in \eqref{spec} is given by
\begin{center}
\begin{equation}\label{hess}
\bar{H}_\ell = \begin{bmatrix}
{h}_{11}&&&\dots & {h}_{1\ell}\\ 
{h}_{21} & h_{22} & \\
&{h}_{32} & {h}_{33} &&  \vdots \\
&& 	 \ddots  & \ddots & \\
&&& {h}_{\ell, \ell-1} & {h}_{\ell,\ell}\\
O &&&& {h}_{\ell+1,\ell}
\end{bmatrix} \in \mathbb{R}^{ (\ell+1)\times \ell}.
\end{equation}
\end{center}
The relation 
\begin{equation}\label{Bcom}
\mathcal{B}=\mathbb{Q}_{\ell+1}\circledast e_1\beta, \;\;\; e_1 = [1,0,\dots,0]^T
\end{equation}
is easily deduced from Algorithm \ref{Alg: 11}. 

\vspace{.3cm}
\begin{algorithm}[H]
\SetAlgoLined
\KwIn{$\mathcal{A}\in\mathbb{R}^{m\times m\times n}$, 
$\mathcal{B}\in\mathbb{R}^{m\times p\times n}$}
 Set $\beta\leftarrow\|\mathcal{B}\|_F$, $\mathcal{Q}_1\leftarrow 
 \frac{1}{\beta}\mathcal{B}$\\ 
\For {$j=1,2,\dots,\ell$}{
$\mathcal{W}\leftarrow\mathcal{A}*\mathcal{Q}_j$\\
 \For {$i=1,2,\dots,j$}{
$h_{ij}\leftarrow\langle \mathcal{Q}_i,\mathcal{W}\rangle$\\
$\mathcal{W}\leftarrow\mathcal{W}-h_{ij}\mathcal{Q}_i$
}
$h_{j+1, j}\leftarrow\|\mathcal{W}\|_F$, $\mathtt{if}~h_{j+1, j}=0$ $\mathtt{stop}$; 
$\mathtt{else}$\\
$\mathcal{Q}_{j+1}\leftarrow\mathcal{W}/h_{j+1,j}$
 }
 \caption{The generalized global t-Arnoldi (GG-tA) process \cite{GIJS}}
 \label{Alg: 11}
\end{algorithm}\vspace{.3cm}

Differently from the t-Arnoldi process, the GG-tA process uses the data tensor 
$\mathcal{B}\in \mathbb{R}^{m \times p \times n}$, $p>1$, and only requires transformation
to and from the Fourier domain in step 3. Each transformation of $\mathcal{A}$ and 
$\mathcal{Q}_j$ to and from the Fourier domain in step 3 costs $\mathcal{O}(m^2n\log(n))$ 
and $\mathcal{O}(mpn\log(n))$ flops, respectively. This step computes $\ell$ matrix-matrix
product of the frontal slices $\mathcal{\widehat{A}}^{(i)}$ and 
$\mathcal{\widehat{Q}}_j^{(i)}$, $i=1,2,\dots,n$, for $\mathcal{O}(\ell m^2 p)$ flops 
each. Hence for $n$ frontal slices, the cost of implementing step 3 in the Fourier domain 
is $\mathcal{O}(\ell m^2 p n)$ flops. The orthogonalization steps $4$-$7$ demands 
$\mathcal{O}( \ell^2 m np)$ flops. Hence, the GG-tA process has a complexity of 
$\mathcal{O}((\ell m^2 + \ell^2 m)np)$ flops in the Fourier domain. 
This cost is the same when the t-Arnoldi and G-tA processes are applied to separately solve
the $p$ minimization problems \eqref{33n}, since solving each one of the $p$ 
minimization problems independently costs $\mathcal{O}((\ell m^2 + \ell^2 m)n)$ flops in 
the Fourier domain.

We use the decomposition \eqref{decom} to determine an approximate solution of the 
Tikhonov minimization problem \eqref{multrhs} in Subsection \ref{sec4.1}, and of the 
minimization problem \eqref{GGEM} in Subsection \ref{sec4.2}.

\subsection{The GG-tAT method for the solution of \eqref{multrhs}}\label{sec4.1}
This subsection describes a modification of the T-global Arnoldi-Tikhonov 
regularization method recently presented by El Guide et al. \cite{GIJS} for the 
approximate solution of \eqref{multrhs} with $\mathcal{L=I}$ to allow a general third 
order tensor regularization operator $\mathcal{L \neq I}$. This modification requires Algorithm 
\ref{Alg: tggqr}. We refer to this modification of the method by El Guide et al. 
\cite{GIJS} as the generalized global tAT (GG-tAT) method. This method is based on first 
reducing $\mathcal{A}$ in \eqref{multrhs} to an upper Hessenberg matrix by carrying out a 
few, say $\ell$, steps of the GG-tA process, which is described by Algorithm 
\ref{Alg: 11}. Differently from the approach of El Guide et al. \cite{GIJS}, who apply a
restarted GG-tA process, determine the regularization parameter by the GCV, and use a 
stopping criterion based on the residual Frobenius norm and a prespecified tolerance that is 
independent of the error in the data tensor, we use the discrepancy principle to determine
the regularization parameter and the number of iterations required by the GG-tA process. 
Then the implementation of the GG-tA process does not required restarts since only a small 
number of iterations are needed.

We compute an approximate solution of \eqref{multrhs} analogously as described in 
Subsection \ref{sec3.11}. Thus, letting $\mathcal{X}=\mathbb{Q}_\ell\circledast y$, and 
using \eqref{decom} and \eqref{Bcom}, the minimization problem \eqref{multrhs} reduces to
\begin{equation}\label{N1}
\min_{y\in\mathbb{R}^\ell}\{\|\mathbb{Q}_{\ell+1}\circledast\bar{H}_\ell\circledast y-
\mathbb{Q}_{\ell+1}\circledast e_1\beta\|^2_F+\mu^{-1}\|\mathcal{L}*\mathbb{Q}_\ell 
\circledast y\|^2_F\},
\end{equation}
where $\beta=\|\mathcal{B}\|_F$. Algorithm \ref{Alg: tggqr} yields the 
GG-tQR factorization 
\begin{equation}
\mathcal{L}*\mathbb{Q}_\ell=\mathbb{Q}_{\mathcal{L},\ell}\circledast R_{\mathcal{L},\ell}
\in\mathbb{R}^{s\times \ell p\times n},
\label{N2}
\end{equation}
where $R_{\mathcal{L},\ell}\in\mathbb{R}^{\ell\times\ell}$ is an upper triangular matrix 
and $\mathbb{Q}_{\mathcal{L},\ell}\in\mathbb{R}^{s\times\ell p\times n}$ has $\ell$ 
orthonormal tensor columns. Substituting \eqref{N2} into \eqref{N1}, and using the 
left-hand side of \eqref{norm2F}, gives  
\begin{equation}
\min_{y\in\mathbb{R}^\ell}\{\|\bar{H}_\ell y-e_1\beta\|^2_2+
\mu^{-1}\|R_{\mathcal{L},\ell}y\|^2_2\}.
\label{N4}
\end{equation}
Typically, the matrix $R_{\mathcal{L},\ell}$ is nonsingular and not very ill-conditioned. 
Then we can express \eqref{N4} as a Tikhonov minimization problem in standard form,
\begin{equation}
\min_{z\in \mathbb{R}^\ell}\{\|\widetilde{H}_\ell z-e_1\beta\|^2_2+\mu^{-1}\|z\|^2_2\},
\label{A66}
\end{equation}
where
\begin{equation}\label{NN}
z:=R_{\mathcal{L},\ell}y, \;\;\; \widetilde{H}_\ell:=\bar{H}R^{-1}_{\mathcal{L},\ell}.
\end{equation}
Similarly as above, we compute $\widetilde{H}_\ell$ by solving $\ell$ linear systems of
equations. The minimization problem \eqref{A66} is analogous to \eqref{4.99}. Its 
solution, $z_{\mu,\ell}$, can be computed fairly stably by solving
\begin{equation}
\min_{z\in \mathbb{R}^\ell}\left\|\begin{bmatrix}
\widetilde{H}_\ell \\
\mu^{-1/2} I 
\end{bmatrix} z - 
\begin{bmatrix}
e_1 \beta\\
0
\end{bmatrix}\right\|_2.
\label{A6}
\end{equation}
The associated approximate solution of \eqref{multrhs} is given by 
\[
\mathcal{\vec{X}}_{\mu,\ell}=\mathbb{Q}_\ell\circledast R^{-1}_{\mathcal{L},\ell} 
z_{\mu,\ell}.
\]

We determine the regularization parameter $\mu$ by the discrepancy principle based on the 
Frobenius norm. This assumes knowledge of a bound 
\[
\|\mathcal{E}\|_F\leq\delta
\]
for the error $\mathcal{E}$ in $\mathcal{B}$. Thus, we choose $\mu>0$ so that the 
solution $z_{\mu,\ell}$ of \eqref{A6} satisfies
\[
\|\widetilde{H}_\ell z_{\mu,\ell}-e_1\beta\|_2=\eta\delta.
\]

Define the function 
\[
\psi_\ell(\mu):=\|\widetilde{H}_\ell z_{\mu,\ell}-e_1 \beta\|_2^2,
\]
where $z_{\mu,\ell}$ solves \eqref{A6}. Manipulations similar to those applied in 
Subsection \ref{sec3.11} show that $\psi_\ell(\mu)$ can be expressed as 
\be \label{prophi}
\psi_\ell(\mu) = \beta^2e_1^T(\mu \widetilde{H}_\ell\widetilde{H}_\ell^T +  I)^{-2}e_1.
\ee
It is readily verified that the 
function $\mu\rightarrow\psi_\ell(\mu)$ is decreasing and convex for $\mu\geq 0$ with 
$\psi_\ell(0)=\beta^2$. 

\begin{prop}\label{propp3.6}
Let $\psi_\ell(\mu)$ be given in \eqref{prophi}. Then 
\begin{equation}\label{propp36}
\lim_{\mu \rightarrow \infty} \psi_\ell(\mu) = \gamma \beta^2,
\end{equation}
where $\gamma>0$ is the square of the $(1,1)$ entry of the $(\ell+1)$st left singular vector 
of $\widetilde{H}_\ell$.
\end{prop}

The infimum of $\psi_\ell(\mu)$ on the right-hand side of \eqref{propp36} typically 
decreases quite rapidly as $\ell$, which is the dimension of the solution subspace, 
increases; see \cite{LU} for a proof \eqref{propp36}.

A similar reasoning as in Subsection \ref{sec3.1} suggests that it may be convenient to 
solve 
\begin{equation}
\psi_\ell(\mu)-\eta^2\delta^2 = 0
\label{A13}
\end{equation}
by Newton's method with initial approximate solution $\mu=0$.

We turn to a matrix analogue of Proposition \ref{prop32}.

\begin{prop}\label{prop3d}
Let $\mu_\ell$ solve \eqref{A13} and let $z_{\mu,\ell}$ be the associated solution of 
\eqref{A66} with $\mu=\mu_\ell$. Let $y_{\mu,\ell}$ and $z_{\mu,\ell}$ be related by 
\eqref{NN}. Then the approximate solution 
$\mathcal{X}_{\mu,\ell}=\mathbb{Q}_\ell\circledast y_{\mu,\ell}$ of \eqref{multrhs}
satisfies
\begin{equation}
\|\mathcal{A}*\mathcal{X}_{\mu,\ell}-\mathcal{B}\|_F^2= 
\beta^2e_1^T(\mu\widetilde{H}_\ell\widetilde{H}_\ell^T+I)^{-2}e_1.
\label{4.12}
\end{equation}
\end{prop}

\noindent
{\it Proof:} Substituting $\mathcal{X}_{\mu,\ell}=\mathbb{Q}_\ell\circledast y_{\mu,\ell}$
into \eqref{4.12}, using \eqref{decom} and \eqref{Bcom}, as well as left-hand side of 
\eqref{norm2F}, gives
\begin{equation*}
\|\mathcal{A}*\mathcal{X}_{\mu,\ell}-\mathcal{B}\|_F^2 
=\|\mathbb{Q}_{\ell+1} \circledast(\bar{H}_\ell\circledast y_{\mu,\ell}-e_1 \beta)\|_F^2 
=\|\bar{H}_\ell y_{\mu,\ell}-e_1\beta\|_2^2 =
\|\widetilde{H}_\ell z_{\mu,\ell}-e_1\beta\|_2^2. \;\; \Box
\end{equation*} 

We refer to the solution method described above as the GG-tAT method. It is implemented 
by Algorithm \ref{Alg: 12}. The method works with all lateral slices 
$\mathcal{\vec{B}}_j$, $j=1,2,\dots, p$, of $\mathcal{B}$ simultaneously.

\vspace{.3cm}
\begin{algorithm}[H]
\SetAlgoLined
\KwIn{$\mathcal{A}$, $\mathcal{B}$, $\delta$, $\mathcal{L}$, $\eta > 1$, 
$\ell_\text{init}=2$}
$\ell\leftarrow\ell_\text{init}$, $\beta\leftarrow\|\mathcal{B}\|_F$, 
$\mathcal{Q}_1\leftarrow\frac{1}{\beta}\mathcal{B}$\\
Compute $\mathbb{Q}_\ell$, $\mathbb{Q}_{\ell+1}$, and $\bar{H}_\ell$  by Algorithm 
\ref{Alg: 11}\\
Determine $R_{\mathcal{L},\ell}$ by computing the GG-tQR factorization of 
$\mathcal{L}*\mathbb{Q}_\ell$ using Algorithm \ref{Alg: tggqr}\\
Compute $\widetilde{H}_\ell\leftarrow\bar{H}_\ell R_{\mathcal{L},\ell}^{-1}$\\
Solve the minimization problem 
\begin{equation*}
\min_{z \in \mathbb{R}^\ell} \| \widetilde{H}_\ell z-e_1\beta\|_2
\end{equation*}
for $z_\ell$\\

\While{$\|\widetilde{H}_\ell z_\ell -  e_1 \beta\|_2  \geq \eta \delta$}{
 $\ell \leftarrow \ell+1$\\
$\mathtt{Go \; to \; step}$ 2
}
Determine the regularization parameter $\mu_\ell$ by the discrepancy principle, i.e., 
compute the zero $\mu_\ell>0$ of 
\[
\varphi_\ell(\mu):=\|\widetilde{H}_\ell z_{\mu,\ell}-e_1\beta\|_2^2-\eta^2\delta^2
\]
and the associated solution $z_{\mu,\ell}$ of
\[
\min_{z \in \mathbb{R}^\ell}\left\| \begin{bmatrix}
\widetilde{H}_\ell \\
\mu_\ell^{-1/2} I 
\end{bmatrix} z - \begin{bmatrix}
e_1 \beta\\
0
\end{bmatrix}\right\|_2
\] \\
Compute $y_{\mu,\ell}\leftarrow R^{-1}_{\mathcal{L},\ell}z_{\mu,\ell}, \;\; 
\mathcal{X}_{\mu,\ell}\leftarrow\mathbb{Q}_\ell \circledast y_{\mu,\ell}$ \\
 \caption{The GG-tAT method for the solution of \eqref{multrhs}}
 \label{Alg: 12}
\end{algorithm} \vspace{.3cm}

\subsection{The GG-tGMRES method for the approximate solution of \eqref{GGEM}}\label{sec4.2}
We describe the generalized global tGMRES (GG-tGMRES) method for the approximate solution
of \eqref{GGEM}. This method works with all lateral slices $\mathcal{\vec{B}}_j$, 
$j=1,2,\dots,p$, of $\mathcal{B}$ simultaneously. A closely related method, referred to as
the T-global GMRES method, recently has been described by El Guide et al. \cite{GIJS}. 
The latter method differs from the GG-tGMRES method in the following ways: it uses a 
restarted GG-tA process and a stopping criterion based on the residual Frobenius norm 
with a prespecified tolerance that is independent of the error in $\mathcal{B}$. The GG-tGMRES 
method uses the discrepancy principle to decide when to terminate the iterations. The 
number of iterations required by this method to satisfy the discrepancy principle 
typically is quite small. Restarting therefore generally is not required. 

Substituting $\mathcal{X} = \mathbb{Q}_\ell \circledast y$ into the right-hand side of 
\eqref{GGEM}, using \eqref{decom} and \eqref{Bcom}, as well as the 
left-hand side of \eqref{norm2F}, gives the reduced minimization problem 
\begin{equation}\label{tggmres}
\min_{y\in \mathbb{R}^\ell}\|\bar{H}_{\ell}y - \beta e_1\|_F.
\end{equation}
The GG-tGMRES method solves \eqref{tggmres} for a value of $\ell$ determined by the 
discrepancy principle and requires that a bound $\delta$ for $\|\mathcal{E}\|_F$ be known,
where $\mathcal{E}$ is the error in $\mathcal{B}$. This method is analogous to the tGMRES 
method described in Subsection \ref{sec3.2}. It is implemented by Algorithm \ref{Alg: 13}. 

\vspace{.3cm}
\begin{algorithm}[H]
\SetAlgoLined
\KwIn{$\mathcal{A}$, $\mathcal{B}$, $\delta$, $\mathcal{L}$, $\eta>1$, 
$\ell_\text{init}=2$}
\KwOut{Approximate solution $\mathcal{X}_\ell$ of \eqref{GGEM}}
$\ell\leftarrow\ell_\text{init}$, $\beta\leftarrow\|\mathcal{B}\|_F$, 
$\mathcal{Q}_1\leftarrow\frac{1}{\beta}\mathcal{B}$\\
Compute $\mathbb{Q}_\ell$, $\mathbb{Q}_{\ell+1}$, and $\bar{H}_\ell$ by Algorithm 
\ref{Alg: 11}\\
Solve the minimization problem 
\begin{equation*}
\min_{y\in\mathbb{R}^\ell}\|\bar{H}_\ell y - e_1\beta\|_2
\end{equation*}
for $y_\ell$\\
\While{$\|\bar{H}_\ell y_\ell-e_1\beta\|_2\geq\eta\delta$}{
 $\ell\leftarrow\ell+1$\\
$\mathtt{Go \; to \; step}$ 2
}

Compute $\mathcal{X}_\ell\leftarrow\mathbb{Q}_\ell\circledast y_\ell$ \\
\caption{The GG-tGMRES method for the solution of \eqref{GGEM}}
\label{Alg: 13}
\end{algorithm} \vspace{.3cm}

\section{Methods Based on the Global t-Arnoldi Process}\label{sec5}
This section discusses the computation of an approximate solution of the tensor Tikhonov 
regularization problems \eqref{LSQ} and \eqref{multrhs}, and of the minimization problems 
\eqref{GEM} and \eqref{GGEM}, with the aid of the global t-Arnoldi (G-tA) process. This  
process is readily implemented by taking $p=1$ in Algorithm \ref{Alg: 11}. We assume that 
$\ell$ is small enough to avoid breakdown. Algorithm \ref{Alg: 14} determines the G-tA 
decomposition
\[
\mathcal{A}*\mathcal{Q}_\ell = \mathcal{Q}_{\ell+1} \circledast \bar{\bar{H}}_\ell, 
\]
where
\begin{equation*}
\mathcal{Q}_j := [\mathcal{\vec{Q}}_1,\mathcal{\vec{Q}}_2,\dots,
\mathcal{\vec{Q}}_j]\in\mathbb{R}^{m\times j\times n},\;\; j\in\{\ell,\ell+1\}.
\end{equation*}

\vspace{.3cm}
\begin{algorithm}[H]
\SetAlgoLined
\KwIn{$\mathcal{A}\in\mathbb{R}^{m\times m\times n}$, 
$\mathcal{\vec{B}}\in\mathbb{R}^{m\times 1\times n}$}
 Set $\beta\leftarrow\|\mathcal{\vec{B}}\|_F$, 
 $\mathcal{\vec{Q}}_1\leftarrow\frac{1}{\beta}\mathcal{\vec{B}}$\\ 
 \For {$j=1,2,\dots,\ell$}{
 $\mathcal{\vec{W}}\leftarrow\mathcal{A}*\mathcal{\vec{Q}}_j$\\
 \For {$i=1,2,\dots,j$}{
$h_{ij}\leftarrow\langle\mathcal{\vec{Q}}_i,\mathcal{\vec{W}}\rangle $\\
$\mathcal{\vec{W}}\leftarrow\mathcal{\vec{W}}-h_{ij}\mathcal{\vec{Q}}_i$
}
$h_{j+1,j}\leftarrow\|\mathcal{\vec{W}}\|_F$, $\mathtt{if}~h_{j+1, j}=0$
$\mathtt{stop}$; $\mathtt{else}$\\
$\mathcal{\vec{Q}}_{j+1}\leftarrow\mathcal{\vec{W}}/h_{j+1, j}$
}
\caption{The global t-Arnoldi (G-tA) process}\label{Alg: 14}
\end{algorithm}

The expressions $\mathcal{A}*\mathcal{Q}_\ell$ and 
$\mathcal{Q}_{\ell+1}\circledast\bar{\bar{H}}_\ell$ are defined similarly to \eqref{spec},
and $\bar{\bar{H}} \in \mathbb{R}^{(\ell+1) \times \ell}$ has a form analogous to 
\eqref{hess}. The tensors $\mathcal{\vec{Q}}_j\in\mathbb{R}^{\ell\times 1\times n}$, 
$j=1,2,\dots, \ell$, generated by Algorithm \ref{Alg: 14} form an orthonormal tensor basis
for the t-Krylov subspace $\mathbb{K}_\ell(\mathcal{A},\mathcal{\vec{B}})$, where the
definition of {\rm t-span} is analogous to \eqref{gKry}. We use the G-tA process to 
determine an approximate solution of the Tikhonov minimization problems \eqref{multrhs} 
and \eqref{LSQ} in Section \ref{sec5.1}.

\subsection{The G-tAT method for the solution of \eqref{multrhs} and \eqref{LSQ}}\label{sec5.1}
We describe a solution method for \eqref{multrhs} that works with each lateral slice 
$\mathcal{\vec{B}}_j$, $j=1,2,\dots,p$, of the data tensor $\mathcal{B}$ independently. 
Thus, one solves \eqref{multrhs} by applying  the global t-product Arnoldi-Tikhonov (G-tAT) method  to the $p$ Tikhonov minimization problems \eqref{33n} separately.  We refer to this solution approach as the G-tAT$_p$ method. It is  implemented by Algorithm \ref{Alg: 15}.

The G-tAT method for the approximate solution of \eqref{LSQ} first reduces $\mathcal{A}$ in 
\eqref{LSQ} to an upper Hessenberg matrix by carrying out a few, say $\ell$, steps of 
the G-tA process described by Algorithm \ref{Alg: 14}. 
Let $\mathcal{\vec{X}}=\mathcal{Q}_\ell \circledast y$. Then following a similar approach as 
in Subsection \ref{sec4.1}, we reduce  \eqref{LSQ} to
\begin{equation}\label{J1}
\min_{y\in \mathbb{R}^\ell}\{ \|\mathcal{Q}_{\ell+1} \circledast \bar{\bar{H}}_\ell 
\circledast y-\mathcal{Q}_{\ell+1}\circledast e_1\beta\|^2_F+\mu^{-1}
\|\mathcal{L}*\mathcal{Q}_\ell \circledast y\|^2_F\}.
\end{equation}
Compute the G-tQR factorization of $\mathcal{L}*\mathcal{Q}_\ell$ by Algorithm 
\ref{Alg: tgqr} to obtain
\begin{equation}\label{tgqr}
\mathcal{L}*\mathcal{Q}_\ell=\mathcal{Q}_{\mathcal{L},\ell}\circledast
\bar{R}_{\mathcal{L},\ell},
\end{equation}
where the tensor $\mathcal{Q}_{\mathcal{L},\ell}\in\mathbb{R}^{s\times\ell\times n}$ has 
$\ell$ orthonormal tensor columns and the matrix 
$\bar{R}_{\mathcal{L},\ell}\in\mathbb{R}^{\ell\times\ell}$ is upper triangular.   

Substitute \eqref{tgqr} into \eqref{J1}, use the right-hand side of \eqref{norm2F}, and 
define 
\[
z := \bar{R}_{\mathcal{L},\ell} y, \;\;\; 
\breve{H}_\ell := \bar{\bar{H}} \bar{R}^{-1}_{\mathcal{L},\ell},
\]
where we assume that the matrix $\bar{R}_{\mathcal{L},\ell}$ is invertible and not very 
ill-conditioned. We obtain the Tikhonov minimization problem in standard form
\[
\min_{z\in\mathbb{R}^\ell}\{\|\breve{H}_\ell z-e_1\beta\|^2_2+\mu^{-1}\|z\|^2_2\}.
\]
This problem can be solved similarly as \eqref{A66}. We refer to this approach of solving 
\eqref{LSQ} as the G-tAT method. It is implemented by Algorithm \ref{Alg: 15} with $p=1$. 
The parameter $\delta_1$ is set to $\delta$ determined by \eqref{errbd}.  When applying Algorithm \ref{Alg: 15} to solve \eqref{multrhs}, the input parameters $\delta_1,\delta_2,\dots,\delta_p$ are determined by \eqref{deltaj}.

\vspace{.3cm}
\begin{algorithm}[H]
\SetAlgoLined
\KwIn{$\mathcal{A}$, $p$, $\mathcal{\vec{B}}_1,\mathcal{\vec{B}}_2,\dots,
\mathcal{\vec{B}}_p$, $\mathcal{L}$, $\delta_1,\delta_2,\dots,\delta_p$, $\eta > 1$, 
$\ell_\text{init}=2$}
\For{$j = 1,2,\dots, p$}{
$\ell\leftarrow\ell_\text{init}$, $\beta\leftarrow\|\mathcal{\vec{B}}_j\|_F$, 
$\mathcal{\vec{Q}}_1\leftarrow\frac{1}{\beta}\mathcal{\vec{B}}_j$\\
Compute $\mathcal{Q}_\ell$, $\mathcal{Q}_{\ell+1}$, and $\bar{\bar{H}}_\ell$ by Algorithm 
\ref{Alg: 14}\\
Determine $\bar{R}_{\mathcal{L},\ell}$ by computing the G-tQR factorization of 
$\mathcal{L}*\mathcal{Q}_\ell$ using Algorithm \ref{Alg: tgqr}\\
Compute $\breve{H}_\ell\leftarrow\bar{\bar{H}}_\ell\bar{R}_{\mathcal{L},\ell}^{-1}$\\
Solve the minimization problem  
\begin{equation*}
\min_{z\in\mathbb{R}^\ell}\|\breve{H}_\ell z-e_1\beta\|_2
\end{equation*}
for $z_\ell$\\
\While{$\|\breve{H}_\ell z_\ell-e_1\beta\|_2\geq\eta\delta_j$}{
$\ell\leftarrow\ell+1$\\
$\mathtt{Go \; to \; step}$ 3
}
Determine the regularization parameter $\mu_\ell>0$ by the discrepancy principle, i.e., by 
computing the zero $\mu_\ell$ of 
\[
\varphi_\ell(\mu):=\|\breve{H}_\ell z_{j,\mu_\ell}-e_1 \beta\|_2^2-\eta^2 \delta_j^2
\]
and the associated solution $z_{j,\mu_\ell}$ of
\[
\min_{z\in\mathbb{R}^\ell}\left\|\begin{bmatrix}
\breve{H}_\ell \\
\mu_\ell^{-1/2} I 
\end{bmatrix} z-\begin{bmatrix}
e_1\beta\\
0
\end{bmatrix}\right\|_2
\] \\
Compute: $y_{j,\mu_\ell}\leftarrow \bar{R}^{-1}_{\mathcal{L},\ell} z_{j,\mu_\ell}$, 
$\mathcal{\vec{X}}_{j,\mu_\ell}\leftarrow\mathcal{Q}_\ell\circledast y_{j,\mu_\ell}$ \\
}
\caption{The G-tAT$_p$ method for the solution of \eqref{multrhs}}\label{Alg: 15}
\end{algorithm}
\vspace{.3cm}

\subsection{The G-tGMRES method for the solution of \eqref{GEM} and 
\eqref{GGEM}}\label{sec5.2}
This subsection describes the global tGMRES (G-tGMRES) method for the approximate 
solution of \eqref{GEM} and \eqref{GGEM}. The G-tGMRES method uses the G-tA process 
described by Algorithm \ref{Alg: 13} and works with a data tensor slice 
$\mathcal{\vec{B}}$ in \eqref{GEM} and one lateral slice of the data tensor 
$\mathcal{B}$ at a time in \eqref{GGEM}. The G-tGMRES method is analogous to the 
GG-tGMRES method of the previous section. 

Substitute $\mathcal{\vec{X}}=\mathcal{Q}_\ell \circledast y$ into \eqref{GEM} and 
proceed similarly as described in Subsection \ref{sec4.2} to obtain the reduced 
minimization problem 
\[
\min_{y\in\mathbb{R}^\ell}\|\bar{\bar{H}}_{\ell}y-\beta e_1\|_2.
\]
We refer to the solution method so defined as the G-tGMRES method. It is implemented by
Algorithm \ref{Alg: 16} with $p=1$.  

We conclude this subsection by describing an algorithm for the approximate solution of 
\eqref{GGEM} based on the G-tGMRES method. This algorithm provides an alternative to the 
GG-tGMRES method of Subsection \ref{sec4.2}. It works with each lateral slice 
$\mathcal{\vec{B}}_j$, $j=1,2,\dots,p$, of the data tensor $\mathcal{B}$ independently. 
Thus, one solves the $p$ minimization problems \eqref{GMRE} separately by the tGMRES method. This approach is 
implemented by Algorithm \ref{Alg: 16} and will be referred to as the G-tGMRES$_p$ method. The parameters $\delta_1,\delta_2,\dots,\delta_p$ 
for the algorithm are determined by \eqref{deltaj}. 

\vspace{.3cm}
\begin{algorithm}[H]
\SetAlgoLined
\KwIn{$\mathcal{A}$, $p$, $\mathcal{\vec{B}}_1,\mathcal{\vec{B}}_2,\dots,
\mathcal{\vec{B}}_p$, $\mathcal{L}$, $\delta_1,\delta_2,\dots,\delta_p$, $\eta>1$,
$\ell_\text{init}= 2$}
\For{$j = 1,2,\dots,p$}{
$\ell\leftarrow\ell_\text{init}$, $\beta\leftarrow\|\mathcal{\vec{B}}_j\|_F$, 
$\mathcal{\vec{Q}}_1\leftarrow\frac{1}{\beta}\mathcal{\vec{B}}_j$\\
Compute $\mathcal{Q}_\ell$, $\mathcal{Q}_{\ell+1}$, and $\bar{\bar{H}}_\ell$ by Algorithm 
\ref{Alg: 14}\\
Solve the minimization problem 
\begin{equation*}
\min_{y_j\in\mathbb{R}^\ell}\|\bar{\bar{H}}_\ell y_j-e_1\beta\|_2
\end{equation*}
for $y_{j,\ell}$\\

\While{$\|\bar{\bar{H}}_\ell y_{j,\ell}-e_1\beta\|_2\geq\eta\delta_j$}{
$\ell\leftarrow\ell+1$\\
$\mathtt{Go \; to \; step}$ 3
}
Compute: $\mathcal{\vec{X}}_{j,\ell}\leftarrow\mathcal{Q}_\ell\circledast y_{j,\ell}$ \\
}
\caption{The G-tGMRES$_p$ method for the solution of \eqref{multrhs}}
 \label{Alg: 16}
\end{algorithm}\vspace{.3cm}
\noindent

\section{Numerical Examples}\label{sec6}

This section illustrates the performance of the methods described in the previous sections
when applied to the solution of several linear discrete ill-posed tensor problems. These 
methods are broadly categorized into two groups: those that involve flattening, i.e., 
reduce the tensor least squares problems \eqref{LSQ}, \eqref{GEM}, \eqref{multrhs} and \eqref{GGEM} to equivalent 
problems involving matrices and vectors, and those that preserve the tensor structure and 
do not involve flattening. We illustrate that it is generally beneficial to preserve the 
multidimensional tensor structure when solving linear discrete ill-posed tensor problems.

Applications to the restoration of (color) images and gray-scale videos are considered. Computed examples show that methods that preserve the natural spatial ordering yield the most accurate 
approximate solutions. In particular, tAT-type methods, such as tAT, tAT$_p$ and 
{\tt nested$\_$}tAT$_p$, give the best approximate solution in all computed examples except in Example \ref{E2}; Table \ref{Tab: 22}.
All computations were carried out in MATLAB 2019b on a Lenovo computer with an Intel Core i3 processor and 4 GB RAM running Windows 10. 

We use the discrepancy principle to determine the regularization parameter(s) and the 
number of steps of the iterative methods in all examples. The ``noise'' tensor 
$\mathcal{E}\in\mathbb{R}^{m\times p\times n}$, which simulates the error in the data 
tensor $\mathcal{B}=\mathcal{B}_{\text{true}}+\mathcal{E}$, is determined by 
its lateral slices $\mathcal{\vec{E}}_j$, $j=1,2,\dots,p$. The entries of these slices 
are normally distributed random numbers with zero mean and are scaled to correspond to a 
specified noise level $\widetilde{\delta}$. Thus,
\begin{equation}
\mathcal{\vec{E}}_j:=\widetilde{\delta} \frac{\mathcal{\vec{E}}_{0,j}}
{\|{\mathcal{\vec{E}}_{0,j}}\|_F}\|\mathcal{\vec{B}}_{\text{true},j}\|_F, \;\;\; 
j = 1,2,\dots,p,
\label{eq: e1}
\end{equation}
where the entries of the error tensors $\mathcal{\vec{E}}_{0,j}$ are $N(0,1)$. For problem
\eqref{1}, we have $p=1$. 

Let $\mathcal{\vec{X}}_{\rm method}$ be the computed approximate solution of \eqref{1} by 
a chosen method. The relative error 
\begin{equation*}
E_\text{method}=\frac{\|\mathcal{\vec{X}}_\text{method}-\mathcal{\vec{X}}_\text{true}\|_F}
{\|\mathcal{\vec{X}}_\text{true}\|_F} 
\end{equation*}
is used to determine the effectiveness of the proposed methods. The relative error for 
problems with a three-mode data tensor $\mathcal{B}$ is determined analogously. 

We let $\mathcal{A}\in\mathbb{R}^{256\times 256\times 256}$ in all computed examples 
unless otherwise stated. The condition number of the frontal slices of $\mathcal{A}$ are 
computed using the MATLAB command $\mathtt{cond}$. We set ${\tt tol}=10^{-12}$ in 
Algorithm \ref{Alg: 01}.

%==============================================================================
\begin{Ex}\label{E1}
This example compares Tikhonov regularization with the regularization tensor
$\mathcal{L}_2\in\mathbb{R}^{255 \times 256 \times 256}$, see \eqref{regop2}, as 
implemented by the tAT$_p$, {\tt nested}$\_$tAT$_p$, G-tAT$_p$ and GG-tAT
methods to the GMRES-type methods described by the tGMRES$_p$, G-tGMRES$_p$, and GG-tGMRES methods. 
Let the matrix 
\[
A_1 = {\tt gravity}(256,1,0,1,d), ~~~d=0.8, 
\]
be generated by the function ${\tt gravity}$ from the Hansen's Regularization Tools \cite{Haa} and define the
prolate matrix $A_2=\mathtt{gallery}('\mathtt{prolate}',256,\alpha)$ in MATLAB. We set 
$\alpha=0.46$. Then $A_2$ is a symmetric positive definite ill-conditioned Toeplitz 
matrix. The tensor $\mathcal{A}$ is defined by its frontal slices
\[
\mathcal{A}^{(i)} = A_1(i,1)A_2, ~~~ i=1,2,\dots,256.
\]

The exact data tensor $\mathcal{B}_{\rm true}\in\mathbb{R}^{256\times3\times 256}$ is 
given by $\mathcal{B}_{\rm true} = \mathcal{A}*\mathcal{X}_\text{true}$, where the exact 
solution $\mathcal{X}_{\rm true}\in\mathbb{R}^{256\times 3\times 256}$ has all entries 
equal to unity. The noise-contaminated right-hand side 
$\mathcal{B}\in\mathbb{R}^{256\times 3\times 256}$ is generated by 
$\mathcal{B}=\mathcal{B}_\text{true}+\mathcal{E}$, where the noise tensor 
$\mathcal{E}\in\mathbb{R}^{256\times 3\times 256}$ is determined according to 
\eqref{eq: e1}. The condition numbers of the slices $\mathcal{A}^{(i)}$ satisfy 
$\mathtt{cond}(\mathcal{A}^{(i)})\geq 1\cdot 10^{16}$ for all $i$. Thus, every slice is 
numerically singular. We take $\eta=1.15$ and determine the regularization parameter(s) 
for Tikhonov regularization by Newton's method. The computed regularization parameters 
and relative errors for different noise levels, as well as the number of iterations 
required to satisfy the discrepancy principle by each method, are displayed in Table 
\ref{Tab 2}. Here and below the table entry ``-'' indicates that the solution method carries out  different numbers of t-Arnoldi steps or computes different values of the regularization parameter for the different lateral slices of $\mathcal{B}$, or that no regularization 
parameter is required.

Table \ref{Tab 2} shows the GG-tAT and GG-tGMRES methods to be the fastest for both 
noise levels, but the tAT$_p$ and {\tt nested}$\_$tAT$_p$ methods, which do not involve 
flattening, yield approximate solutions of higher accuracy for both noise levels. 
The tAT$_p$ method determines the most accurate approximations of 
$\mathcal{X}_\text{true}$ and requires the most CPU time for both noise levels. The tGMRES$_p$ method yields the worst quality solution for both noise levels. In general, the 
quality of approximate solution is higher for Tikhonov regularization than for GMRES-type 
methods. This depends on the use of the regularization operator $\mathcal{L}_2$ by the former
methods.

\begin{table}[h!]
\begin{center}
\begin{tabular}{ccccccc} 
\hline
Noise level &Method &$\ell$ & $\mu_\ell$ & Relative error& CPU time (secs) \\
\hline
\multirow{7}{2em}{$10^{-3}$} & tAT$_p$ &-&- &2.09e-03 &12.53 \\ 
&{\tt nested}$\_$tAT$_p$ &3 &-& 2.23e-03 & 8.46  \\ 
&tGMRES$_p$&-&-  & 8.94e-01 & 7.67  \\ 
& G-tAT$_p$&-&-  & 6.20e-03 & 10.59 \\ 
&G-tGMRES$_p$&-&-  & 7.57e-03 & 7.16  \\ 
&GG-tAT&3& 7.13e-02 & 6.20e-03 & 5.50  \\ 
&GG-tGMRES &3&- & 7.57e-03 & 2.77 \\  \hline
\multirow{7}{2em}{$10^{-2}$} & tAT$_p$ &-&- &7.90e-03 &5.97 \\ 
&{\tt nested}$\_$tAT$_p$ &2 &-& 1.13e-02 & 4.82  \\ 
&tGMRES$_p$&-&-  & 4.71e+00 & 3.28  \\ 
&G-tAT$_p$&-&-  & 1.18e-02 & 4.76 \\ 
&G-tGMRES$_p$&-&-  & 2.37e-02 & 3.08  \\ 
&GG-tAT&2&3.09e-02  & 1.18e-02 & 2.31  \\ 
&GG-tGMRES &2&- & 2.37e-02 & 1.10 \\  \hline
\end{tabular}
\end{center}\vspace{-.5cm}
\caption{\small Results for Example \ref{E1}.}
\label{Tab 2}
\end{table}
\end{Ex}

\begin{Ex}\label{E2}
This example implements Example \ref{E1} analogously by taking $\mathcal{L=I}$, $d = 0.025$ to generate $A_1$, and determines the regularization parameter(s) by Newton's method with $\eta = 1.1$. The condition numbers of $\mathcal{A}^{(i)}$ are as described above. The relative errors for different noise levels and the CPU times are displayed in Table \ref{Tab 21}.

\begin{table}[h!]
\begin{center}
\begin{tabular}{ccccccc} 
\hline
Noise level &Method &$\ell$ & $\mu_\ell$ & Relative error& CPU time (secs) \\
\hline
\multirow{7}{2em}{$10^{-3}$} &  tAT$_p$ &-&- &6.69e-03 &14.25 \\ 
&{\tt nested}$\_$tAT$_p$ &3 &-& 4.35e-03 & 12.37  \\ 
&tGMRES$_p$&-&-  & 2.11e-02 & 8.12  \\ 
&G-tAT$_p$&-&-  & 5.65e-03 & 12.91 \\ 
&G-tGMRES$_p$&-&-  & 5.65e-03 & 7.47  \\ 
&GG-tAT&3& 3.28e-01 & 5.65e-03 & 6.54  \\ 
&GG-tGMRES &3&- & 5.65e-03 & 2.85 \\  \hline
\multirow{7}{2em}{$10^{-2}$} &  tAT$_p$ &-&- &4.10e-02 &6.47 \\ 
&{\tt nested}$\_$tAT$_p$ &2 &-& 2.59e-02 & 5.43  \\ 
&tGMRES$_p$&-&-  & 1.07e-01 & 3.31  \\ 
&G-tAT$_p$&-&-  & 2.46e-02 & 5.11 \\ 
&G-tGMRES$_p$&-&-  & 2.47e-02 & 3.01  \\ 
&GG-tAT&2&3.30e-02  & 2.46e-02 & 2.54  \\ 
&GG-tGMRES &2&- & 2.47e-02 & 1.16 \\  \hline
\end{tabular}
\end{center}\vspace{-.5cm}
\caption{\small Results for Example \ref{E2}.}
\label{Tab 21}
\end{table}

Table \ref{Tab 21} shows that the GG-tAT and GG-tGMRES methods that involve flattening
are the fastest for both noise levels. The {\tt nested}$\_$tAT$_p$ method, which does not
involve flattening and is based on nested t-Krylov subspaces, yields the most accurate 
approximate solutions. The G-tAT$_p$ and GG-tAT methods with Tikhonov regularization 
determine approximate solutions of almost the same quality as the GMRES-type methods 
implemented by the G-tGMRES$_p$ and GG-tGMRES methods for both noise levels. The 
tGMRES$_p$ method yields approximate solutions of least accuracy for both noise levels. 
For the solution methods that do not involve flattening (implemented by the tAT$_p$, 
{\tt nested}$\_$tAT$_p$, and tGMRES$_p$ methods), the quality of the computed approximate 
solutions is higher when Tikhonov regularization is applied.

We finally compare the tAT and G-tAT methods to the tGMRES and G-tGMRES methods. The
exact solution is the tensor column 
$\mathcal{\vec{X}}_{\rm true}\in\mathbb{R}^{256\times 1\times 256}$ with all entries equal
to unity. The noise-contaminated right-hand side 
$\mathcal{\vec{B}}\in\mathbb{R}^{256\times 1\times 256}$ is generated by 
$\mathcal{\vec{B}}=\mathcal{\vec{B}}_\text{true}+\mathcal{\vec{E}}$, where the noise 
tensor $\mathcal{\vec{E}}\in\mathbb{R}^{256\times 1\times 256}$ is generated as described
above. Table \ref{Tab: 22} shows the number of iterations required to satisfy the 
discrepancy principle by each method, the regularization parameters as well as the 
relative errors and CPU times for both noise levels.

\begin{table}[h!]
\begin{center}
\begin{tabular}{cccccccc} 
\hline 
\multicolumn{1}{c}{Noise level} & \multicolumn{1}{c}{Method}&\multicolumn{1}{c}{$\ell$} 
&\multicolumn{1}{c}{$\mu_\ell$} & \multicolumn{1}{c}{Relative error}& 
\multicolumn{1}{c}{CPU time (secs)}
\\ \hline 
\multirow{4}{4em}{$10^{-3}$} &tAT &3& 9.87e-01 & 8.40e-03 & 14.10   \\ 
&G-tAT &3& 7.25e-01 & 5.96e-03 & 13.72 \\ 
&tGMRES &3& - & 2.80e-02 & 3.10  \\ 
&G-tGMRES &3& - & 5.99e-03 & 2.90   \\ \hline
\multirow{4}{4em}{$10^{-2}$}&tAT &3& 5.54e-02 & 4.37e-02 & 3.67  \\ 
&G-tAT &2& 7.35e-02 & 2.46e-02 & 3.20  \\ 
&tGMRES &2&- & 1.45e-01 & 1.80 \\ 
&G-tGMRES &2& - & 2.56e-02 &1.00  \\ \hline
\end{tabular} 
\end{center} \vspace{-.5cm}
\caption{\small Results for Example \ref{E2}.}
\label{Tab: 22}
\end{table}

We see from Table \ref{Tab: 22} that the quality of approximate solution improves when
using Tikhonov regularization. The G-tGMRES and tGMRES methods are the fastest, but the 
tGMRES method yields approximate solutions of least quality for both noise levels. The 
G-tAT and G-tGMRES methods, which matricize the tensor equation, yield the most accurate 
solutions for both noise levels. This is the only one of our examples in which 
matricizing is beneficial for the quality of the computed solutions. In our experience
this situation is quite rare.
\end{Ex}

The remainder of this section discusses image and video restoration problems. We use the
bisection method to determine the regularization parameter over a chosen interval. The 
blurring operator $\mathcal{A}$ is constructed similarly as described in \cite{KKA} by 
using the function ${\tt blur}$ from \cite{Haa}. We determine the quality of restorations by each method using the relative error defined above, and Peak Signal-to-Noise Ratio (PSNR) defined by
\[
{\rm PSNR} = 10 {\rm log}_{10}\bigg( \frac{{\rm MAX}_{\mathcal{X}_{\rm true}}}{\sqrt{\rm MSE}} \bigg), ~{\rm where}~ {\rm MSE} = \frac{1}{mpn} \sum_{i=1}^m \sum_{j=1}^p \sum_{k=1}^n \big(\mathcal{X}_{\rm true} (i,j,k) - \mathcal{X}_{\rm method}(i,j,k)\big)^2
\]
denotes the Mean Square Error and ${\rm MAX}_{\mathcal{X}_{\rm true}}$ is the maximum of all the pixel values of the true image represented by $\mathcal{X}_{\rm true}\in \mathbb{R}^{m \times p \times n}$. The computation of MSE and ${\rm MAX}_{\mathcal{X}_{\rm true}}$ are carried out by using the MATLAB commands,
\[
{\rm MSE} = {\tt 1/(m*p*n)*sum(sum(sum((\mathcal{X}_{\rm true} - \mathcal{X}_{\rm method}).^2)))},
\]
\[
{\rm MAX}_{\mathcal{X}_{\rm true}} = {\tt max(max(max(\mathcal{X}_{\rm true})))}.
\]
For the problem \eqref{1}, discussed in Example \ref{Ex3}, we use $p=1$.

\begin{Ex}\label{Ex3} {(2D image restoration problem)}
This example illustrates the advantage of preserving the tensor structure when solving tensor linear discrete ill-posed problem. Specifically, we show that the tAT method which avoids flattening (matricization and vectorization) of the tensor equation \eqref{1} yields the best quality restorations for both noise levels independently of the regularization operators used. 

We discuss the performance of the tAT and G-tAT methods with the 
regularization tensors $\mathcal{L=I}$, and 
$\mathcal{L}=\mathcal{L}_1\in\mathbb{R}^{298\times 300\times 300}$ defined by
\eqref{regop2}, and compare these methods to the standard Arnoldi-Tikhonov
(AT) regularization method with regularization matrix $L=I$ described in
\cite{LR}, (standard) GMRES, tGMRES and G-tGMRES methods when applied to the restoration
of {\tt Telescope}\footnote{
\url{
https://github.com/jnagy1/IRtools/blob/master/Extra/test_data/HSTgray.jpg}} 
image of size $300 \times 300$ pixels that have been contaminated by blur and noise. The
AT and GMRES methods compute an approximate solution of the linear system of equations
\be \label{mtxe}
(A_1 \otimes A_2)x = b,
\ee
where $\otimes$ denotes the Kronecker product; the block matrix 
$A_1 \otimes A_2 \in \mathbb{R}^{300^2 \times 300^2}$ represents the blurring operator.
The right-hand side $b\in\mathbb{R}^{300^2}$ is the vectorized available blur- and 
noise-contaminated image $B\in\mathbb{R}^{300 \times 300}$.This vector is contaminated
by $e\in \mathbb{R}^{300^2}$, which represents (unknown) noise; it is a vectorization of
the noise matrix $E\in\mathbb{R}^{300 \times 300}$. We would like to determine an 
approximation of the ``true'' blur- and noise-free image 
$X_{\rm true}\in \mathbb{R}^{300 \times 300}$ or its vectorized form 
$x_{\rm true}\in\mathbb{R}^{300^2}$. The circulant matrix $A_1$ and Toeplitz matrix $A_2$ are generated with the MATLAB commands
\begin{equation}\label{mtxgen}
\begin{split}
\mathtt{z}_1 =
\mathtt{[exp(-([0:band-1].^2)/(2\sigma^2)),zeros(1,N-band)]},~~~
A_2 = \frac{1}{\sigma \sqrt{2\pi}}\mathtt{toeplitz(z_1)},\\
{\tt z_2} = {\tt [z_1(1) ~fliplr(z_1(end-length(z_1)+2:end))]}, ~~~ A_1 = \frac{1}{\sigma \sqrt{2\pi}}
{\tt toeplitz(z_1,z_2)},
\end{split}
\end{equation}
with $N = 300$, $\sigma=3$ and ${\tt band}=9$. By exploiting the circulant structure of 
$A_1 \otimes A_2$ and using the {\tt fold}, {\tt unfold}, and {\tt twist} operators, the 
2D deblurring problem \eqref{mtxe} can be formulated as the following 3D deblurring 
problem
\be \label{tmtx}
\mathcal{A*\vec{X}} = \mathcal{\vec{B}},
\ee
where $\mathcal{\vec{X}} = {\tt twist}(X)$, $\mathcal{\vec{B}} = {\tt twist}(B)$, and 
$\mathcal{\vec{E}} = {\tt twist}(E)$. The frontal slices
$\mathcal{A}^{(i)} \in\mathbb{R}^{300\times 300}$, $i = 1,2,\dots,300$, of
the blurring operator $\mathcal{A}\in\mathbb{R}^{300\times 300\times 300}$
are generated by folding the first block column of $A_1 \otimes A_2$, i.e.,
\be \label{constr}
\mathcal{A}^{(i)} = A_1(i,1)A_2, ~~~ i=1,2,\dots,300.
\ee
The computed condition numbers of $\mathcal{A}^{(i)}$ are
$\mathtt{cond}(\mathcal{A}^{(i)})=1.6\cdot 10^5$ for $i=1,2,\ldots,9$, and
$\mathtt{cond}(\mathcal{A}^{(i)})$ is ``infinite'' for $i\geq 10$. We let $\eta=1.1$ in
\eqref{discr} and determine the regularization parameter by the bisection method over the
interval $[10^{1},10^7]$.

The true image $\mathtt{Telescope}$ of size $300 \times 300$ is shown on the left-hand 
side of Figure \ref{Fig: 1}. For the matrix problem \eqref{mtxe}, this image is stored as 
a vector $x_{\rm true} \in \mathbb{R}^{300^2}$ and blurred by $A_1 \otimes A_2$, while for
the tensor problem \eqref{tmtx}, it is stored as
$\mathcal{\vec{X}}_{\text{true}}\in\mathbb{R}^{300\times 1\times 300}$
using the $\mathtt{twist}$ operator and blurred by the tensor
$\mathcal{A}$. The blurred and noisy image represented by $b$ is shown in Figure \ref{Fig: 1} (middle) using MATLAB {\tt reshape} command.

\begin{figure}[!htb]
\hspace{-1cm}
\minipage{0.43\textwidth}
\includegraphics[width=\linewidth]{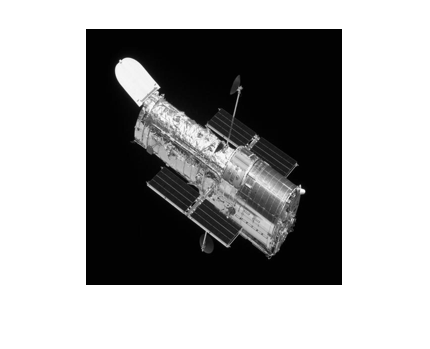} %\vspace{-.6cm}
\endminipage\hfill \hspace{-2cm}
\minipage{0.43\textwidth}
\includegraphics[width=\linewidth]{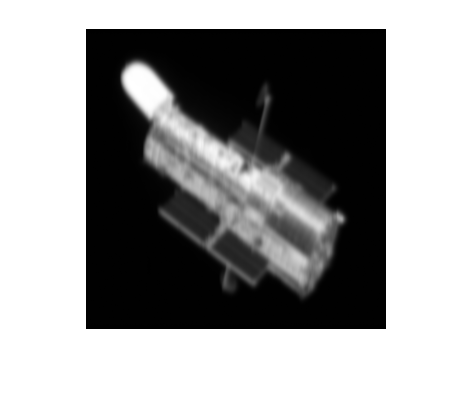}
\endminipage\hfill \hspace{-2cm}
\minipage{0.43\textwidth}
\includegraphics[width=\linewidth]{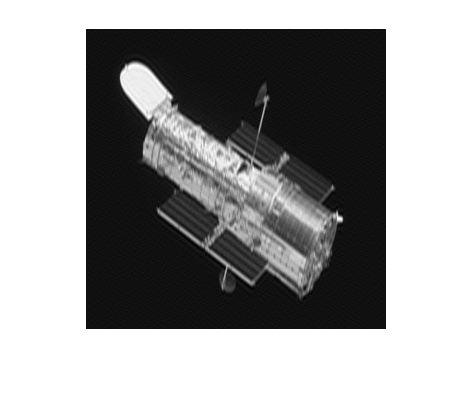}
\endminipage\hfill \hspace{-2cm} \vspace{-.9cm}
\caption{\small True image (left), blurred and noisy image (middle) with noise level $\widetilde{\delta} = 10^{-3}$, and restored image by the tAT (right) method after $8$ iterations.}\label{Fig: 1}
\end{figure}

The restored images determined by the tAT, G-tAT and tGMRES methods are accessed using the
$\mathtt{squeeze}$ operator and displayed in Figures \ref{Fig: 1} and \ref{Fig: 2} for the noise level 
$\widetilde{\delta}=10^{-3}$. Similarly, the restored image computed by the GMRES method 
is displayed in Figure \ref{Fig: 2} (middle) using MATLAB {\tt reshape} command. 

Table \ref{Tab: 2} shows the computed regularization parameters,relative errors and PNSR for the noise levels $10^{-2}$ and $10^{-3}$, as well as CPU times. As can be expected, the quality of the computed restorations improves when the noise level is smaller.  
The tGMRES method requires the least CPU time for $\widetilde{\delta}=10^{-3}$ and yields 
the worst restorations for both noise levels. Independent of the choice of $\mathcal{L}$,
Tikhonov regularization implemented by the tAT method determines restorations of the 
highest quality. The G-tAT and G-tGMRES methods, which involve flattening, demand the most 
CPU times and require the most iterations for both noise levels. The GMRES and G-tGMRES 
methods require the same number of iterations and yield the same quality restorations for 
both noise levels. Similar observations can be made for the AT and G-tAT methods when the
regularization operator is the identity matrix and identity tensor,
respectively.

\begin{figure}[!htb]
\hspace{-1cm}
\minipage{0.43\textwidth}
\includegraphics[width=\linewidth]{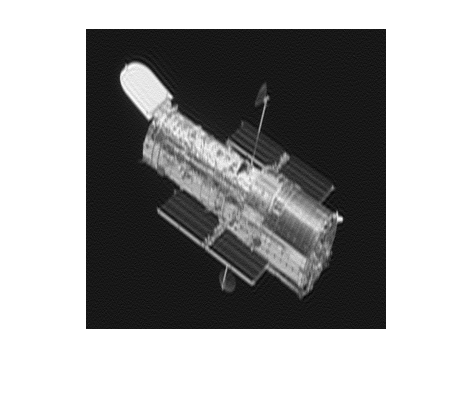} %\vspace{-.6cm}
\endminipage\hfill \hspace{-2cm}
\minipage{0.43\textwidth}
\includegraphics[width=\linewidth]{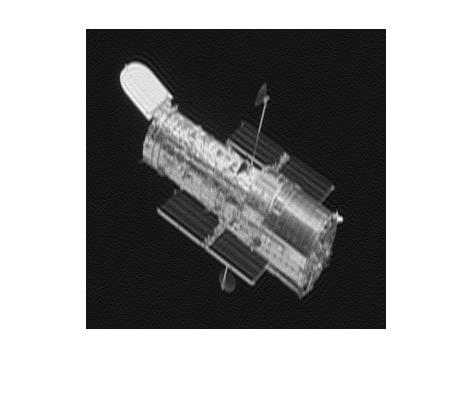}
\endminipage\hfill \hspace{-2cm}
\minipage{0.43\textwidth}
\includegraphics[width=\linewidth]{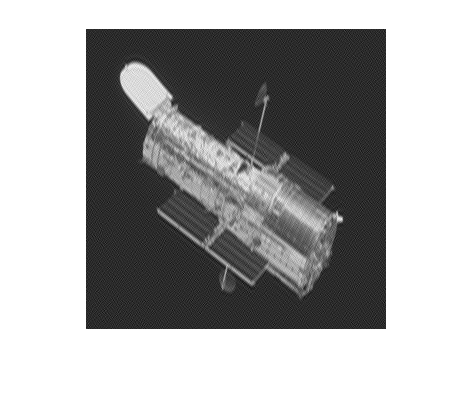}
\endminipage\hfill \hspace{-2cm} \vspace{-.9cm}
\caption{\small Restored images by the G-tAT (left), GMRES (middle), and tGMRES (right) 
methods after $51$, $51$, and $8$ iterations, respectively, for the noise level
$\widetilde{\delta}=10^{-3}$.}\label{Fig: 2}
\end{figure}

\begin{table}[h!]
\begin{center}
\begin{tabular}{cccccccc}
\hline
\multicolumn{1}{c}{$\mathcal{L}$} &\multicolumn{1}{c}{Noise level} &
\multicolumn{1}{c}{Method}&\multicolumn{1}{c}{$\ell$}
&\multicolumn{1}{c}{$\mu_\ell$} & \multicolumn{1}{c}{PSNR}&\multicolumn{1}{c}{Relative error}& \multicolumn{1}{c}{CPU time (secs)}
\\ \hline
\multirow{4}{1em}{$\mathcal{L}_1$} &\multirow{2}{2em}{$10^{-3}$} &tAT &8&
2.27e+04 &29.09& 1.19e-01 & 35.14 \\
&&G-tAT &51& 3.18e+04 &28.04& 1.34e-01 & 977.90 \\
\cmidrule(lr){2-8}
&\multirow{2}{2em}{$10^{-2}$}&tAT &3& 4.43e+01 &26.81& 1.53e-01 & 7.23 \\
&&G-tAT &12& 3.98e+02 &25.30& 1.84e-01 & 63.50 \\ \hline
\multirow{4}{1em}{$\mathcal{I}$} &\multirow{2}{2em}{$10^{-3}$} &tAT &8&
9.26e+04 &29.05& 1.19e-01 & 27.39 \\
&&G-tAT &51& 1.11e+05 &28.04& 1.34e-01 & 910.96 \\
\cmidrule(lr){2-8}
&\multirow{2}{2em}{$10^{-2}$}&tAT &3& 1.34e+03 &26.99& 1.51e-01 & 5.41 \\
&&G-tAT &12& 1.86e+03 &25.21& 1.86e-01 & 52.29 \\ \hline
&\multirow{1}{2em}{$10^{-3}$} &AT &51& 1.11e+05 &28.04& 1.34e-01 & 60.05 \\
&\multirow{1}{2em}{$10^{-2}$}&AT &12& 1.90e+03 &25.21& 1.86e-01 & 2.76 \\
\cmidrule(lr){2-8}
&\multirow{3}{2em}{$10^{-3}$} &GMRES &51& - &27.97& 1.35e-01 &59.91 \\
&&tGMRES &8& - &20.28& 2.03e-01 & 24.43 \\
&&G-tGMRES &51& - &27.97& 1.35e-01 & 898.32 \\
\cmidrule(lr){2-8}
&\multirow{3}{2em}{$10^{-2}$}&GMRES &12& - &24.94& 1.91e-01 &2.64\\
&&tGMRES &3&- &17.74& 4.39e-01 & 3.40 \\
&&G-tGMRES &12& - &24.94& 1.91e-01 &50.91 \\ \cmidrule(lr){2-8}
\end{tabular}
\end{center} \vspace{-.5cm}
\caption{\small Results for Example 6.3.}
\label{Tab: 2}
\end{table}
\end{Ex}

\begin{Ex}{(Color image restoration)}\label{E4} 
This example is concerned with the restoration of color images using the same 
regularization operators as in Example \ref{Ex3}. We seek to determine an approximate 
solution of the image deblurring problem
\be \label{bmtx}
(A_1 \otimes A_2)X = B,
\ee 
where the desired unavailable blur- and noise-free image 
$X_{\rm true}\in \mathbb{R}^{300^2 \times 3}$ is the matricized three-channeled image 
$\mathcal{X}_{\rm true}\in\mathbb{R}^{300\times 300\times 3}$. The right-hand side 
$B \in \mathbb{R}^{300^2 \times 3}$ in \eqref{bmtx} is generated by $B = (A_1 \otimes A_2)X_{\rm true} +  E $, where
the unknown noise in the matrix $B$ is represented by $E\in \mathbb{R}^{300^2 \times 3}$, 
which is the matricized ``noise'' tensor 
$\mathcal{E}\in\mathbb{R}^{300\times 300\times 3}$. The blurring matrices $A_1$ and $A_2$
are defined by \eqref{mtxgen}  in Example \ref{Ex3} with $N = 300$, $\sigma=3$ and ${\tt band}=12$. By the 
same reasoning as in Example \ref{Ex3}, we formulate \eqref{bmtx} as the 3D image 
deblurring problem
 \be \label{btmtx}
\mathcal{A*X} = \mathcal{B},
\ee 
where the blurring tensor $\mathcal{A}\in\mathbb{R}^{300\times 300\times 300}$ is constructed by \eqref{constr} in Example \ref{Ex3}.  The computed condition numbers of the frontal slices of $\mathcal{A}$ are 
$\mathtt{cond}(\mathcal{A}^{(i)})=7.6\cdot 10^8$ for $i=1,2,\dots,12$, and 
$\mathtt{cond}(\mathcal{A}^{(i)})$ is ``infinite'' for $i\geq 13$. We determine the regularization 
parameter(s) by the bisection method over the interval $[10^{-5},10^7]$.  The discrepancy principle is used with the parameter $\eta=1.1$. The (standard) global GMRES (G-GMRES) and (standard) global 
Arnoldi-Tikhonov (GAT) methods for \eqref{bmtx} are based on the global Arnoldi process applied by Huang et al. \cite{HRY}. We compare the performance of  these methods to the tAT$_p$, 
{\tt nested}$\_$tAT$_p$, G-tAT$_p$, GG-tAT, tGMRES$_p$, G-tGMRES$_p$, and GG-tGMRES 
methods for the solution of \eqref{btmtx}.  

The original (blur- and noise-free) 
{\tt flower}\footnote{\url{http://www.hlevkin.com/TestImages}} image shown on the 
left-hand side of Figure \ref{Fig: 3} is stored as a tensor 
$\mathcal{X}_\text{true}\in\mathbb{R}^{300\times 3\times 300}$. It is blurred using the 
tensor $\mathcal{A}$. Thus, 
$\mathcal{{B}}_{\text{true}}=\mathcal{A}*\mathcal{{X}}_{\text{true}}\in
\mathbb{R}^{300\times 3\times 300}$ represents the blurred but noise-free image 
associated with $\mathcal{X}_{\text{true}}$. The ``noise'' tensor 
$\mathcal{{E}}\in\mathbb{R}^{300\times 3\times 300}$ is generated as described by 
\eqref{eq: e1} with noise level $\widetilde{\delta} = 10^{-3}$ and added to 
$\mathcal{B}_{\text{true}}$ to obtain the blurred and noisy image $\mathcal{{B}}$ 
shown in Figure \ref{Fig: 3} (middle). The latter image is accessed by using the {\tt multi$\_$squeeze}
operator. 

\begin{figure}[!htb] 
\hspace{-1cm}
\minipage{0.42\textwidth}
\includegraphics[width=\linewidth]{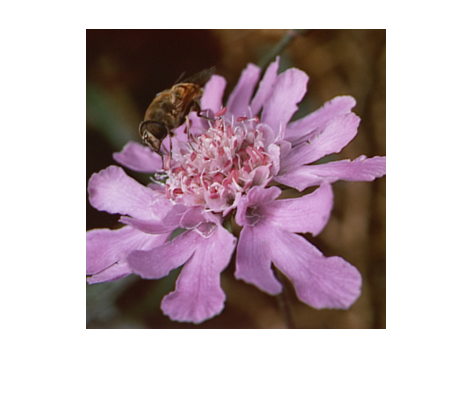} %\vspace{-.6cm} 
\endminipage\hfill \hspace{-2cm}
\minipage{0.42\textwidth}
\includegraphics[width=\linewidth]{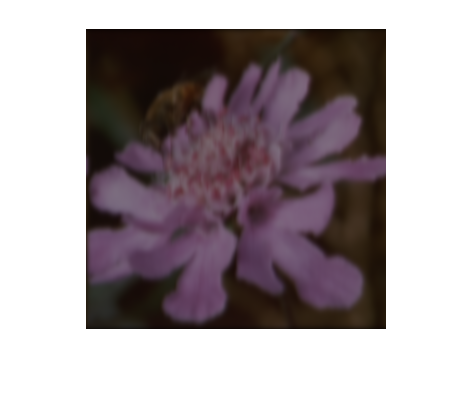} %\vspace{-.6cm} 
\endminipage\hfill \hspace{-2cm}
\minipage{0.42\textwidth}
\includegraphics[width=\linewidth]{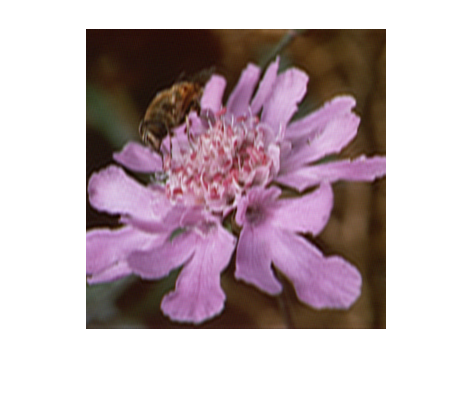} 
\endminipage\hfill \hspace{-2cm}  
\vspace{-.9cm}
\caption{\small True image (left),  blurred and noisy image (middle) with noise level $\widetilde{\delta} = 10^{-3}$, and restored image determined by {\tt nested}$\_$tAT$_p$ (right) after $9$ iterations.}\label{Fig: 3}
\end{figure} 

\begin{figure}[!htb] 
\hspace{-1cm}
\minipage{0.42\textwidth}
\includegraphics[width=\linewidth]{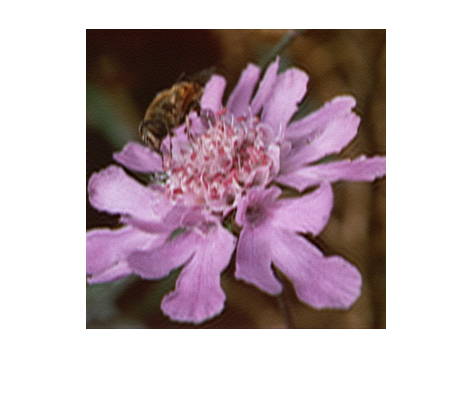} %\vspace{-.6cm} 
\endminipage\hfill \hspace{-2cm}
\minipage{0.42\textwidth}
\includegraphics[width=\linewidth]{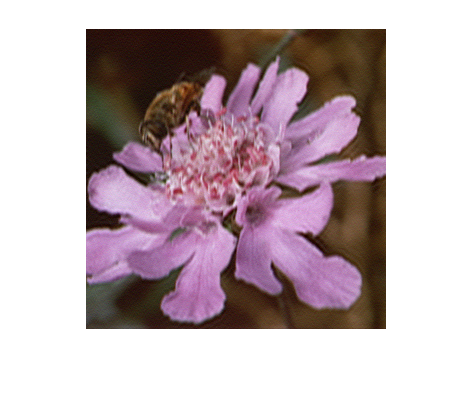} %\vspace{-.6cm} 
\endminipage\hfill \hspace{-2cm}
\minipage{0.42\textwidth}
\includegraphics[width=\linewidth]{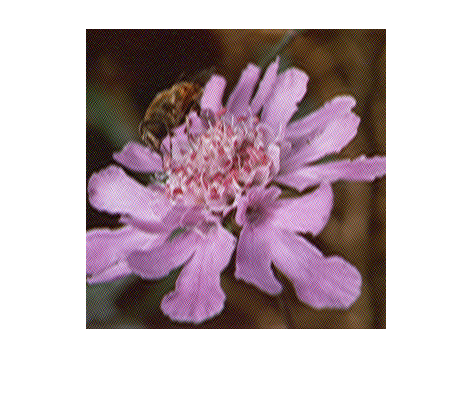} 
\endminipage\hfill \hspace{-2cm}  
\vspace{-.9cm}
\caption{\small Restored images determined by GG-tAT (left) after $32$ 
iterations, and G-GMRES (middle) after $32$ iterations and the tGMRES$_p$ (right), 
for the noise level $\widetilde{\delta} = 10^{-3}$.}\label{Fig: 4}
\end{figure}

The restored images determined by the {\tt nested}$\_$tAT$_p$, GG-tAT, G-GMRES, and tGMRES methods are
displayed in Figures \ref{Fig: 3} and \ref{Fig: 4}. Relative errors and PSNR as well as CPU times are shown in Table 
\ref{Tab: 4}. The tAT method gives restorations of the highest quality, followed by the 
{\tt nested}$\_$tAT$_p$ method. These methods do not involve flattening. Solution methods that 
involve flattening such as the G-tAT$_p$, GG-tAT, G-tGMRES$_p$, and GG-tGMRES methods 
require the most CPU time for both noise levels. The GAT, GG-tAT, G-GMRES, and GG-tGMRES 
methods require the same number of iterations, which are more than the number of 
iterations used by the {\tt nested}$\_$tAT$_p$ method for both noise levels. The tGMRES$_p$ 
method yields restorations of the worst quality for both noise levels. The GG-tGMRES 
method that works with the whole data tensor at a time yields the same quality 
restorations as the G-GMRES method for both noise levels. The same conclusion can be drawn
for the GG-tAT and GAT methods when the regularization operator is the identity tensor and 
the identity matrix, respectively. The quality of restorations by the G-tAT$_p$ and GG-tAT methods improves significantly with the use of the regularization operator $\mathcal{L}_1$ for both noise levels.

\begin{table}[h!]
\begin{center}
\begin{tabular}{cccccccc} 
\hline
$\mathcal{L}$&Noise level &Method &$\ell$ & $\mu_\ell$ &PSNR& Relative error& CPU time (secs) \\ \hline 
\multirow{8}{1em}{$\mathcal{L}_1$}&\multirow{4}{2em}{$10^{-3}$} & tAT$_p$ &-&- & 30.56 &5.85e-02 & 93.08 \\ 
&&{\tt nested}$\_$tAT$_p$ &9 &-& 30.56 & 5.86e-02 & 64.79  \\
&&G-tAT$_p$&-&-  & 29.47 & 6.64e-02 & 1187.47 \\ 
&&GG-tAT&32& 7.34e+03 & 29.43 & 6.67e-02 & 894.05  \\
\cmidrule(lr){2-8}
&\multirow{4}{2em}{$10^{-2}$} & tAT$_p$ &-&- & 27.20 &8.62e-02 &21.07 \\ 
&&{\tt nested}$\_$tAT$_p$ &4 &-& 25.90 & 1.00e-01 & 21.51  \\
&&G-tAT$_p$&-&-  & 25.20 & 1.09e-01 & 101.27 \\ 
&&GG-tAT&9&1.51e+02  & 25.22 & 1.08e-01 & 70.43  \\ \hline

\multirow{8}{1em}{$\mathcal{I}$}&\multirow{4}{2em}{$10^{-3}$} & tAT$_p$ &-&-& 30.67 &5.78e-02 & 75.66 \\ 
&&{\tt nested}$\_$tAT$_p$ &9 &-& 30.69 & 5.77e-02 & 57.02  \\
&&G-tAT$_p$&-&-  & 29.44 & 6.66e-02 & 1114.94 \\ 
&&GG-tAT&32& 3.64e+04 & 29.40 & 6.69e-02 & 430.23  \\
\cmidrule(lr){2-8}
&\multirow{4}{2em}{$10^{-2}$} & tAT$_p$ &-&- & 27.66 &8.18e-02 &16.64 \\ 
&&{\tt nested}$\_$tAT$_p$ &4 &-& 26.26 & 9.60e-02 & 22.88\\
&&G-tAT$_p$&-&-  & 24.96 & 1.12e-01 & 82.66\\ 
&&GG-tAT&9&1.15e+03  & 24.87 & 1.13e-01 & 33.77 \\ \hline

&\multirow{1}{2em}{$10^{-3}$} &GAT&32&3.63e+04  & 29.40 & 6.69e-02 & 96.86 \\
&\multirow{1}{2em}{$10^{-2}$} &GAT&9&1.15e+03  & 24.87 & 1.13e-01 & 6.22 \\ \cmidrule(lr){2-8}

&\multirow{4}{2em}{$10^{-3}$} &G-GMRES &32&- & 29.33 & 6.75e-02 & 98.89 \\  
&&tGMRES$_p$&-&-  & 19.23 & 2.16e-01 & 69.72  \\  
&&G-tGMRES$_p$&-&-  & 29.36 & 6.73e-02 & 1105.67 \\ 
&&GG-tGMRES &32&- & 29.32 & 6.75e-02 & 425.55 \\
\cmidrule(lr){2-8} 
&\multirow{4}{2em}{$10^{-2}$} &G-GMRES &9&- & 24.56 & 1.17e-01 & 5.21 \\ 
&&tGMRES$_p$&-&-  & 12.75 & 4.55e-01 & 10.45  \\ 
&&G-tGMRES$_p$&-&-  & 24.78 & 1.14e-01 & 82.35  \\
&&GG-tGMRES &9&- & 24.56 & 1.17e-01 & 33.07 \\  
\cmidrule(lr){2-8}
\end{tabular}
\end{center}\vspace{-.5cm}
\caption{\small Results for Example \ref{E4}.}
\label{Tab: 4}
\end{table}

\end{Ex}

%=================================example 4.5===================================
\begin{Ex}{(Video restoration)}\label{E5}
This example considers the restoration of the first six consecutive frames of the 
$\mathtt{Xylophone}$ video from MATLAB. Each video frame is in the $\mathtt{MP4}$ format 
and has $240\times 240$ pixels.

\begin{figure}[!htb] 
\hspace{-1cm}
\minipage{0.43\textwidth}
\includegraphics[width=\linewidth]{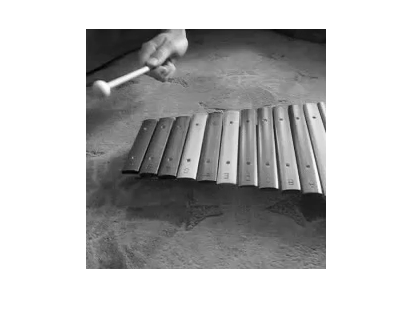} %\vspace{-.6cm} 
\endminipage\hfill \hspace{-2cm}
\minipage{0.43\textwidth}
\includegraphics[width=\linewidth]{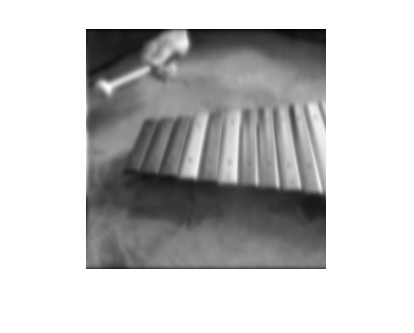} %\vspace{-.6cm} 
\endminipage\hfill \hspace{-2cm}
\minipage{0.43\textwidth}
\includegraphics[width=\linewidth]{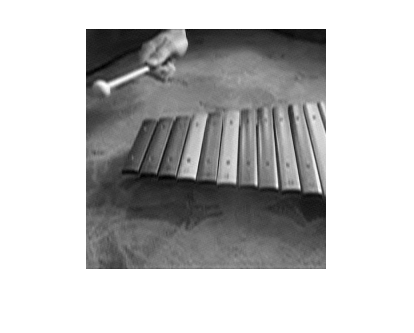} 
\endminipage\hfill \hspace{-2cm}  \vspace{-.9cm}
\caption{\small True image (left), blurred and noisy image with noise level 
$\widetilde{\delta}=10^{-3}$ (middle), and restored image determined by $8$ iterations
with the tAT$_p$ method (right).}\label{Fig: 5}
\end{figure}

The first six blur- and noise-free frames are stored as a tensor 
$\mathcal{X}_\text{true}\in\mathbb{R}^{240\times 6\times 240}$ using the 
$\mathtt{multi}\_\mathtt{twist}$ operator. They are blurred by the tensor
$\mathcal{A}\in\mathbb{R}^{240\times 240\times 240}$, which is generated similarly as
in Example \ref{Ex3} with its frontal slices determined by 
\[
\mathcal{A}^{(i)} = A_2(i,1)A_2, ~~~ i = 1,2,\dots, n, ~~~ N = 240,~~ \sigma=2.5~~ {\rm and}~~ {\tt band}=12.
\]
The condition numbers of the 
frontal slices of $\mathcal{A}$ are $\mathtt{cond}(\mathcal{A}^{(i)})=1.4\cdot 10^7$ for 
$i=1,2,\dots,12$. The condition numbers of the remaining frontal slices are ``infinite''. 

We use the regularization operator
$\mathcal{L}=\mathcal{L}_2\in\mathbb{R}^{239\times 240\times 240}$ and
determine the  regularization parameter(s) by the bisection method over the interval 
$[10^{-5},10^7]$ using the discrepancy principle with $\eta=1.1$. The blurred and noisy 
frames are generated by 
$\mathcal{B} = \mathcal{A}*\mathcal{X}_\text{true} + \mathcal{E}  \in 
\mathbb{R}^{240 \times 6 \times 240}$  with the ``noise'' tensor
$\mathcal{E}  \in \mathbb{R}^{240 \times 6 \times 240}$ defined by \eqref{eq: e1}. 

The true third frame is displayed in Figure \ref{Fig: 5} (left), and the blurred and 
noisy third frame is shown in Figure \ref{Fig: 5} (middle) using the $\mathtt{squeeze}$ 
operator. Similarly, the restored images of the third frame determined by the G-tAT$_p$,
{\tt nested}$\_$tAT$_p$, G-tGMRES, and tGMRES methods are shown in Figures \ref{Fig: 5} 
and \ref{Fig: 6}. 

\begin{figure}[!htb] 
\hspace{-1cm}
\minipage{0.43\textwidth}
\includegraphics[width=\linewidth]{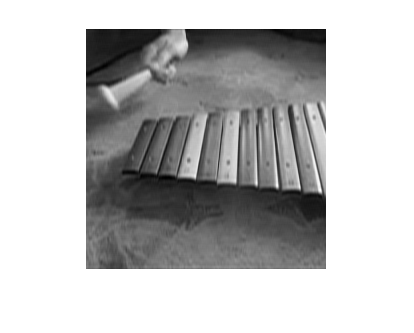} %\vspace{-.6cm} 
\endminipage\hfill \hspace{-2cm}
\minipage{0.43\textwidth}
\includegraphics[width=\linewidth]{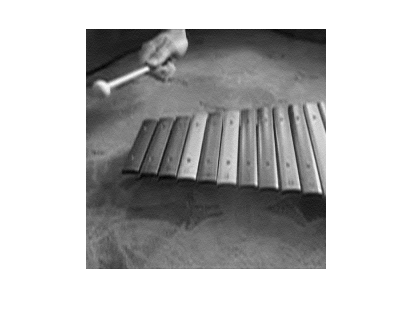} %\vspace{-.6cm} 
\endminipage\hfill \hspace{-2cm}
\minipage{0.43\textwidth}
\includegraphics[width=\linewidth]{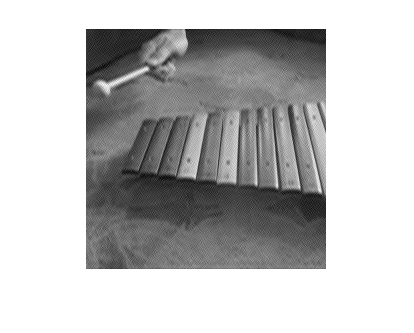} 
\endminipage\hfill \hspace{-2cm}  \vspace{-.9cm}
\caption{\small Restored images by the {\tt nested}$\_$G-tAT$_p$ (left), G-tGMRES$_p$
(middle), and tGMRES$_p$ (right) for $\widetilde{\delta}=10^{-3}$.}
\label{Fig: 6}
\end{figure}

The relative errors, PSNR and CPU times are displayed in Table \ref{Tab: 5}. The tAT$_p$ and {\tt nested}$\_$tAT$_p$ methods, which do not involve flattening, are seen to yield 
restorations of the highest quality for all noise levels. The tGMRES method is
the fastest for $\widetilde{\delta}=10^{-2}$ and $10^{-3}$, but gives the worst restorations for all noise levels. Solution methods 
that involve flattening, such as tAT$_p$, GG-tAT and G-tGMRES$_p$ and GG-tGMRES methods, are the slowest for $\widetilde{\delta}=10^{-3}$.

\begin{table}[h!]
\begin{center}
\begin{tabular}{ccccccc} 
\hline
Noise level &Method &$\ell$ & $\mu_\ell$& PSNR & Relative error& CPU time (secs) \\
\hline
\multirow{7}{2em}{$10^{-3}$} & tAT$_p$ &-&- &34.07&4.15e-02 &52.53 \\ 
&{\tt nested}$\_$tAT$_p$ &10 &-&33.81& 4.27e-02 & 46.40  \\ 
&G-tAT$_p$&-&- &33.27 & 4.54e-02 & 463.08 \\ 
&GG-tAT&22& 7.46e+02 &33.24& 4.56e-02 & 203.42  \\ 
&tGMRES$_p$&-&- &27.21& 9.13e-02 & 31.88  \\ 
&G-tGMRES$_p$&-&- &33.17& 4.60e-02 & 406.93 \\ 
&GG-tGMRES &22&- &33.21& 4.58e-02 & 97.66 \\  \hline
\multirow{7}{2em}{$10^{-2}$} & tAT$_p$ &-&-&30.75&6.07e-02 &20.28\\ 
&{\tt nested}$\_$tAT$_p$ &3 &-&25.64& 1.09e-01 & 18.01  \\ 
&G-tAT$_p$&-&-&27.22& 9.12e-02 & 69.21\\
&GG-tAT&8&2.16e+02 &27.22& 9.12e-02 & 25.74  \\
&tGMRES$_p$&-&-&15.67& 3.45e-01 & 8.23  \\  
&G-tGMRES$_p$&-&-&26.82& 9.55e-02 & 52.78 \\ 
&GG-tGMRES &8&-&26.82& 9.55e-02 &11.93 \\  \hline
\multirow{7}{2em}{$10^{-1}$} & tAT$_p$ &-&-&24.69&1.22e-01 &13.84\\ 
&{\tt nested}$\_$tAT$_p$ &2 &-&21.25& 1.81e-01 & 16.94  \\ 
&G-tAT$_p$&-&- &21.17& 1.83e-01 & 5.24\\ 
&GG-tAT&2&1.04e+01 &21.17& 1.83e-01 & 1.60  \\ 
&tGMRES$_p$&-&- &0.45& 1.99e+00 & 3.46  \\ 
&G-tGMRES$_p$&-&-&19.21& 2.29e-01 & 3.16 \\
&GG-tGMRES &2&-&19.21& 2.29e-01 &0.68 \\  \hline
\end{tabular}
\end{center}\vspace{-.5cm}
\caption{\small Results for Example \ref{E5}.}
\label{Tab: 5}
\end{table}

\end{Ex}

\section{Conclusion}\label{sec7}
This paper extends the standard Arnoldi iteration for matrices to third order tensors and 
describes several algorithms based on this extension for solving linear discrete ill-posed 
problems with a t-product structure. The solution methods are based on computing a few 
steps of the extended Arnoldi process, which is referred to as the t-Arnoldi process. The 
global t-Arnoldi and generalized global t-Arnoldi processes also are considered. 
Differently from the t-Arnoldi process, the latter processes involve flattening. Both 
Tikhonov regularization and regularization by truncated iteration are illustrated. The 
latter gives rise to an extension of the standard GMRES method, referred to as the tGMRES 
and global tGMRES methods. The discrepancy principle is used to determine the number of 
iterations with the t-Arnoldi, global t-Arnoldi, and generalized global t-Arnoldi 
processes, as well as the regularization parameter in Tikhonov regularization and the 
number of iterations by the Arnoldi-type and GMRES-type methods. The effectiveness of the
proposed methods is illustrated by applications to image and video restorations. Solution 
methods such as tAT, tAT$_p$, and {\tt nested$\_$}tAT$_p$ that avoid matricization or 
vectorization of discrete ill-posed problems for tensors show great promise in terms of 
speed and quality of the computed restorations determined by their relative errors and PSNR when compared to solution methods that matricize or vectorize.

\section*{Acknowledgment}
The authors would like to thank the referees for comments that led to improvements of the
presentation. Research by LR was supported in part by NSF grant DMS-1720259.


\begin{thebibliography}{35} \label{sec:bib}

\bibitem{BJNR}
F. P. A. Beik, K. Jbilou, M. Najafi-Kalyani, and L. Reichel, Golub-Kahan bidiagonalization
for ill-conditioned tensor equations with applications, Numer. Algorithms, 84 (2020), pp.
1535--1563.

\bibitem{BNR} 
F. P. A. Beik, M. Najafi-Kalyani, and L. Reichel, Iterative Tikhonov regularization of 
tensor equations based on the Arnoldi process and some of its generalizations, Appl. 
Numer. Math., 151 (2020), pp. 425--447.

%\bibitem{B} 
%K. Braman, Third-order tensors as linear operators on a space of matrices, Linear Algebra Appl. 433 (2010) 1241–1253.

\bibitem{BPR}
A. Buccini, M. Pasha, and L. Reichel, Generalized singular value decomposition with 
iterated Tikhonov regularization, J. Comput. Appl. Math., 373 (2020), Art. 112276.

\bibitem{CLR}
D. Calvetti, B. Lewis, and L. Reichel, On the regularizing properties of the GMRES method,
Numer. Math., 91 (2002), pp. 605--625.

\bibitem{CMRS} 
D. Calvetti, S. Morigi, L. Reichel, and F. Sgallari, Tikhonov regularization and the 
L-curve for large, discrete ill-posed problems, J. Comput. Appl. Math., 123 (2000), 
pp. 423--446.

\bibitem{CR}
D. Calvetti and L. Reichel, Tikhonov regularization of large linear problems, BIT Numer. Math., 
43 (2003), pp. 263--283.

\bibitem{DMR}
M. Donatelli, D. Martin, and L. Reichel, Arnoldi methods for image deblurring with 
anti-reflective boundary conditions, Appl. Math. Comput., 253 (2015), pp. 135--150.

\bibitem{GIJB} 
M. El Guide, A. El Ichi, K. Jbilou, and F. P. A Beik, Tensor GMRES and Golub-Kahan
bidiagonalization methods via the Einstein product with applications to image and video
processing, https://arxiv.org/pdf/2005.07458.pdf

\bibitem{IGJ}
A. El Ichi, M. El Guide and K. Jbilou, Discrete cosine transform LSQR and GMRES methods for multidimensional ill-posed problems, March 2021. https://arxiv.org/pdf/2103.11847.pdf

\bibitem{GIJS} 
M. El Guide, A. El Ichi, K. Jbilou, and R. Sadaka, Tensor Krylov subspace methods via the
T-product for color image processing, June 2020. https://arxiv.org/pdf/2006.07133.pdf

\bibitem{EANK} 
G. Ely, S. Aeron, N. Hao, and M. E. Kilmer, 5d and 4d pre-stack seismic data completion 
using tensor nuclear norm (TNN), SEG International Exposition and Eighty-Third Annual 
Meeting at Houston, TX, 2013.

\bibitem{EHN} 
H. W. Engl, M. Hanke, and A. Neubauer, Regularization of Inverse Problems, Kluwer, 
Dordrecht, 1996.

\bibitem{FRR}
C. Fenu, L. Reichel, and G. Rodriguez, GCV for Tikhonov regularization via global 
Golub-Kahan decomposition, Numer. Linear Algebra Appl., 23 (2016), pp. 467--484.

\bibitem{GNR}
S. Gazzola, P. Novati, and M. R. Russo, On Krylov projection methods and Tikhonov 
regularization, Electron. Trans. Numer. Anal., 44 (2015), pp. 83--123.

%\bibitem{GGV} 
%D. F. Gleich, C. Greif, and J. M. Varah, The power and Arnoldi methods in the algebra of
%circulants, Numer. Linear Algebra Appl., 2012, dio 10.1002/nla.1845.

\bibitem{GHW} 
G. H. Golub, M. Heath, and G. Wahba, Generalized cross-validation as a method for choosing 
a good ridge parameter, Technometrics, 21 (1979), pp. 215--223.

\bibitem{H2n}  
P. C. Hansen, Analysis of discrete ill-posed problems by means of the L-curve, SIAM Rev.,
34 (1992), pp. 561--580.

\bibitem{Hn} 
P. C. Hansen, Rank-Deficient and Discrete Ill-Posed Problems, SIAM, Philadelphia, 1998.

\bibitem{Haa} 
P. C. Hansen, Regularization tools version 4.0 for MATLAB 7.3. Numer. Algorithms, 46 
(2007), pp. 189--194.

\bibitem{HKBH} 
N. Hao, M. E. Kilmer, K. Braman, and R. C. Hoover, Facial recognition using tensor-tensor 
decompositions, SIAM J. Imaging Sci., 6 (2013), pp. 437--463.

\bibitem{HRY}
G. Huang, L. Reichel, and F. Yin, On the choice of subspace for large-scale Tikhonov 
regularization problems in general form, Numer. Algorithms, 81 (2019), pp. 33--55.

%\bibitem{LBC} 
%Y. Lin, L. Bao and Y. Cao, Augmented Arnoldi-Tikhonov regularization methods for solving large-scale linear ill-posed systems, Mathematical Problems in Engineering, 2013, pp. 1–8.

\bibitem{KKA} 
E. Kernfeld, M. Kilmer, and S. Aeron, Tensor-tensor products with invertible linear 
transforms, Linear Algebra Appl., 485 (2015), pp. 545--570.

\bibitem{KBH} 
M. Kilmer, K. Braman, and N. Hao, Third order tensors as operators on matrices: A 
theoretical and computational framework, Tufts University, Department of Computer 
Science, Tech. Rep., January 2011.

\bibitem{KBHH} 
M. E. Kilmer, K. Braman, N. Hao, and R.  C. Hoover, Third-order tensors as operators on 
matrices: A theoretical and computational framework with applications in imaging, SIAM 
J. Matrix Anal. Appl., 34 (2013), pp. 148--172.

\bibitem{KM} 
M. E. Kilmer and C. D. Martin, Factorization strategies for third-order tensors, Linear 
Algebra Appl., 435 (2011), pp. 641--658.

\bibitem{Ki}
S. Kindermann, Convergence analysis of minimization-based noise level-free parameter 
choice rules forlinear ill-posed problems, Electron. Trans. Numer. Anal., 38 (2011), pp. 
233--257.

\bibitem{KR}
S. Kindermann and K. Raik, A simplified L-curve method as error estimator. Electron. 
Trans. Numer. Anal., 53 (2020), pp. 217--238.

\bibitem{KB} 
T. G. Kolda and B. W. Bader, Tensor decompositions and applications, SIAM Rev., 51 (2009), 
pp. 455--500.

\bibitem{LR} 
B. Lewis and L. Reichel, Arnoldi-Tikhonov regularization methods, J. Comput. Appl. Math.,
226 (2009), pp. 92--102.

\bibitem{FCLLY}
C. Lu, J. Feng, Y. Chen, W. Liu, Z. Lin, and S. Yan, Tensor robust principal component 
analysis with a new tensor nuclear norm, IEEE Trans. Pattern Anal. Mach. Intell., 42 
(2020), pp. 925--938, doi:10.1109/TPAMI.2019.2891760.

\bibitem{Lund} 
K. Lund. The tensor t-function: a definition for functions of third-order tensors. ArXiv preprint, arXiv:1806.07261, 2018.

\bibitem{MSL}  
C. D. Martin, R. Shafer, and B. LaRue, An order-$p$ tensor factorization with applications
in imaging. SIAM J. Sci. Comput., 35 (2013), pp. A474--A490.

\bibitem{Ne}
A. Neubauer, Augmented GMRES-type versus CGNE methods for the solution of linear ill-posed
problems, Electron. Trans. Numer. Anal., 51 (2019), pp. 412--431.

\bibitem{MQW1}
Y. Miao, L. Qi and Y. Wei, T-Jordan Canonical Form and T-Drazin Inverse Based on the T-Product. Commun. Appl. Math. Comput. (2020). https://doi.org/10.1007/s42967-019-00055-4

\bibitem{MQW2}
Y. Miao, L. Qi and Y. Wei, Generalized tensor function via the tensor singular value decomposition based on the T-product. Linear Algebra and its Applications, 590 (2020) 258-303.

%\bibitem{NKH} 
%E. Newman, M. Kilmer, and L. Horesh, Image classification using local tensor singular
%value decompositions, 2017 IEEE 7th International Workshop on Computational Advances
%in Multi-Sensor Adaptive Processing (CAMSAP), Curacao, 2017, pp. 1-5, doi: 10.1109/CAM-
%SAP.2017.8313137.

\bibitem{RS}
L. Reichel and A. Shyshkov, A new zero-finder for Tikhonov regularization, BIT Numer. 
Math., 48 (2008), pp. 627--643.

\bibitem{LU}  
L. Reichel and U. O. Ugwu, The tensor Golub-Kahan-Tikhonov method applied to the solution 
of ill-posed problem with a t-product structure, 2020, submitted for publication. 

\bibitem{RR} 
L. Reichel and G. Rodriguez, Old and new parameter choice rules for discrete ill-posed 
problems, Numer. Algorithms, 63 (2013), pp. 65--87.

\bibitem{Sa}
Y. Saad, Iterative Methods for Sparse Linear Systems, 2nd ed., SIAM, Philadelphia, 2003.

\bibitem{SS} 
Y. Saad and M. H. Schultz, GMRES: a generalized minimal residual method for solving 
nonsymmetric linear systems, SIAM J. Sci. Stat. Comput., 7 (1986), pp. 856--869.

%\bibitem{SE1} 
%B. Savas and L. Eldén, Krylov subspace methods for tensor computations, Technical Report LiTH-MAT-R-2009-02-SE, Linköpings University, Februrary 2009. (Available from: http://www.mai.liu.se/ laeld/tensorKrylov.pdf).

\bibitem{SE2} 
B. Savas and L. Eld\'en, Krylov-type methods for tensor computations, Linear Algebra 
Appl., 438 (2013), pp. 891--918.

\bibitem{SKH} 
S. Soltani, M. E. Kilmer, and P. C. Hansen, A tensor-based dictionary learning approach to tomographic image reconstruction, BIT Numer. Math., 56 (2015), pp. 1425--1454.

\bibitem{TB} 
L. N. Trefethen and D. Bau III, Numerical Linear Algebra, SIAM, Philadelphia, 1997. 

\bibitem{ZEAHK} 
Z. Zhang, G. Ely, S. Aeron, N. Hao, and M. E. Kilmer, Novel methods for multilinear data 
completion and de-noising based on tensor-svd, In 2014 IEEE Conference on Computer Vision 
and Pattern Recognition, CVPR 2014, Columbus, OH, USA, June 23-28, 2014, pp. 3842--3849.
\end{thebibliography}
\end{document}